\documentclass[10pt]{amsart}
\theoremstyle{plain}
\usepackage{amsmath, amsfonts, amsthm, amssymb,amscd,bbold}
\usepackage{graphics}
\usepackage{color}
\usepackage{epsfig}
\usepackage{hyperref}
\usepackage[all]{xy}
\graphicspath{ {Figures/} }

\title[Quiver algebras and bordered Floer homology]{On Khovanov-Seidel quiver algebras and bordered Floer homology}
\author{Denis Auroux}
\thanks{DA was partially supported by NSF grant DMS-1007177.}
\address{UC Berkeley, Department of Mathematics, 970 Evans Hall \# 3840, Berkeley CA 94720, USA}
\email{auroux@math.berkeley.edu}

\author{J. Elisenda Grigsby}
\thanks{JEG was partially supported by a Viterbi-endowed MSRI postdoctoral fellowship and NSF grant number DMS-0905848.}
\address{Boston College; Mathematics Department; 301 Carney Hall; Chestnut Hill, MA 02467, USA}
\email{grigsbyj@bc.edu}

\author{Stephan M. Wehrli}
\thanks{SMW was partially supported by an MSRI postdoctoral fellowship.}
\address{Syracuse University; Mathematics Department; 215 Carnegie; Syracuse, NY 13244, USA}
\email{smwehrli@syr.edu}

\theoremstyle{plain}
\newtheorem{theorem}{Theorem}[section]
\newtheorem{construction}[theorem]{Construction}
\newtheorem{lemma}[theorem]{Lemma}
\newtheorem{proposition}[theorem]{Proposition}

\newtheorem*{theoremA}{Theorem~\ref{thm:filtdiffmodgen}}
\theoremstyle{definition}
\newtheorem{notation}[theorem]{Notation}
\newtheorem{definition}[theorem]{Definition}
\newtheorem{remark}[theorem]{Remark}

\newcommand{\Szabo}{{Szab{\'o}} }

\newcommand{\betatwoleft}{\ensuremath{\left(\beta_k^{Kh}\right)_{(1|1|0)}}}
\newcommand{\betatworight}{\ensuremath{\left(\beta_k^{Kh}\right)_{(0|1|1)}}}

\newcommand{\Z}{\ensuremath{\mathbb{Z}}}
\newcommand{\F}{\ensuremath{\mathbb{F}}}

\newcommand{\C}{\ensuremath{\mathbb{C}}}

\newcommand{\boldSigma}{\ensuremath{\mbox{\boldmath $\Sigma$}}}
\newcommand{\cG}{\ensuremath{\mathcal{G}}}
\newcommand{\cB}{\ensuremath{\mathcal{B}}}

\newcommand{\cI}{\ensuremath{\mathcal{I}}}
\newcommand{\cJ}{\ensuremath{\mathcal{J}}}

\newcommand{\cF}{\ensuremath{\mathcal{F}}}

\newcommand{\cM}{\ensuremath{\mathcal{M}}}
\newcommand{\Ztwo}{\ensuremath{\mathbb{Z}/2\mathbb{Z}}}

\newcommand{\tildeotimes}{\ensuremath{\widetilde{\otimes}}}

\newcommand{\Hom}{{\rm Hom}}
\newcommand{\Id}{\ensuremath{\mbox{\textbb{1}}}}


\begin{document}
\bibliographystyle{plain}

\begin{abstract} We discuss a relationship between Khovanov- and Heegaard Floer-type homology theories for braids.  Explicitly, we define a filtration on the bordered Heegaard-Floer homology bimodule associated to the double-branched cover of a braid and show that its associated graded bimodule is equivalent to a similar bimodule defined by Khovanov and Seidel.\end{abstract}
\maketitle

\section{Introduction}\label{sec:Intro}
The low-dimensional topology community has been energized in recent years by the introduction of a wealth of so-called {\em homology-type} invariants.  These invariants are defined by associating to a topological object (for example, a link or a $3$--manifold) an abstract chain complex whose quasi-isomorphism class--hence, homology--is an invariant of the object.

One obtains such invariants from two apparently unrelated points of view: 
\begin{enumerate}
	\item algebraically, via the higher representation theory of quantum groups, and 
	\item geometrically/analytically, via symplectic geometry and gauge theory.
\end{enumerate}

Although the invariants themselves share a number of formal properties, finding explicit connections between the two viewpoints has proven challenging.

A striking success in this direction is a result of Ozsv{\'a}th and \Szabo relating the $\Ztwo$ versions of Khovanov homology and Heegaard Floer homology:

\begin{theorem} \cite{MR2141852} Let $L \subset S^3$ be a link and $\overline{L} \subset S^3$ denote its mirror.  There exists a spectral sequence whose $E^2$ term is $\widetilde{Kh}(\overline{L})$, the reduced Khovanov homology of the mirror of $L$, and whose $E^\infty$ term is $\widehat{HF}(\boldSigma(L))$, the Heegaard-Floer homology of the double-branched cover of $L$.
\end{theorem}

This result has generated applications in a number of directions (see, e.g., \cite{GT0412184}, \cite{GT08071341}, \cite{GT08082336}).  It also served as inspiration for Kronheimer and Mrowka's  construction of an analogous spectral sequence from Khovanov homology to a version of instanton knot homology, yielding a proof that Khovanov homology detects the unknot \cite{KhDetectUnknot}.

The aim of the present paper is to move toward a more ``atomic" understanding of the Ozsv{\'a}th-\Szabo spectral sequence and its sutured generalizations (\cite{GT07060741, GT08071432, AnnularLinks, SurfDecomp}).  In particular, viewing a link in $S^3$ as the closure of a braid, we can ask whether there are appropriate Khovanov-type (algebraic) and Heegaard-Floer-type (geometric/analytic) invariants associated to braids such that the Ozsv{\'a}th-\Szabo spectral sequence emerges as an algebraic consequence of a relationship between these invariants.  

Such a description would not only be of theoretical interest.  Ozsv{\'a}th-Szab{\'o}'s original description of the above spectral sequence involves holomorphic polygon counts in Heegaard multi-diagrams.  Since these counts are tricky to carry out in practice, finding ways to perform them combinatorially should prove valuable, especially in light of subsequent work of Baldwin \cite{MR2822178}  (see also L. Roberts \cite{ GT08082817}) proving that the terms of the Ozsv{\'a}th-\Szabo spectral sequence are themselves link invariants.  

We should at this point remark that recent work of Lipshitz-Ozsv{\'a}th-Thurston, in \cite{GT10110499} and its sequel, does precisely this.  In addition, Szab{\'o} \cite{GT10104252} has constructed a combinatorial filtration on the Khovanov cube of resolutions associated to a link diagram that he conjectures yields the original Ozsv{\'a}th-\Szabo spectral sequence.

In the present paper, we address a slightly different question from a substantially different direction.  First, we focus not on the original Ozsv{\'a}th-\Szabo spectral sequence but rather on (a direct summand of) one of its sutured generalizations \cite{GT07060741, AnnularLinks}.  Second, we take as our starting point a paper of Khovanov-Seidel \cite{MR1862802}, which explores a concrete instance of Kontsevich's homological mirror symmetry conjecture \cite{MR1403918}.  The constructions found there, when combined with work of the first author \cite{GT10014323}, lead naturally to a new view on the filtered complexes appearing in \cite{GT07060741, AnnularLinks}.

Explicitly, given a braid $\sigma \subset D^2 \times I$, we consider the closure of the braid, not in the three-ball but in the solid torus (viewed as a product sutured annulus, $A \times I$).  Associated to the resulting {\em annular link}  are Khovanov-type and Heegaard-Floer-type invariants connected by a {\em sutured} spectral sequence \cite{MR2113902, GT07060741, AnnularLinks} that splits along an extra grading measuring ``wrapping" around the $S^1$ factor.\footnote{This extra grading has a natural interpretation on the Khovanov side in terms of $U_q(\mathfrak{sl}_2)$ weight space decompositions and on the Heegaard-Floer side in terms of relative Spin$^c$ structures.  See \cite{JacoFest} for more details.}  In \cite{HochHom}, building on work in \cite{GT10030598}, we obtain a similar spectral sequence in the ``next-to-top" graded piece as the Hochschild homology of a filtered $A_\infty$ bimodule associated to the original braid, $\sigma$. 

The purpose of the present paper is to give an explicit combinatorial construction of this filtered $A_\infty$ bimodule.  
Informally, the resulting spectral sequence interpolates between the ``open" Khovanov- and Heegaard-Floer-type invariants of a braid $\sigma \subset D^2 \times I$ just as the sutured spectral sequence interpolates between the analogous ``closed" invariants of its closure, $\hat{\sigma} \subset A \times I$.

More precisely:
\begin{enumerate}
	\item On the algebraic side, we show how to use ideas of Khovanov-Seidel in \cite{MR1862802} to construct an $A_\infty$ bimodule, $\mathcal{M}_\sigma^{Kh}$, via Yoneda imbedding of a distinguished collection of objects in the derived category of a quiver algebra.
	\item On the geometric/analytic side, we use the bordered Floer homology package of Lipshitz-Ozsv{\'a}th-Thurston in \cite{GT08100687, GT10030598} to construct an $A_\infty$ bimodule, $\mathcal{M}_\sigma^{HF}$, the $1$--strand CFDA bimodule associated to the mapping class $\hat{\sigma}$ obtained as the double-branched cover of $\sigma \subset D^2 \times I$.
\end{enumerate}

Letting $\Id$ denote the identity braid of the same index as $\sigma$, we prove:

\begin{theoremA} There exists a filtration on $\mathcal{M}^{HF}_\sigma$ whose associated graded bimodule is quasi-isomorphic, as an ungraded $A_\infty$ bimodule over $\left\{\mbox{gr}(\cM_{\Id}^{HF})=\cM_{\Id}^{Kh}\right\}$, to $\mathcal{M}_\sigma^{Kh}$.
\end{theoremA} 

In particular, for each braid there exists a spectral sequence connecting the Khovanov-Seidel (algebraic) bimodule to the Lipshitz-Ozsv{\'a}th-Thurston (geometric/analytic) one.  Moreover, these ``open" spectral sequences can be defined without reference to holomorphic curves.  In fact, our construction is based on a remarkably simple toy model (Lemma \ref{lem:ZtwoequivS1}): a filtered complex interpolating between the cohomology of $S^1$ and the cohomology of $S^0$ (both over \Ztwo) coming from a \Ztwo--equivariant cochain complex for $S^1$.  This toy model was, in turn, inspired by work of Seidel and Smith \cite{SeidelSmithInvolution}.

We pause here to emphasize some key points. First, the algebraic objects appearing in \cite{MR1862802} do themselves admit a geometric interpretation in terms of the Fukaya category of a certain Lefschetz fibration (cf. Section \ref{sec:BKhFukaya}). Those readers familiar with \cite{SeidelBook} may therefore prefer to perform the Section \ref{sec:Khalgmod} calculations geometrically. We have opted instead to work entirely in the algebraic setting, using symplectic geometry only as motivation. Although this has surely increased the paper's length, we hope it has simultaneously increased its accessibility to non-geometers.

This accessibility is essential, as the algebraic version of the Khovanov-Seidel construction has a beautiful representation-theoretic interpretation. Explicitly, the Khovanov-Seidel algebra is a special case (for $k=1$) of a family of algebras $A^{k,n-k}$, introduced by Chen-Khovanov~\cite{QA0610054} and independently by Stroppel~\cite{StroppelAlgebra}, giving rise to a categorification of the $U_q(\mathfrak{sl}_2)$ Reshetikhin-Turaev invariant for tangles. These algebras can also be identified with endomorphism algebras of projective generators of certain blocks $\mathcal{O}_{k,n-k}$ of category $\mathcal{O}$.  We conjecture that Theorem~\ref{thm:filtdiffmodgen} admits a generalization which, for every $n$-strand braid $\sigma$, provides a relationship between the $k$-strand part of the Lipshitz-Ozsv\'ath-Thurston bimodule associated to $\widehat{\sigma}$ and a Khovanov-type bimodule defined over the Ext-algebra of the direct sum of all standard $A^{k,n-k}$-modules.

We end by remarking that the construction of our filtration required a choice of a common ``basis" of generators for the relevant Fukaya categories. One natural choice is made in \cite{MR1862802} (corresponding in the geometric setting to Lagrangians where all but one is compact and in the algebraic setting to Luzstig's canonical basis for a tensor product representation), while another equally natural choice is made in  \cite{GT08100687}, as reinterpreted in \cite{GT10014323} (corresponding to non-compact Lagrangians and the standard basis for a tensor product representation). We work with the latter, {\em noncompact} basis because both  ($k=1$) algebras in the noncompact case are formal (see Lemma \ref{lem:Bformal} and \cite[Prop. 3.6]{GT10014323}) while the bordered Floer algebra corresponding to the compact basis is not \cite[Chp. 20]{SeidelBook}, \cite{LekiliPerutz}. 

The paper is organized as follows:

In Section \ref{sec:Algprelim}, we establish notation and collect a number of useful definitions and elementary algebraic results.

In Section \ref{sec:Khalgmod}, we describe the topological input needed for the algebraic constructions in the remainder of the paper.  After reviewing the key points in \cite{MR1862802}, we proceed to the construction and description of 
\begin{itemize}
	\item an algebra, $B^{Kh}$, associated to a marked disk $D_m$ equipped with a specific {\em basis of curves} and 
	\item a module, $\mathcal{M}_\sigma^{Kh}$, associated to each braid $\sigma$, decomposed as a product of elementary Artin generators.
\end{itemize}
We conclude the section with a brief discussion of the Fukaya-theoretic interpretation of $B^{Kh}$ and $M_\sigma^{Kh}$.

In Section \ref{sec:HFalgmod}, we turn to the construction and description of the analogous bordered Floer algebra $B^{HF}$ and bimodules $\mathcal{M}_\sigma^{HF}$, using the same topological input.

In Section \ref{sec:Filtalg}, we describe a natural filtration on $B^{HF}$ whose associated graded algebra is isomorphic to $B^{Kh}$.  Our construction is based on a simple ``toy model" (Lemma \ref{lem:ZtwoequivS1}).

In Section \ref{sec:Filtdiffmod}, we describe a filtration on $\mathcal{M}_\sigma^{HF}$ whose associated graded homology bimodule is quasi-isomorphic to $\mathcal{M}_\sigma^{Kh}$.  We proceed by choosing a decomposition \[\sigma = \sigma_{k_1}^\pm \cdots \sigma_{k_n}^\pm\] of $\sigma$ as a product of elementary Artin generators, explicitly constructing a filtration on $\mathcal{M}_{\sigma_k^\pm}^{HF}$ for each elementary generator, then realizing $\mathcal{M}_\sigma^{HF}$ as the (filtered) $A_\infty$ tensor product of the elementary bimodules $\mathcal{M}_{\sigma_{k_1}^\pm}^{HF}, \ldots, \mathcal{M}_{\sigma_{k_n}^\pm}^{HF}$.

In Section \ref{sec:Example}, we describe an example highlighting the nontriviality of the filtration on $\mathcal{M}_\sigma^{HF}$.

\subsection{Acknowledgements} We are grateful to Tony Licata, Robert Lipshitz, Peter Ozsv{\'a}th, Catharina Stroppel, and Dylan Thurston for a great number of interesting conversations, and to the MSRI semester-long program on Homology Theories of Knots and Links for making these conversations possible.  We would also like to thank Joshua Sussan for bringing to our attention that some of the algebraic results of Section 3 (in particular, Lemma \ref{lem:Bformal}) were independently obtained by Angela Klamt and Catharina Stroppel in \cite{Klamt} and \cite{KlamtStroppel}. Many thanks are also due to the excellent referee and editor, whose insightful suggestions greatly improved the manuscript. Finally, we are indebted to John Baldwin, who helped us find the example described in Section \ref{sec:Example}.

\section{Algebraic preliminaries} \label{sec:Algprelim}

In this section, we establish some basic facts about filtered $A_\infty$ algebras and modules.
We assume throughout that we are working over the field $\F = \Z/2\Z$.  In addition, many of the spaces we discuss will be graded either by $\Z$, in which case we say it is {\em graded}, or by $\Z^2$, in which case we say it is {\em bigraded}. The {\em (co)homological} grading always appears first.

\begin{notation} If $V$ is a bigraded vector space, i.e. \[V = \bigoplus_{i,j \in \Z} V_{(i,j)},\] and $k_1, k_2 \in \Z$, then $V[k_1]\{k_2\}$ will denote the vector space whose first (homological) grading has been shifted {\em down} by $k_1$ and whose second (internal) grading has been shifted {\em up} by $k_2$.\footnote{The difference in shift conventions for homological versus internal gradings is unfortunate, but standard in the literature. In particular, they coincide with those in \cite{MR1862802}. The reader should be warned, however, that \cite{GT08100687} uses a different convention, since their differential maps {\em decrease} rather than {\em increase} homological grading.}  Explicitly, \[\left(V[k_1]\{k_2\}\right)_{(i,j)} \cong V_{(i+k_1,j-k_2)}.\]
\end{notation}

We omit the standard definitions of $A_\infty$ algebras, modules, morphisms and homotopies, instead referring to Keller's expository papers: \cite{MR1854636},\cite{MR2258042}. Other excellent references are Seidel's book: \cite{SeidelBook} (though the reader should be warned that Seidel's ordering conventions (cf. Eqn 1.1) for $A_\infty$ morphisms differ from ours), the thesis of Lef{\`e}vre-Kasegawa \cite{kenji}, and Chapter 2 of \cite{GT08100687}. All $A_\infty$ modules we consider will be over homologically unital algebras (c-unital, in the terminology of \cite{SeidelBook}), and morphisms between homologically unital algebras must be homologically unital.

\begin{remark}
The algebraically defined modules we study here are, in fact, {\em strictly} unital (cf. Remark \ref{rmk:Unitality}),  but the geometrically defined ones need not be. 
\end{remark}

Let $n \in \Z^+$ and $n_1, n_2 \in \Z^{\geq 0}$. We shall use the notation $m_n$ to refer to the $n$th structure map \[m_n: {\bf A}^{\otimes n} \rightarrow {\bf A}[2-n],\] of an $A_\infty$ algebra ${\bf A}$ and the notation $m_{(n_1|1|n_2)}$ to refer to the $(n_1|1|n_2)$ structure map \[m_{(n_1|1|n_2)}: {\bf A}^{\otimes n_1} \otimes {\bf M} \otimes {\bf B}^{\otimes n_2} \rightarrow {\bf M}[2-(n_1+1+n_2)]\] of a bimodule ${\bf M}$ admitting a left (resp., right) $A_\infty$ action by the $A_\infty$ algebra ${\bf A}$ (resp., ${\bf B}$).

If ${\bf A}$ is ungraded but otherwise satisfies all of the conditions of an $A_\infty$ algebra, we call ${\bf A}$ an {\em ungraded $A_\infty$ algebra}.

A graded (resp., ungraded) $A_\infty$ algebra satisfying $m_n = 0$ for all $n > 2$ is a {\em differential graded algebra (dga)} (resp., a {\em differential algebra}) with differential $\partial := m_1$ and multiplication $m_2$. The terminology is completely analogous for graded and ungraded $A_\infty$ and differential modules.

If ${\bf M}$ and ${\bf N}$ are $A_\infty$ bimodules, we will refer to the map \[f_{(n_1|1|n_2)}: {\bf A}^{\otimes n_1} \otimes {\bf M} \otimes {\bf B}^{\otimes n_1} \rightarrow {\bf N}[1-(n_1+1+n_2)]\] associated to an $A_\infty$ morphism $f$ as the ``{\em $(n_1|1|n_2)$ term of $f$}."  In addition, we will use the terminology ``{\em $(n_1|1|n_2)$ $A_\infty$ relation}" to refer to the $A_\infty$ relation corresponding to $n_1$ left inputs and $n_2$ right inputs.  For example, the $(1|1|0)$ $A_\infty$ relation for an $A_\infty$ morphism $f: {\bf M} \rightarrow {\bf N}$ is given by:
\[f_{(1|1|0)}(m_1 \otimes \Id + \Id \otimes m_{(0|1|0)}) + f_{(0|1|0)} m_{(1|1|0)} = m_{(1|1|0)}(\Id \otimes f_{(0|1|0)}) + m_{(0|1|0)}f_{(1|1|0)}.\]


If $f_{(n_1|1|n_2)} = 0$ for all $n_1,n_2 >0$, then we say that $f = f_{(0|1|0)}: {\bf A} \rightarrow {\bf B}$ is a {\em strict morphism} of $A_\infty$ modules.  In particular, a strict morphism $f: {\bf A} \rightarrow {\bf B}$ of differential (graded) algebras is a chain map intertwining the multiplication, $m_2$.

An $A_\infty$ morphism $f$ is said to be a {\em quasi-isomorphism} if $f_1$ induces an isomorphism on homology.

Homological perturbation theory allows one to {\em transfer} $A_\infty$ structures along certain morphisms. Although the situation of particular interest to us is the transfer of an $A_\infty$ structure along a chain homotopy equivalence $p: {\bf A} \rightarrow H_*({\bf A})$ as in \cite{MR720689, Merkulov, KontsevichSoibelman}, such a transfer can be performed in much greater generality. See \cite{MR2287133} (and the related discussion in \cite[Sec. (1i)]{SeidelBook}). A nice account is also given in \cite[Thm. 2.1]{Berglund}. The {tree formulas} for this transferred structure are summarized in the following proposition:

	


\begin{proposition} \label{prop:InducedAinfty} Let 
\[\xymatrix{{\bf A} \ar@(dl,ul)[]^h \ar@/^/[r]^-p & H_*({\bf A}) \ar@/^/[l]^-\iota}\] be a contraction of a chain complex, ${\bf A}$, onto its homology, $H_*({\bf A})$. In other words, $p$ and $\iota$ are chain maps and $h$ is a chain homotopy satisfying: \begin{equation} \label{eqn:Ainftymap}
\iota p = \mbox{Id} + \partial h + h\partial, \hskip 20pt p \iota = \mbox{Id}. 
\end{equation}

Suppose ${\bf A}$ is further endowed with a (not necessarily unital) $A_\infty$ structure extending the differential structure, i.e., multiplication maps \[m_n^A: {\bf A}^{\otimes n} \rightarrow {\bf A}[2-n]\] with $ m_1^A = \partial$ and satisfying the $A_\infty$ relations.

Then $H_*({\bf A})$ admits a (not necessarily unital) $A_\infty$ algebra structure such that
\begin{enumerate}
	\item $m_1 =0$ and $m_2$ is induced from $m_2^A$,
	\item there are $A_\infty$ quasi-isomorphisms $p': {\bf A} \rightarrow H_*({\bf A})$ and $\iota': H_*({\bf A}) \rightarrow {\bf A}$, and an $A_\infty$ homotopy $h': {\bf A} \rightarrow {\bf A}$ extending $p, \iota, h$.
\end{enumerate}

The $n$th $A_\infty$ multiplication \[m_n: (H_*({\bf A}))^{\otimes n} \rightarrow H_*({\bf A})[2-n]\] is given by \[m_n := \sum_T m_n^T\] where the sum ranges over all planar rooted trees $T$ with $n$ leaves and $m_n^T$ is defined by applying the $T$--shaped diagram with 
\begin{enumerate}
	\item leaves labeled with $\iota$,
	\item interior edges labeled with $h$,
	\item vertices labeled with the multiplication maps $m_i$ in the algebra ${\bf A}$, and
	\item root labeled with $p$
\end{enumerate}
to an element of $\left(H_*({\bf A})\right)^{\otimes n}$.
\end{proposition}

\begin{figure}
\begin{center}
\resizebox{3in}{!}{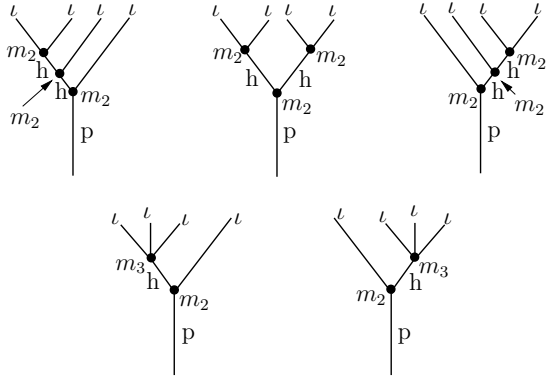}
\end{center}
\caption{The full collection of rooted trees with $4$ inputs specifying the multiplication $m_4$ described by Proposition \ref{prop:InducedAinfty}.}
\label{fig:RootedTrees}
\end{figure}

 See Figure \ref{fig:RootedTrees} for an enumeration of all such rooted trees $T$ specifying the multiplication $m_n$ when $n=4$.  The resulting ``transferred" $A_\infty$ structure on $H_*({\bf A})$ is unique (independent of the choice of $p, \iota, h$) up to non-unique $A_\infty$ isomorphism.

\begin{remark} If ${\bf M}$ is an $A_\infty$ module, then the induced $A_\infty$ structure on $H_*({\bf M})$ is constructed exactly as described in Proposition \ref{prop:InducedAinfty}, where the leaves and root of each tree have been labeled with $H_*({\bf M})$ rather than $H_*({\bf A})$, where appropriate.
\end{remark}

\begin{remark}\label{rmk:Unitality}
The condition $\iota p = \mbox{Id} + \partial h + h \partial$ in the statement of Proposition \ref{prop:InducedAinfty} is all that is needed to transfer the $A_\infty$ structure from ${\bf A}$ to $H_*({\bf A})$, while the extra condition $p\iota = \mbox{Id}$ ensures that the two structures are quasi-isomorphic. Moreover, although Proposition \ref{prop:InducedAinfty} as stated is a result about {\em non-unital} $A_\infty$ algebras (and non-unital modules over them), in the cases of interest in the present work (specifically, Lemmas \ref{lem:Bformal} and \ref{lem:PkHom}), Proposition \ref{prop:InducedAinfty} yields quasi-isomorphisms of {\em strictly unital} algebras (resp., modules).

Note also that if $H_*({\bf A})$ is finite-dimensional, the condition $p\iota = \mbox{Id}$ is a consequence of the condition $\iota p = \mbox{Id} + \partial h + h \partial$, hence may be omitted.
\end{remark}

\begin{definition} \label{defn:MinModelFormal} An $A_\infty$ structure on $H_*({\bf A})$ constructed as in Proposition \ref{prop:InducedAinfty} is called a {\em minimal model} of ${\bf A}$. An $A_\infty$ algebra is said to be {\em formal} if a minimal model can be chosen so that $m_n = 0$ for all $n>2$.
\end{definition}

Henceforth, whenever we refer to {\em the} minimal model, $H_*({\bf A})$, for ${\bf A}$ an $A_\infty$ algebra, we shall always assume it has been endowed with the structure provided by Proposition \ref{prop:InducedAinfty} for suitable maps $\iota, p, h$.


\begin{definition} Let ${\bf A}$ be a homologically unital $A_\infty$--algebra.  The {\em derived category ${\rm D}_\infty({\bf A})$} is the category with objects homologically unital $A_\infty$--modules (left, right, or bi--, depending on the context) and morphisms $A_{\infty}$--homotopy classes of $A_{\infty}$--morphisms.
\end{definition}

\begin{remark}
Since every $A_{\infty}$ quasi-isomorphism has an inverse up to homotopy (see \cite[Lemma 10.12.2.2]{bernstein}), passing to the derived category has the effect of making $A_{\infty}$ quasi-isomorphisms invertible.
\end{remark}

\begin{definition} \label{defn:FiltAlg} A (graded or ungraded) {\em filtered $A_\infty$ algebra} ${\bf A}$ is a (graded or ungraded) $A_\infty$ algebra equipped with a sequence of subsets, for $i \in \Z$:
\[0 \subseteq \ldots \subseteq \mathcal{F}_i \subseteq \mathcal{F}_{i+1} \subseteq \ldots \subseteq {\bf A}\] that are compatible with the $A_\infty$ structure in the following sense:
\[m_n\left(\mathcal{F}_{i_1} \otimes \ldots \otimes \mathcal{F}_{i_{n}}
\right) \subseteq \mathcal{F}_{i_1 + \ldots + i_n}.\]
If $m_n = 0$ for all $n>2$, ${\bf A}$ is a (graded or ungraded) {\em filtered differential algebra}.
(Graded or ungraded) {\em filtered $A_\infty$ modules} and {\em filtered differential modules} are defined analogously.
\end{definition}

Note that the compatibility of the filtration with the multiplicative structure ensures that if ${\bf A}$ is a filtered $A_\infty$ algebra, the associated graded algebra $\bigoplus_{i} \mathcal{F}_i/\mathcal{F}_{i-1}$ is a well-defined (graded or ungraded) $A_\infty$ algebra, and if ${\bf M}$ is a filtered $A_\infty$ module over a filtered $A_\infty$ algebra ${\bf A}$, then the associated graded module $\bigoplus_i \mathcal{F}_i/\mathcal{F}_{i-1}$ is a well-defined $A_\infty$ module over the associated graded algebra of ${\bf A}$.

\begin{definition}\label{defn:bounded}
A filtered $A_\infty$ algebra ${\bf A}$ (resp., module ${\bf M}$) is said to be {\em bounded} if there exist $n < N \in \mathbb{Z}$ such that $0 = \mathcal{F}_n$ and ${\bf A} = \mathcal{F}_{N}({\bf A})$ (resp., ${\bf M} = \mathcal{F}_{N}({\bf M})$). 
\end{definition}

\begin{notation}\label{defn:Filtshift}
If ${\bf M}$ is a filtered $A_\infty$ module and $k \in \Z$, ${\bf M}\{k\}$ will denote the filtered $A_\infty$ module whose filtration has been shifted by $k$.  Explicitly, 
\[\cF_n\left({\bf M}\{k\}\right) := \cF_{n-k}\left({\bf M}\right).\]
\end{notation}

A filtration on an $A_\infty$ algebra (resp., module) induces a spectral sequence in the standard way, and if the filtered complex is bounded this spectral sequence converges in a finite number of steps.  Furthermore, each page of the corresponding spectral sequence has the structure of an $A_\infty$ algebra (resp., module), by Proposition \ref{prop:InducedAinfty}.  We will call the homology of the associated graded complex,  $\bigoplus_{i \in \Z} \mathcal{F}_i/\mathcal{F}_{i-1}$, the {\em associated graded homology algebra} (resp., the {\em associated graded homology module}) and the homology of the total complex (i.e., the $E^\infty$ page of this spectral sequence) the {\em total homology algebra} (resp., the {\em total homology module}).

If {\bf M} is a filtered left $A_\infty$ {\bf A}-module, and {\bf N} is a filtered right $A_\infty$ {\bf B}-bimodule, then ${\bf M} \otimes {\bf N}$ inherits a filtration (and, hence, the structure of a filtered $A_\infty$ {\bf A}-{\bf B} bimodule in the sense of Definition \ref{defn:FiltAlg}) via: $a \otimes b \in \mathcal{F}_{m+n}\left({\bf M} \otimes {\bf N}\right)$ if $a \in \mathcal{F}_m\left({\bf M}\right)$ and $b \in \mathcal{F}_n\left({\bf N}\right)$.  

Similarly, the $A_\infty$ tensor product of filtered $A_\infty$ bimodules naturally inherits the structure of a filtered $A_\infty$ bimodule:

\begin{lemma} \label{lem:Ainftensorprod} Let ${\bf M}, {\bf N}$ be two filtered $A_\infty$ bimodules over a filtered $A_\infty$ algebra ${\bf A}$.  Then the $A_\infty$ tensor product, with underlying vector space:
\[{\bf M} \widetilde{\otimes} {\bf N} := \bigoplus_{n=0}^\infty {\bf M} \otimes {\bf A}^{\otimes n} \otimes {\bf N}\]  inherits the structure of a filtered $A_\infty$ bimodule as follows:
\[\mathcal{F}_\ell({\bf M} \widetilde{\otimes} {\bf N}) := \bigoplus_{n=0}^\infty\left[\bigoplus_{i + j_1 + \ldots + j_n + k = \ell} \mathcal{F}_i({\bf M}) \otimes \mathcal{F}_{j_1}({\bf A}) \otimes \ldots \otimes \mathcal{F}_{j_n}({\bf A}) \otimes \mathcal{F}_k({\bf N})\right]\]

\end{lemma}

\begin{proof}  Since ${\bf M}, {\bf N}$ are filtered $A_\infty$ bimodules, the multiplications
\begin{eqnarray*}
	m_{(0|1|i)}: {\bf M} \otimes {\bf A}_1 \otimes \ldots \otimes {\bf A}_i &\rightarrow& {\bf M}\\
	m_{(i|1|0)}: {\bf A}_{n-i+1} \otimes \ldots \otimes {\bf A}_n \otimes {\bf N} &\rightarrow& {\bf N}\\
	m_i: {\bf A}_{\ell} \otimes \ldots \otimes {\bf A}_{\ell+i-1} &\rightarrow& {\bf A}
\end{eqnarray*}
contributing to the differential on the complex all respect the filtration in the sense of Definition \ref{defn:FiltAlg}.  The same is true for the higher multiplications on the complex, for the same reason. 
\end{proof}

\begin{definition} \label{defn:Filtmorph}
An $A_\infty$ morphism $f: {\bf M} \rightarrow {\bf N}$ between two filtered $A_\infty$ modules is said to be {\em filtered} if \[f_{(n_1|1|n_2)}\left(\cF_{i_1} \otimes \ldots \otimes \cF_{i_{n_1+n_2+1}}\right) \subseteq \cF_{i_1 + \ldots + i_{n_1+n_2+1}}.\]
\end{definition}

\begin{definition} \label{defn:MC} Let ${\bf A}$ be a filtered $A_\infty$ algebra,  and $f:{\bf M} \rightarrow {\bf N}$ a filtered $A_\infty$ morphism between filtered ${\bf A}$--modules ${\bf M}$ and ${\bf N}$.   Let $m_{(n_1|1|n_2)}^M$ (resp., $m_{(n_1|1|n_2)}^N$) denote the $A_\infty$ multiplication maps for ${\bf M}$ (resp., for ${\bf N}$).  

Then the {\em mapping cone of $f$}, denoted $MC(f)$, is the filtered $A_\infty$ ${\bf A}$--module with underlying \F--vector space ${\bf M}[1] \oplus \left({\bf N}\right)$, $A_\infty$ multiplication maps:
\[m_{(n_1|1|n_2)} := \left(\begin{array}{cc}
	m^M_{(n_1|1|n_2)} & 0\\
	f_{(n_1|1|n_2)} & m^N_{(n_1|1|n_2)}
	\end{array}\right).\]
and filtration given by:
\[\cF_n(MC(f)) := \{(a,b) \in MC(f)\,\,\vline\,\, a \in \cF_n({\bf M}) \mbox{ and } b \in \cF_n({\bf N})\}.\]
\end{definition}

The following lemma will be useful in the proof of Theorem \ref{thm:filtdiffmodgen}.

\begin{lemma} \label{lem:assgrtensor} Let ${\bf M} \widetilde{\otimes} {\bf N}$ be the filtered $A_\infty$ bimodule (over the filtered algebra ${\bf A}$) obtained as the $A_\infty$ tensor product of the two filtered $A_\infty$ bimodules ${\bf M}$ and ${\bf N}$ as in Lemma \ref{lem:Ainftensorprod}.  Let $\mbox{gr}(-)$ denote the associated graded $A_\infty$ module of $-$.

Then $\mbox{gr}({\bf M}) \widetilde{\otimes}_{\mbox{gr}({\bf A})} \mbox{gr}({\bf N}) = \mbox{gr}({\bf M} \widetilde{\otimes}_{\bf A} {\bf N})$ as $A_\infty$ bimodules over $gr({\bf A})$.
\end{lemma}

\begin{proof} 

The well-defined chain map 
\begin{eqnarray*}
	\xymatrix{\mbox{gr}({\bf M}) \widetilde{\otimes}_{\mbox{gr}({\bf A})} \mbox{gr}({\bf N}) \ar[r]^{\hskip12pt \Phi} &\mbox{gr}({\bf M} \widetilde{\otimes}_{\bf A}{\bf N})} \\
\end{eqnarray*}

sending 
\begin{eqnarray*}
	[x] \otimes [a_1] \otimes \ldots \otimes [a_n] \otimes [y] &\in& \frac{\mathcal{F}_i}{\mathcal{F}_{i-1}}({\bf M}) \otimes \frac{\mathcal{F}_{j_1}}{\mathcal{F}_{j_1-1}}({\bf A}) \otimes \ldots \otimes \frac{\mathcal{F}_{j_n}}{\mathcal{F}_{j_n-1}}({\bf A}) \otimes \frac{\mathcal{F}_{k}}{\mathcal{F}_{k-1}}({\bf N})\\
	&\subseteq& \mbox{gr}({\bf M}) \widetilde{\otimes}_{\mbox{gr}({\bf A})} \mbox{gr}({\bf N}).
\end{eqnarray*}

to \[\left[x \otimes a_1 \otimes \ldots \otimes a_n \otimes y\right] \in \frac{\mathcal{F}_{I}}{\mathcal{F}_{I-1}}({\bf M} \widetilde{\otimes}_{\bf A} {\bf N})\]  (where in the above $I := i+ j_1 \ldots + j_n + k$), is an isomorphism of chain complexes.
This map is well-defined, since any other representative, $x' \otimes a_1' \otimes \ldots \otimes a_n' \otimes y' \in {\bf M} \widetilde{\otimes}_{\bf A} {\bf N},$ of $[x] \otimes [a_1] \otimes \ldots \otimes [a_n] \otimes [y]$ will differ from $x \otimes a_1 \otimes \ldots \otimes a_n \otimes y$ by an element in $\mathcal{F}_{I-1}$, by the definition of the filtration on ${\bf M}\widetilde{\otimes}_{\bf A} {\bf N}$.

Similarly, we send an equivalence class $[x \otimes a_1 \otimes \ldots \otimes a_n \otimes y] \in \mbox{gr}({\bf M} \widetilde{\otimes}_{\bf A} {\bf N})$ to the uniquely-specified equivalence class 
\begin{eqnarray*}
	\Psi([x \otimes a_1 \otimes \ldots \otimes a_n \otimes y]) &:=& [x] \otimes [a_1] \otimes \ldots \otimes [a_n] \otimes [y] \\
	&\in& \mbox{gr}({\bf M}) \widetilde{\otimes}_{\mbox{gr}({\bf A})} \mbox{gr}({\bf N}).
\end{eqnarray*}

These maps are well-defined, and the differentials on the two sides can be easily seen to agree.
Furthermore, the differentials on $\mbox{gr}({\bf M} \widetilde{\otimes}_{\bf A} {\bf N})$ and $\mbox{gr}({\bf M}) \widetilde{\otimes}_{\mbox{gr}({\bf A})} \mbox{gr}({\bf N})$ agree, by the same argument above applied to the image of the differential of a representative $x \otimes a_1 \otimes \ldots \otimes a_n \otimes y \in {\bf M} \widetilde{\otimes}_{\bf A} {\bf N}$.
\end{proof}

\subsection{Formality of dg algebras and modules}
The following technical lemmas provide sufficient (but not necessary) conditions for formality of an $A_\infty$ module.

\begin{lemma}  \label{lem:FormalCond} Let ${\bf A}$ be a differential (graded) algebra (resp., let ${\bf M}$ be a differential (graded) module over ${\bf A}$), and let $\iota, p, h$ be maps satisfying the conditions in Proposition \ref{prop:InducedAinfty}.  If, in addition,
\begin{enumerate}
	\item $h^2 = h\iota = 0$, and
	\item $m_2^A(\iota \otimes \iota)({\bf A}^{\otimes 2}) \subseteq \iota({\bf A})$ (resp., $m_{(n_1|1|n_2)}^M(\iota \otimes \iota)({\bf A}^{\otimes n_1} \otimes {\bf M} \otimes {\bf A}^{\otimes n_2}) \subseteq \iota({\bf M})$ whenever $n_1+1+n_2 = 2$),
\end{enumerate}
then ${\bf A}$ is formal (resp., ${\bf M}$ is formal).

Furthermore, $\iota: {\bf A} \rightarrow H_*({\bf A})$ (resp., $\iota: {\bf M} \rightarrow H_*({\bf M})$) is a strict $A_\infty$ quasi-isomorphism.
\end{lemma}

\begin{proof}
In the interest of brevity, we give the argument for the case of ${\bf A}$ a differential (graded) algebra, leaving the completely analogous proof in the case of ${\bf M}$ a differential (graded) module to the reader.

Each tree $T$ contributing to the definition of \[m_n: \left(H_*({\bf A})\right)^{\otimes n} \rightarrow H_*({\bf A})\] for $n>2$ 
yields the $0$ map, since each such tree $T$ involves a product of terms in ${\bf A}$, 
 at least one of which is either:
\begin{itemize}
	\item of the form $h\circ m_2^{{A}} \circ (\iota \otimes \iota)$ (if $T$ is trivalent) 
	or
	\item of the form $m_n^{{A}}(\iota \otimes \ldots \otimes \iota)$, 
	for $n>2$  (if $T$ is not trivalent).
\end{itemize}
In both cases, such a term is $0$ in ${\bf A}$ 
by assumption, hence the corresponding map is $0$, implying formality of ${\bf A}$.

To see that $\iota: {\bf A} \rightarrow H_*({\bf A})$ is a {\em strict} quasi-isomorphism, we refer to \cite[Thm. 2.1]{Berglund}, which  tells us that $\iota_n$ can be defined recursively as \[\iota_n := \sum_{\substack{i_1 + \ldots + i_r = n\\r>1}} hm_r^A\left(\iota_{i_1} \otimes \ldots \otimes \iota_{i_r}\right).\]

Assumptions (1) and (2), combined with the assumption that $m_r^A = 0$ for $r > 2$, now allow us to conclude inductively that $\iota_n = 0$ for $n \geq 2$, as desired.


\end{proof}

\begin{lemma} \label{lem:FormalCond2} Let {\bf M} be a differential (graded) module over an algebra {\bf A}, and let
\[\iota_M: H_*({\bf M}) \rightarrow {\bf M}, \hskip 10pt p_M: {\bf M} \rightarrow H_*({\bf M}), \hskip 10pt h_M: {\bf M} \rightarrow {\bf M}\] satisfy the conditions in Proposition \ref{prop:InducedAinfty}.  Suppose in addition that 
\begin{enumerate}
	\item $h_M^2 = p_Mh_M = 0$, and
	\item $\mbox{Im}(h_M)$ and $\mbox{Im}(m^M_{(0|1|0)})$ are both submodules of {\bf M} over {\bf A} (i.e., left or/and right multiplication by an element of {\bf A} preserves Im$\left(h_M\right)$ and Im$(m^M_{(0|1|0)})$).
\end{enumerate}
Then {\bf M} is formal, and the projection map $p_M:{\bf M} \rightarrow H_*({\bf M})$ is a strict quasi-isomorphism.
\end{lemma}
\begin{proof}
We give the proof in the case that ${\bf M}$ is a differential (graded) bimodule over ${\bf A}$.  If Assumption (2) holds only under left (resp., right) multiplication, then $p_M$ will be a strict quasi-isomorphism of left (resp., right) {\bf A}--modules.

Since ${\bf A}$ is an algebra, $m^A_{n} = 0$ unless $n=2$, and ${\bf A}$ is trivially $A_\infty$ isomorphic to its homology.  Choosing $\iota_A: H_*({\bf A}) \rightarrow {\bf A}$ and $p_A: {\bf A} \rightarrow H_*({\bf A})$ to be the identity morphism, and $h_A: {\bf A} \rightarrow {\bf A}$ to be the zero morphism, we now claim that any tree $T$ contributing to the definition of \[m_{(n_1|1|n_2)}: {\bf A}^{\otimes n_1} \otimes (H_*({\bf M})) \otimes {\bf A}^{\otimes n_2} \rightarrow H_*({\bf M})\] is zero if $n_1+n_2+1>2$.  This follows because:
\begin{itemize} 
	\item If $T$ is trivalent then it corresponds to a summand of the form $p_M\circ h_M (m)$, since Im($h_M$) is an ${\bf A}$--bimodule.  Such a term is zero by Assumption (1).
	\item If $T$ is not trivalent then it involves a product with at least one term of the form:
	\[m^M_{(n_1'|1|n_2')}(\iota \otimes \ldots \otimes \iota) \mbox{ (resp., } m^A_n(\iota \otimes \ldots \otimes \iota))\] for $n_1'+n_2'+1 >2$ \mbox{ (resp., $n>2$)}, which is zero since ${\bf M}$ is a dg module (resp., since ${\bf A}$ is an algebra).
\end{itemize}


To see that $p_M$ is a strict quasi-isomorphism, we again appeal to \cite[Thm. 2.1]{Berglund}, which gives recursive definitions for $\left(p_M\right)_{(n_1|1|n_2)}$ in terms of $p$, $m$, and an auxiliary morphism $h^{[n_1|1|n_2]}$, defined recursively in terms of $p, \iota, h$.



Assumptions (1) and (2) allow us to conclude:
\begin{eqnarray*}
	\left(p_M\right)_{(1|1|0)} &:=& p_{(0|1|0)} \circ m_{(1|1|0)} \circ (1 \otimes h)\\
	&=& 0
\end{eqnarray*}
and
\begin{eqnarray*}
	\left(p_M\right)_{(0|1|1)} &:=& p_{(0|1|0)} \circ m_{(0|1|1)} \circ (h \otimes 1)\\
	&=& 0.
\end{eqnarray*}

Combined with the fact that $m_n^A = 0$ for all $n \neq 2$ and $m^M_{(n_1|1|n_2)} = 0$ for all $n_1 + 1+n_2>2$, $(p_M)_{(n_1|1|n_2)}$ is then identically $0$ by induction for all $(n_1+1+n_2) \geq 2$, as desired.


\end{proof}

\section{Khovanov-Seidel Hom algebras and bimodules} \label{sec:Khalgmod}
In this section, we construct dg bimodules following Khovanov-Seidel in \cite{MR1862802}.  We begin by describing the topological data needed for the construction of both the Khovanov-Seidel bimodules and their bordered Floer analogues (described in Section \ref{sec:HFalgmod}).

We emphasize that although we have chosen to describe the Khovanov-Seidel objects from a purely algebraic viewpoint, they also admit a beautiful Fukaya-theoretic description (cf. Section \ref{sec:BKhFukaya}). Readers familiar with \cite[Sec. 6]{MR1862802} and \cite[Chp. 20]{SeidelBook} will likely benefit from keeping this geometric picture in mind.

\subsection{Topological data: Basis of curves}\label{sec:pictures} 

Let $D_m$ denote the unit disk in the complex plane, equipped with a set,  \[\Delta := \left\{-1 +\frac{2(j+1)}{m+2} \in D_m \subset \C \,\,\vline\,\, j= 0, \ldots, m\right\},\] of $m+1$ points equally distributed along the intersection of the real axis with $D_m$.  Label by ${\bf j}$ the point at position $-1 + \frac{2(j+1)}{m+2}$.  

By a {\em curve in $D_m$} we shall always mean the image of a smooth imbedding $\gamma:[0,1] \rightarrow D_m$ which is transverse to $\partial D_m$ and satisfies $\gamma^{-1}(\partial D_m \cup \Delta) = \{0,1\}.$

\begin{definition} \label{def:AdmCurve} A {\em $\partial$--admissible curve in $D_m$} is a curve in $D_m$ for which $\gamma(0) = -1$ and $\gamma(1) \in \Delta$.
\end{definition}

A $\partial$--admissible curve is a particular type of admissible curve in the sense of \cite[Sec. 3b]{MR1862802}.  Two $\partial$--admissible curves $c_1$ and $c_2$ are said to be isotopic if there is a homotopy between $c_1$ and $c_2$ through $\partial$--admissible curves.  

\begin{notation} \label{not:bigrading} Associated to any curve, $c \subset D_m$, is a canonical section of the interior of $c$ to the real projectivization of the tangent bundle of $D_m\setminus\Delta$.  By choosing a lift of this section to a particular $\Z^2$ cover as described in \cite[Sec. 3d]{MR1862802}, one assigns a bigrading to $c$.  We shall denote by $\widetilde{c}$ the data of a curve $c \subset D_m$ equipped with such a choice of bigrading.
\end{notation}

\begin{definition} \cite[Sec. 3a]{MR1862802} Two curves $c_0, c_1 \subset D_m$ are said to have {\em minimal geometric intersection} if they satisfy the following conditions:
\begin{itemize}
	\item $c_0$ and $c_1$ intersect transversely,
	\item $c_0 \cap c_1 \cap \partial D_m = \emptyset$, and
	\item If $z_- \neq z_+$ are two points in $c_0 \cap c_1$ not both in $\Delta$, $\alpha_0 \subset c_0$ and $\alpha_1 \subset c_1$ are two arcs with endpoints $z_-, z_+$ such that $\alpha_0 \cap \alpha_1 = \{z_-,z_+\}$, and $K$ is the connected component of $D_m - (c_0 \cup c_1)$ bounded by $\alpha_0 \cup \alpha_1$, then if $K$ is topologically an open disk, it must contain at least one point of $\Delta$.  Informally, we say there are no ``trivial bigons" among the connected components of $D_m - (c_0 \cup c_1)$.
\end{itemize}
\end{definition}

\begin{figure}
\begin{center}
\resizebox{2in}{!}{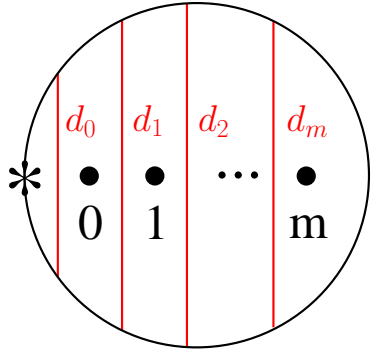}
\end{center}
\caption{The curves $d_j$, for $j=0, \ldots, m$, are the intersections of the lines $Re(z) = \left(-1-\frac{1}{m+2}\right) + \frac{2(j+1)}{m+2}$ with the unit disk in $\mathbb{C}$.  By convention, the distinguished point, labeled by a $*$, at $-1 \in \partial D_m$, is the left endpoint for all $\partial$--admissible curves in $D_m$.}
\label{fig:DiskDm}
\end{figure}

\begin{definition} \cite[Sec. 3e]{MR1862802}
Let $d_0, \ldots, d_m \subset D_m$ be the curves pictured in Figure \ref{fig:DiskDm}.  A $\partial$--admissible curve in $D_m$ is said to be in {\em normal form} if it has minimal geometric intersection with $d_j$ for each $j=0, \ldots, m$.
\end{definition}

\begin{definition}
A {\em basis of $\partial$--admissible curves in $D_m$} is a set,
$\cB = \left\{c_0, \ldots, c_m\right\},$ of $\partial$--admissible curves satisfying the conditions:
\begin{itemize}
	\item If $\gamma_j: [0,1] \rightarrow D_m$ is the imbedding whose image is $c_j$, then $\gamma(1) = {\bf j} \in \Delta$ (the right endpoint of $c_j$ is {\bf j}), and
	\item $c_i \cap c_j = \{-1\}$ if $i\neq j$ (distinct curves $c_i$ and $c_j$ intersect only at their left endpoints).
\end{itemize}

If we, furthermore, specify a lift of each curve, $c_j \in \cB$, to a bigraded curve, $\widetilde{c}_j$, we say that we have a {\em basis, $\widetilde{\cB} = \left\{\widetilde{c}_0, \ldots, \widetilde{c}_m\right\}$, of $\partial$--admissible bigraded curves in $D_m$}.
\end{definition}

Unless otherwise specified, from this point forward whenever we write that $\widetilde{\cB}$ is a {\em basis}, we shall always mean that $\widetilde{\cB}$ is a {\em basis of $\partial$--admissible bigraded curves in normal form in $D_m$}.  Two bases $\cB = \{\widetilde{c}_0, \ldots \widetilde{c}_m\}$ and $\cB'= \{\widetilde{c}_0', \ldots, \widetilde{c}_m'\}$ are said to be {\em equivalent} if there exists an isotopy $\widetilde{c}_i \rightarrow \widetilde{c}_i'$ for each $i = 0, \ldots, m$ through $\partial$--admissible bigraded curves in normal form.

As in \cite{MR1862802}, we let $\mathcal{G} = \mbox{Diff}(D_m, \partial D_m;\Delta)$ denote the group of diffeomorphisms $f$ of $D_m$ satisfying $f|_{\partial D_m} = \mbox{Id}$ and $f(\Delta) = \Delta$ and note that there is a canonical identification of $\pi_0(\cG)$ with $B_{m+1}$, the Artin braid group on $m+1$ generators.  Under this correspondence, (isotopy classes of) $\partial$--admissible curves are sent to (isotopy classes of) $\partial$--admissible curves.  Moreover, an (equivalence class of) basis $\widetilde{\cB}$ is sent to an (equivalence class of) basis $\sigma(\widetilde{\cB})$, after suitably reordering the curves in $\sigma(\widetilde{\cB})$.

 \subsection{The ring $A_m$ and a braid group action on $D^b(A_m)$} \label{sec:Am}
In  \cite{MR1862802}, Khovanov-Seidel associate to a braid, $\sigma \in B_{m+1}$, a bimodule over a quiver algebra, $A_m$ (defined below).  In this subsection, we explain how their construction yields a family of algebras and bimodules, one for each choice of basis.  
Our end goal is the construction of a particular algebra, $B^{Kh}$, and a bimodule, $\cM_\sigma^{Kh}$ over $B^{Kh}$, from the data of a particular such basis, $\widetilde{\mathcal{Q}}$.

\begin{figure}
\begin{center}
\resizebox{2in}{!}{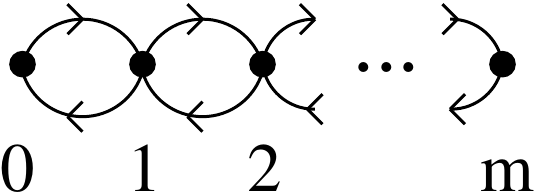}
\end{center}
\caption{}
\label{fig:QuiverDm}
\end{figure}

We begin by reviewing the original construction of Khovanov-Seidel in \cite{MR1862802}.  Let $\Gamma_m$ be the oriented graph (quiver) whose vertices are labeled ${\bf 0}, \ldots, {\bf m}$ and whose edges are shown in Figure \ref{fig:QuiverDm}.  Recall that, given any oriented graph $\Gamma$, one defines its path ring as the vector space over $\mathbb{F}$ freely generated by the set of all finite-length paths in $\Gamma$, where multiplication is given by concatenation, and the product of two non-composable paths is set to $0$.  The ring $A_m$ is then defined as a quotient of the path ring of $\Gamma_m$ by the collection of relations \[(i-1|i|i+1) =0, \hskip 10pt (i+1|i|i-1)=0, \hskip 10pt (i|i+1|i)=(i|i-1|i), \hskip 10pt (0|1|0)=0\] for each $0 < i < m$.  In the above, following  \cite{MR1862802}, we have labeled each path in $\Gamma_m$ by the complete ordered tuple of vertices it traverses.  So, for instance, $(i-1|i|i+1)$ denotes the path that starts at vertex $i-1$, moves right to $i$, then right again to $i+1$.  The path ring of $\Gamma_m$ is further endowed with a grading by setting $\mbox{deg}(i) = \mbox{deg}(i|i+1) = 0$ and $\mbox{deg}(i|i-1) = 1$ for all $i$.  This grading descends to the quotient, $A_m$, since the relations defining $A_m$ are homogeneous with respect to the grading.\footnote{This {\em internal grading} corresponds to the second of the two gradings discussed in Notation \ref{not:bigrading}.  Note that this grading is {\em not} the grading by path length which appears in \cite{QA0610054, StroppelAlgebra} and corresponds to the $j$ (quantum) grading of \cite{MR1740682}.  See Remark \ref{rmk:gradings}.}

Note that the collection $\{(i)|i \in 0,\ldots,m\}$ of constant paths are mutually orthogonal idempotents, and $\sum_{i=0}^m (i)$ is the identity in $A_m$.  There are corresponding decompositions of $A_m$ as a direct sum of projective left modules $A_m = \bigoplus_{i=0}^m A_m(i)$ (resp., projective right-modules $A_m = \bigoplus_{i=0}^m (i)A_m$).  As in  \cite{MR1862802}, we denote $A_m(i)$ (resp., $(i)A_m$) by $P_i$ (resp., ${}_i{P}$).  Note that $P_i$ (resp., ${}_i{P}$) is the set of all paths ending at $i$ (resp., beginning at $i$).

To streamline notation, we henceforth assume that we have fixed $m \geq 0 \in \Z$, and let $A$ denote the algebra $A_m$.

\begin{figure}
\begin{center}
\resizebox{2in}{!}{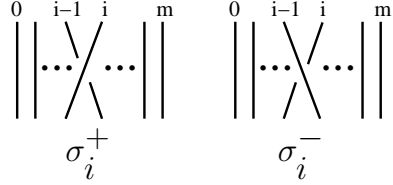}
\end{center}
\caption{The elementary Artin generators, $\sigma_i^\pm$}
\label{fig:Braidgen}
\end{figure}

Khovanov-Seidel go on to associate to each braid $\sigma \in B_{m+1}$ an element of $D^b(A)$, the bounded derived category of $A$--bimodules, by associating to each elementary Artin braid generator $\sigma_i^{\pm 1}$ (pictured in Figure \ref{fig:Braidgen}) a dg bimodule $\cM_{{\sigma_i^\pm}}$ and to each braid, $\sigma := {\sigma_{i_1}}^\pm \cdots {\sigma_{i_k}}^\pm$, decomposed as a product of elementary braid words, the dg bimodule \[\cM_\sigma = \cM_{\sigma_{i_1}\pm} \otimes_{A} \ldots \otimes_{A} \cM_{\sigma_{i_k}^\pm}.\]  They then verify that any two decompositions of $\sigma$ as a product of elementary Artin braid generators give rise to quasi-isomorphic complexes, and hence $\cM_\sigma$ gives rise to a well-defined element in $D^b(A)$.

\subsection{The dg algebra $B$ and the algebra $B^{Kh}$}
Now, suppose we are given the data of a $\partial$--admissible bigraded curve in normal form.   Khovanov-Seidel show, in \cite[Sec. 4]{MR1862802}, how to use this data to construct a bounded complex of bigraded projective left modules over the algebra $A$.  Furthermore, a basis, $\widetilde{\cB}$, of such curves yields a dga via Yoneda imbedding (cf. \cite[Sec. 2.6]{MR2258042}).  Recall:

\begin{definition} Let $(\mathcal{C}_1,\partial_1), (\mathcal{C}_2,\partial_2)$ be two bounded dg left modules over an algebra ${\bf A}$.  Then the {\em Hom complex of the pair $(\mathcal{C}_1, \mathcal{C}_2)$}, denoted $\Hom_{\bf A}(\mathcal{C}_1, \mathcal{C}_2)$, is the bounded complex whose generators are left module morphisms, $F:\mathcal{C}_1 \rightarrow \mathcal{C}_2$, and whose differential, $D$, is given by 
\[D(F) := \partial_{2}F + F\partial_1.\]
\end{definition}

\begin{construction} \label{const:Homalgebra} Let $\widetilde{\cB} = \{\widetilde{c}_0, \ldots, \widetilde{c}_m\}$ be a basis, and let $L(\widetilde{c}_j)$ be the bounded complex of projective $A$--modules associated to $\widetilde{c}_j$, for each $j=0, \ldots, m$.  Then the direct sum,
\[\bigoplus_{i,j=0}^m \Hom_{A}(L(\widetilde{c}_i),L(\widetilde{c}_j)),\] is a dga, with multiplication given by composition of $A$--module morphisms.  We will refer to $\bigoplus_{i,j=0}^m \Hom_{A}(L(\widetilde{c}_i),L(\widetilde{c}_j))$ as the {\em Hom algebra associated to $\widetilde{\cB}$}.
\end{construction}

\begin{figure}
\begin{center}
\resizebox{2in}{!}{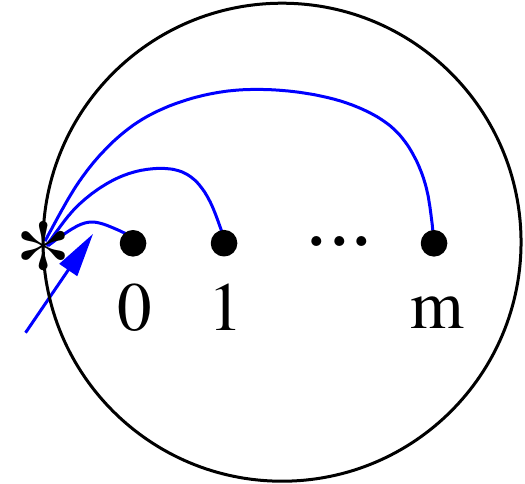}
\end{center}
\caption{The basis $\mathcal{Q} = \left\{q_0, \ldots, q_m\right\}$}
\label{fig:Qbasis}
\end{figure}

We focus in the present paper on the Hom algebra associated to the basis $\widetilde{\mathcal{Q}} = \left\{\widetilde{q}_0, \ldots, \widetilde{q}_m\right\}$ given by (a particular lift of) the collection of curves pictured in Figure \ref{fig:Qbasis}.\footnote{We expect that results similar to those described in Theorems \ref{thm:filtdiffalg} and \ref{thm:filtdiffmodgen} hold for other choices of basis, but we do not address that here.}

Applying the construction of \cite[Sec. 4a]{MR1862802}, we associate to $\widetilde{q}_j$ the dg module:
\[Q_j := \xymatrix{0 \ar[r] & P_0 \ar[r]^{\cdot (0|1)} & P_1 \ar[r]^{\cdot (1|2)} & \ldots \ar[r]^{\cdot (j-1|j)} & P_j \ar[r] & 0},\] where the differential map ``$\,\,\cdot (i-1|i)$" denotes ``right multiplication by the element $(i-1|i)$."  By fixing a lift of the tangent vector to the curve $q_0$ at a point near ${\bf 0} \in \Delta$ and declaring this lift to correspond to bigrading $(0,0)$, we obtain a ``canonical" bigrading on $Q_j$ satisfying the property that the bigrading of the idempotent $(i) \in P_i$ is $(i,0)$.\footnote{In the language of \cite{BrundanStroppel}, $Q_j$ is a projective resolution of the {\em standard} module associated to the length $m+1$ weight $\lambda = (\vee\ldots\vee\wedge\vee\ldots\vee)$, where the lone $\wedge$ is in the $(j \in \{0,\ldots, m\})$th position.}


\begin{notation} \label{defn:KhHomalgmod} We shall denote by $B$ the Hom algebra associated to $\widetilde{\mathcal{Q}}$:
\[\bigoplus_{i,j=0}^m\Hom_{A}\left(Q_i,Q_j\right)\] and by $B^{Kh}$ its homology, $H_*(B)$, considered as an $A_\infty$ algebra via the construction in Proposition \ref{prop:InducedAinfty}.

\end{notation}

We will eventually be interested in $D_\infty(B^{Kh})$--in particular, a braid group action on this category--so we now devote some time to describing the structure of $B$ and $B^{Kh}$.

\begin{notation} Let $R_{\cI}$ be a bounded complex of elementary projective left $A$--modules (e.g., one obtained from an admissible curve in normal form in $D_m$ as explained in \cite[Sec. 4]{MR1862802}): 
\[R_{\cI} = 0 \rightarrow P_{i_0}\{s_0\} \rightarrow \ldots \rightarrow P_{i_N}\{s_N\} \rightarrow 0.\]

Suppose further that $P_{i_0}\{s_0\}$ is in (co)homological grading $0$.  Then we will use the notation ${}_{\cI}{R}$ to denote the following bounded complex of elementary projective {\em right} $A$--modules:
\[{}_\cI{R} := 0 \leftarrow {}_{i_0}P\{-s_0\}[0] \leftarrow \ldots \leftarrow  {}_{i_N}P\{-s_N\}[-N] \leftarrow 0,\]
where, if a map $P_{i_j} \rightarrow P_{i_{j+1}}$ in $R_{\cI}$ is given by right multiplication by a path $\gamma \in A$, then the corresponding map ${}_{i_j}P \leftarrow {}_{i_{j+1}}P$ in ${}_{\cI}{R}$ is given by {\em left} multiplication by $\gamma$.
\end{notation}

\begin{lemma} \label{lem:Homcpx} Let $R_{\cI}, S_{\cJ}$ be bounded complexes of elementary projective left $A$--modules as above.  Then $\Hom_A(R_{\cI},S_{\cJ}) \cong {}_\cI{R} \otimes_A S_\cJ$.
\end{lemma}

\begin{proof} Each element $\phi \in \Hom_A(R_\cI, S_\cJ)$ can be decomposed as a sum of left $A$--module maps $\phi_{k,\ell}: P_{i_k}\{s_k\} \rightarrow P_{j_\ell}\{s_\ell\}$, each of which is uniquely determined by the image, $\phi_{k,\ell}(i_k)$, of the idempotent, $(i_k)$.  We therefore obtain an isomorphism 
\[\Hom_A(R_\cI,S_\cJ) \rightarrow {}_\cI{R} \otimes_A S_{\cJ}\] of $\mathbb{F}$--vector spaces identifying $\phi$ with the element, $\sum_{k,\ell} \left((i_k) \otimes \phi_{k,\ell}(i_k)\right)$.

To see that the Hom complex differential $D(\phi) := \phi d_{\cI} + d_{\cJ}\phi$ on the left matches the tensor product differential on the right, we simply note that if $\phi = \sum_{k,\ell} \phi_{k,\ell} \in \Hom_A(R_\cI, S_\cJ)$, then for each pair, $(k,\ell)$, $\phi_{k,\ell} d_\cI$ is obtained by pre- (i.e., left-) (resp., $d_\cJ \phi_{k,\ell}$ is obtained by post- (i.e., right-)) multiplying $\phi_{k,\ell}$ by a path $\gamma_k$ (resp., $\gamma_{\ell}$).  This is precisely the induced differential on the tensor product complex ${}_\cI R \otimes_A S_{\cJ}$.

\end{proof}

\begin{lemma} \label{lem:degree}
Let $R_\cI, S_\cJ$ be two bigraded bounded complexes of projective modules obtained from admissible bigraded curves in normal form as explained in \cite[Sec. 4]{MR1862802}.  Then the differential on $\Hom_A(R_\cI, S_\cJ)$ has degree $(1,0)$.
\end{lemma}

\begin{proof}
By definition, the differential on each of $R_\cI$, $S_\cJ$ has degree $(1,0)$, implying that the differential on ${}_\cI R$ and, hence, the differential on \[\Hom_A(R_\cI, S_\cJ) = {}_\cI R \otimes_A S_\cJ,\] has degree $(1,0)$ as well.
\end{proof}

The following lemma was also obtained independently by Klamt and Stroppel.  Compare \cite[Thms. 5.7, 7.3]{Klamt} and \cite[Thms. 5.3, 6.7]{KlamtStroppel}.

\begin{lemma} \label{lem:twoisomcpx} \label{lem:Bformal} The dg algebra $B:= \bigoplus_{i,j=0}^m \Hom_A\left(Q_i,Q_j\right)$ is formal.  Furthermore, the algebra \[B^{Kh} := H_*(B)\] has the following explicit description:  \[B^{Kh} := \bigoplus_{i,j = 0}^m {}_iB^{Kh}_j, \mbox{ with}\] 
\[{}_iB^{Kh}_j := 
		\left\{\begin{array}{cl}
			0 & \mbox{if $i<j$,}\\
			\mbox{Span}_\F\langle {}_i\Id_j \rangle & \mbox{if $i=j$, and}\\
			\mbox{Span}_\F\langle {}_i\Id_j, {}_i{\bf x}_j \rangle & \mbox{if $i>j$,}
		\end{array}\right.\]
where the bigradings on generators are given by: 
\begin{eqnarray*}
	gr\left({}_i\Id_j\right) &=& (0,0) \,\,\,\,\mbox{ for all $i\geq j$,}\\
	gr\left({}_i{\bf x}_j\right) &=& (-1,1) \,\,\,\,\mbox{ for all $i > j$}.
\end{eqnarray*} and the multiplication is given by:
\begin{eqnarray*}
	m_2({}_i\Id_j \otimes {}_j\Id_k) &:=& {}_i\Id_k\\
	m_2({}_i\Id_j \otimes {}_j{\bf x}_k) &:=& {}_i{\bf x}_k\\
	m_2({}_i{\bf x}_j \otimes {}_j\Id_k) &:=& {}_i{\bf x}_k\\
	m_2({}_i{\bf x}_j \otimes {}_j{\bf x}_k) &=& 0
\end{eqnarray*}
(As usual, $m_2: {}_iB^{Kh}_j \otimes {}_kB^{Kh}_\ell \rightarrow {}_iB^{Kh}_\ell$ is identically $0$ when $j \neq k$.)
\end{lemma}

\begin{proof}
We know from \cite[Prop. 4.9] {MR1862802} that as an {\F}--vector space, ${}_iB^{Kh}_j$ is free of rank  $0$ when $i<j$, $1$ when $i=j$, and $2$ when $i>j$.

Indeed, one sees by direct calculation that when $i<j$ the chain complex splits as the direct sum of two acyclic subcomplexes. When $i=j$, the chain complex splits in a similar fashion, but the first of the two complexes has homology generated by $(0) + \ldots + (j)$ and the second is acyclic. 
When $i > j$, the chain complex again splits, but now both subcomplexes have non-trivial homology, the first generated by $(0) + \ldots + (j)$, and the second generated by $(1|0) + \ldots + (j+1|j)$.

Denote the first (resp., second) subcomplex by $\mathcal{C}_{\Id}$ (resp., by $\mathcal{C}_{\bf x})$.

Proposition \ref{prop:InducedAinfty} now guarantees that $B^{Kh}:= H_*\left(B\right)$ admits an $A_\infty$ structure quasi-isomorphic to $B$, which we may describe explicitly once we have maps $p: {}_iB_j \rightarrow {}_iB^{Kh}_j$, $\iota: {}_iB^{Kh}_j \rightarrow {}_iB_j$, $h: {}_iB_j \rightarrow {}_iB_j$ satisfying the assumptions of Proposition \ref{prop:InducedAinfty}. We describe these maps in the case $i>j$, leaving the completely analogous cases $i \leq j$ to the reader.  

The inclusion map $\iota$ is the \F--linear extension of: 
\begin{itemize}
	\item $\iota({}_i\Id_j) := (0) + \ldots + (j)$,
	\item $\iota({}_i{\bf x}_j) := (1|0) + \ldots + (j+1|j)$.
\end{itemize}

With respect to the bases:
\[\{(0) + \ldots + (k) \,\,|\,\, 0 \leq k \leq j\} \cup \{(k-1|k) \,\,|\,\, 0 \leq k \leq j\}\] for $\mathcal{C}_{\Id}$,  and: \[\{(1|0) + \ldots + (k+1|k)\,\,|\,\,0 \leq k \leq j\} \cup \{(k|k-1|k)\,\,|\,\,0 \leq k \leq j\}\] for $\mathcal{C}_{\bf x}$, the projection map $p$ is the \F--linear extension of: 
\[p(\phi) := \left\{\begin{array}{cl} 
			{}_i\Id_j & \mbox{if $\phi = (0) + \ldots + (j)$,}\\
			{}_i{\bf x}_j & \mbox{if $\phi = (1|0) + \ldots + (j+1|j)$, and}\\
			0 & \mbox{otherwise.}\end{array}\right.\]
The homotopy map $h$ is the \F--linear extension of: 
\[h(\phi) := \left\{\begin{array}{cl}
			\left(m_1^B\right)^{-1}(\phi) & \mbox{if $\phi \in$ Im($m_1^B$)}\\
			0 & \mbox{otherwise,}\end{array}\right.\]
where in the above, $\left(m_1^B\right)^{-1}(\phi)$ is defined to be the (unique) basis element $\phi'$ satisfying $\partial(\phi') = \phi$.

One can now either see directly that $B$ is formal by applying Lemma \ref{lem:FormalCond} or simply note that the sum of the two gradings associated to each element in $B^{Kh}$ is $0$. As each structure map $m_n: \left(B^{Kh}\right)^{\otimes n} \rightarrow B^{Kh}$ is degree $(2-n)$ on this sum, nontrivial multiplications are only possible when $n=2$.

Verification that the bigradings and multiplication are as stated is a straightforward calculation.
\end{proof}

\begin{remark} \label{rmk:BKhMatrix} The algebra $B^{Kh}$ is isomorphic to the algebra of lower triangular $(m+1)\times (m+1)$ matrices over $\mathbb{F}[{\bf x}]/({\bf x}^2)$ with only $0$'s and $1$'s on the main diagonal:
\[
B^{Kh}\cong
\left\{
\left.
\left(\begin{matrix}
d_0        & 0      & \ldots         & 0     \\
\phi_{1,0} & d_1    & \ddots         & \vdots\\
\vdots     & \ddots & \ddots         & 0     \\
\phi_{m,0} & \ldots & \phi_{m,m-1} & d_m  
\end{matrix}\right)
\right |\,d_i\in\{0,1\}
\right\}\subset M_{m+1}(\mathbb{F}[{\bf x}]/({\bf x}^2))
\]

We define an algebra isomorphism by sending the generator ${}_i\mathbb{1}_j\in{}_iB^{Kh}_j$ (resp., ${}_i{\bf x}_j\in B^{Kh}_j$) to the $(m+1)\times(m+1)$ matrix whose only nonzero matrix entry is a $1$ (resp., an ${\bf x}$), located in row number $i$ and column number $j$ (where we assume that rows and columns are numbered from $0$ to $m$).
\end{remark}

We close our discussion of $B^{Kh}$ with a technical lemma that will prove useful in our construction of the braid group action on $D_\infty\left(B^{Kh}\right)$ (in particular, in the proof of Proposition \ref{prop:MCbetak}).

\begin{lemma} \label{lem:iotapB} Let $\iota: B^{Kh} \rightarrow B$, $p: B \rightarrow B^{Kh}$, and $h:B \rightarrow B$ be the $\mathbb{F}$--linear transformations defined in the proof of Lemma \ref{lem:Bformal}.  The $A_\infty$ morphism of $B^{Kh}$--modules, $\iota_B: B^{Kh} \rightarrow B$, given by \[\left(\iota_B\right)_{(n_1|1|n_2)} := \left\{\begin{array}{cl} 
							\iota & \mbox{if $n_1=n_2=0$, and}\\
							0 & \mbox{otherwise.}\end{array}\right.\]
is a quasi-isomorphism.  Furthermore, there exists an $A_\infty$ quasi-isomorphism of $B^{Kh}$--modules, $p_B: B \rightarrow B^{Kh}$, whose first few terms are given by: 
\[\left(p_B\right)_{(n_1|1|n_2)} :=\left\{\begin{array}{cl}
							p & \mbox{if $n_1=n_2=0$,}\\
							0 & \mbox{if $n_1=1$ and $n_2=0$, and}\end{array}\right.\]
$\left(p_B\right)_{(0|1|1)}: B \otimes B^{Kh} \rightarrow B^{Kh}$ is the bilinear map satisfying 
\[\left(p_B\right)_{(0|1|1)}(a \otimes b) := \]
\begin{itemize}
\item ${}_i\Id_k \,\,\,\,\,$ if $a = (\ell|\ell+1) \in {}_iB_j$ with $i<j$,  $k\leq \ell \leq i$, and $b = {}_j\Id_k \in {}_jB^{Kh}_k$ with $j>k$, $i \geq k$,
\item ${}_i{\bf x}_k \,\,\,\,\,$ if $a = (\ell|\ell+1) \in{}_iB_j$ with $i<j$,  $k+1\leq \ell \leq i$, and $b = {}_j{\bf x}_k \in {}_jB^{Kh}_k$ with $j>k$, $i > k$, and
\item ${}_i{\bf x}_k \,\,\,\,\,$ if $a = (\ell|\ell-1|\ell) \in{}_iB_j$ with $i\leq j$,  $k+1\leq \ell \leq i$, and $b = {}_j\Id_k \in {}_jB^{Kh}_k$ with $j>k$, $i > k$.
\item $0 \,\,\,\,\,$ for all other basis elements $a \in B$, $b \in B^{Kh}$ in the proof of Lemma \ref{lem:Bformal}.
\end{itemize}

\end{lemma}

\begin{proof}
Let $m_{(n_1|1|n_2)}$ denote the structure maps for $B$ and $m^{Kh}_{(n_1|1|n_2)}$ denote those (induced by Proposition \ref{prop:InducedAinfty}) for $B^{Kh}$, both considered as $B^{Kh}$--bimodules.

Recall that the ``Transfer Theorem" \cite[Thm. 2.1]{Berglund} tells us how to extend $\iota,p$ to $A_\infty$ quasi-isomorphisms.  Explicitly, one defines \[\left(\iota_B\right)_{(0|1|0)} := \iota, \hskip 15 pt \left(p_B\right)_{(0|1|0)} := p\] and constructs higher terms of $\iota_B, p_B$ satisfying the $A_\infty$ relations for morphisms.  Since $\iota, p$ induce isomorphisms on homology, $\iota_B$ and $p_B$ will then yield $A_\infty$ quasi-isomorphisms $B \leftrightarrow B^{Kh}$.

We begin by calculating the higher terms of $\iota_B$.  But here our work is already done, since $\iota,p,$ and $h$ satisfy the assumptions of Lemma \ref{lem:FormalCond}, hence $\left(\iota_B\right)_{(n_1|1|n_2)} = 0$ for all $(n_1+1+n_2) > 1$, as desired.



We now move to the calculation of the higher terms of $p_B$.  

{\flushleft {\em Computation of $\left(p_B\right)_{(1|1|0)}$:}}

Here we note that $ph = 0$, and Im($h$) and Im($m_1^B$) are both {\em left} $B^{Kh}$ submodules, so an application of Lemma \ref{lem:FormalCond2}, implies that $p: B \rightarrow B^{Kh}$ is a {\em left} module map (and, hence, we can extend $p$ to a left $A_\infty$ morphism with no higher left $A_\infty$ terms).  In particular, $\left(p_B\right)_{(1|1|0)} := 0$, as desired.

{\flushleft {\em Computation of $\left(p_B\right)_{(0|1|1)}$:}}

Unfortunately, Im($h$) and Im($m_1^B$) are not {\em right} $B^{Kh}$ submodules, so we will have to work harder here.  The Transfer Theorem (\cite[Thm. 2.1]{Berglund}), combined with remarks in the proof of Lemma \ref{lem:FormalCond2}, tells us that \[\left(p_B\right)_{(0|1|1)} := p \circ m_{(0|1|1)}\circ (h \otimes 1).\]

We now claim that if ${}_ia_j \in {}_iB_j$ and ${}_jb_k \in {}_jB^{Kh}_k$, then $\left(p_B\right)_{(0|1|1)} ({}_ia_j\otimes {}_jb_k) =0$ unless the triple $i,j,k$ satisfies the property that $i \leq j$, $j>k$, and $i \geq k$.  We can see this by a case-by-case analysis  (see the table below, which describes $\left(p_B\right)_{(0|1|1)}$ in the various cases).  For example, if $j<k$ (first column of table) then ${}_jb_k= 0$, and if $i<k$ (first entry in second column), then $p_{(0|1|0)} := 0$.  In both cases, we then have $\left(p_B\right)_{(0|1|1)}({}_ia_j \otimes {}_jb_k) = 0$.  On the other hand, when  $i > j \geq k$ or $i=j=k$ (the remaining entries in the table except the top two in the third column), we notice that \[m_{(0|1|1)}(Im(h) \otimes {}_jb_k) \subseteq Im(h).\]  Since $ph = 0$, we have $\left(p_B\right)_{(0|1|1)}=0$ in these cases as well.

We are therefore left to compute $\left(p_B\right)_{(0|1|1)}$ when $i \leq j, j>k,$ and $i \geq k$ (the starred entries of the table).  There are three subcases.

\[
\begin{array}{c|c|c|c}
	\left(p_B\right)_{(0|1|1)}\left({}_ia_j \otimes {}_jb_k\right)& j<k & j=k & j>k\\
	\hline
	i<j & 0 & 0 & *\\
	\hline
	i=j & 0 & 0 & *\\
	\hline
	i>j & 0 & 0 & 0\\
\end{array}\]

{\flushleft {\bf Case 1: $i<j$, $j>k$, and $i=k$}}  

Here, we notice that for basis elements ${}_ia_j, {}_jb_k$, we have \[\left(p_B\right)_{(0|1|1)}({}_ia_j \otimes {}_jb_k) \neq 0\] iff ${}_ia_j = (i|i+1)$ and ${}_jb_k = {}_j\Id_k.$

In this case,
\begin{eqnarray*}
	\left(p_B\right)_{(0|1|1)}({}_ia_j \otimes {}_jb_k) &:=& p\left[(0) + \ldots + (i)\right]\\
	&=& {}_i\Id_k.
\end{eqnarray*} 

{\flushleft {\bf Case 2: $i<j$, $j>k$, and $i>k$}} 

Again, we notice that for basis elements ${}_ia_j, {}_jb_k$, we have \[\left(p_B\right)_{(0|1|1)}({}_ia_j \otimes {}_jb_k) \neq 0\] iff either
\begin{itemize}
	\item ${}_ia_j = (\ell|\ell+1)$ for $k \leq \ell \leq i$ and ${}_jb_k = {}_j\Id_k$, in which case 
		\[\left(p_B\right)_{(0|1|1)}({}_ia_j \otimes {}_jb_k) := p\left[(0)+ \ldots + (k)\right] = {}_i\Id_k.\]
	\item ${}_ia_j = (\ell|\ell+1)$ for $k+1\leq \ell \leq i$ and ${}_jb_k = {}_j{\bf x}_k$, in which case \[\left(p_B\right)_{(0|1|1)}({}_ia_j \otimes {}_jb_k) := p\left[(1|0)+ \ldots + (k+1|k)\right] = {}_i{\bf x}_k.\]

	\item ${}_ia_j = (\ell+1|\ell|\ell+1)$ for $k \leq \ell \leq i$ and ${}_jb_k = {}_j\Id_k$, in which case \[\left(p_B\right)_{(0|1|1)}({}_ia_j \otimes {}_jb_k) := p\left[(1|0)+ \ldots + (k+1|k)\right] = {}_i{\bf x}_k.\]

\end{itemize}

{\flushleft {\bf Case 3: $i=j>k$}} 

An analysis similar to the previous cases allows us to conclude that \[p_{(0|1|1)}({}_ia_j \otimes {}_jb_k)=0\] on basis elements ${}_ia_j,{}_jb_k$ except when ${}_ia_j = (\ell|\ell-1|\ell)$ for $k+1 \leq \ell \leq i$ and ${}_jb_k = {}_j\Id_k$.  In these cases, we have $p_{(0|1|1)}\left[{}_ia_j \otimes {}_jb_k\right] = {}_i{\bf x}_k.$

Armed with the above calculations, we define $p_{(0|1|1)}: {}_iB_j \otimes {}_jB^{Kh}_k \rightarrow {}_iB_k$ in the case $i \leq j, j>k, i\geq k$ to be the unique bilinear map assigning the values above to the basis elements described and $0$ to all other basis elements.  The desired conclusion follows.





\end{proof}

\subsection{A braid group action on $D_\infty(B^{Kh})$}
Khovanov-Seidel's braid group action on $D(A)$, the derived category of dg modules over the dg algebra $A$, induces a braid group action on $D_\infty(B^{Kh})$, via the following:

\begin{proposition} \label{prop:ABderequiv} There is an equivalence of triangulated categories \[D(A) \leftrightarrow D(B) \leftrightarrow D_\infty(B^{Kh}).\]
\end{proposition}

\begin{proof}
Borrowing notation from \cite[Sec. 2.4.1]{GT10030598}, let $\mathcal{D}_{\infty,\infty}(-)$ denote the category whose objects are strictly unital $A_\infty$ modules over ``$-$" and whose morphisms are $A_\infty$ homotopy classes of $A_\infty$ morphisms. Since $B^{Kh}$ is a strictly unital minimal $A_\infty$ algebra, \cite[Cor. 3.3.1.3]{kenji} and \cite[Prop. 3.3.1.8]{kenji} imply the equivalence of $\mathcal{D}_{\infty,\infty}(B^{Kh})$ and $D_\infty(B^{Kh})$, and the argument given in \cite[Proposition 2.4.1]{GT10030598} implies the equivalence of $D(B)$ (there denoted $\mathcal{D}_{H,qi}(B)$) and $\mathcal{D}_{\infty,\infty}(B)$.

Moreover, $\mathcal{D}_{\infty,\infty}(B)$ and $\mathcal{D}_{\infty,\infty}(B^{Kh})$ are equivalent triangulated categories, since \cite[Prop. 2.4.10]{GT10030598} tells us that an $A_\infty$ quasi-isomorphism $\phi\colon B \rightarrow B^{Kh}$ (whose existence is guaranteed by Lemma \ref{lem:Bformal}) induces two mutually quasi-inverse functors $\mbox{Induct}_{\phi}\colon D_\infty(B)\rightarrow D_\infty(B^{Kh})$ and $\mbox{Rest}_{\phi}\colon D_\infty(B^{Kh}))\rightarrow D_\infty(B)$.\footnote{Note that although \cite[Prop. 2.4.10]{GT10030598} is formulated for categories of $A_{\infty}$ right modules, similar statements also hold for categories of $A_{\infty}$ left modules and $A_{\infty}$ bimodules; see \cite{GT10030598} for details.} 


To see that $D(A) \leftrightarrow D(B)$, we will show that the functors $\mathcal{F}: D(A) \rightarrow D(B)$ and $\mathcal{G}:D(B) \rightarrow D(A)$ given by
\begin{eqnarray*}
	\mathcal{F}(M) &:=& Q^*\otimes_A M=\Hom_A(Q,M)\\
	\mathcal{G}(N) &:=&
Q\otimes_B N
\end{eqnarray*}
where $Q:=\bigoplus_{i=0}^m Q_i$ and $Q^*:=\Hom_A(Q,A)=\bigoplus_{i=0}^m{}_iQ$ are well-defined mutually inverse equivalences of triangulated categories.

Since each ${}_iQ\subset Q^*$ is a complex of projective right modules over $A$, the functor $Q^*\otimes_A-$ is exact, so $\mathcal{F}$ is clearly well-defined.
To prove that $\mathcal{G}$ is also well-defined, we will show that the right dg $B$-module $\Hom_A(P_i,Q)\subset Q=\Hom_A(A,Q)=\bigoplus_{i=0}^m\Hom_A(P_i,Q)$ is homotopy equivalent to a semi-free dg $B$-module, and so tensoring with this dg $B$-module is exact.

Let $MC({}_{i}\Id_{i-1})$ denote the mapping cone of the chain map
${}_i\Id_{i-1}\colon Q_i\rightarrow Q_{i-1}$ defined by ${}_i\Id_{i-1}:=(0)+\ldots+(i-1)\in\Hom_A(Q_i,Q_{i-1})$. There is an $A$--linear chain map $\iota\colon P_i\rightarrow MC({}_i\Id_{i-1})$ given by the inclusion of $P_i$ into $Q_i$, and an $A$--linear chain map $p\colon MC({}_i\Id_{i-1})\rightarrow P_i$ given by
\[p(\phi) := \left\{\begin{array}{ll}
				\phi&\mbox{if $\phi\in P_i\subset Q_i$, and}\\
				-\phi(i-1|i)&\mbox{if $\phi\in P_{i-1}\subset Q_{i-1}$, and}\\
				0 & \mbox{otherwise.}
			\end{array}\right.\]
We leave it to the reader to verify that
\[p\iota=\mbox{Id} \hskip 20pt \mbox{and} \hskip 20pt \iota p = \mbox{Id} + \partial h + h\partial,\]
where $\partial$ is the differential in $MC({}_i\Id_{i-1})$ and $h\colon MC({}_i\Id_{i-1})\rightarrow MC({}_i\Id_{i-1})$ is the $A$--linear map  $h:={}_{i-1}1_i\colon Q_{i-1}\rightarrow Q_i$ defined by ${}_{i-1}\Id_i:=(0)+\ldots+(i-1)\in\Hom_A(Q_{i-1},Q_i)$.

Thus $P_i$ is homotopy equivalent to the mapping cone of the chain map ${}_i\Id_{i-1}\colon Q_i\rightarrow Q_{i-1}$, and consequently, $\Hom_A(P_i,Q)$ is homotopy equivalent to the mapping cone of the induced chain map  ${}_{i-1}f_{i}\colon{}_{i-1}B\rightarrow{}_iB$, where ${}_iB:=\Hom_A(Q_i,Q)=({}_i\Id_i)B$. Since $MC({}_{i-1}f_{i})$ is semi-free (because ${}_{i-1}B$ and ${}_{i}B$ are semi-free),  the functor $(\Hom_A(P_i,Q)\otimes_B-)\cong (MC({}_{i-1}f_{i})\otimes_B-)$ is exact, as desired.

It remains to show that the functors $\mathcal{F}$ and $\mathcal{G}$ are inverses of each other. Clearly, the composition $\mathcal{F}\circ\mathcal{G}$ is isomorphic to the identity functor of $D(B)$ because $Q^*\otimes_A Q\cong\Hom_A(Q,Q)=B$ by Lemma~\ref{lem:Homcpx}. To show that the composition $\mathcal{G}\circ\mathcal{F}$ is isomorphic to the identity functor of $D(A)$, we will show that the map
$$
\psi \colon Q\otimes_B Q^*\longrightarrow A
$$
defined by $\psi(q\otimes f):=f(q)\in A$ for $f\in Q^*$ and $q\in Q$ is an isomorphism of dg bimodules. We first note that the differential in $Q\otimes_B Q^*$ is trivial because the differential in $Q$ (resp., $Q^*$) is given by right (resp., left) multiplication with the element \[b:=\sum_{i=1}^m \big((0|1)+\ldots+(i-1|i)\big)\in\bigoplus_{i=1}^m\Hom_A(Q_i,Q_i)\subset B;\] and so the differential in $Q\otimes_B Q^*$ is equal to $b\otimes\mbox{Id}+\mbox{Id}\otimes b=2(b\otimes\mbox{Id})=0$. Since the differential in $A$ is trivial as well, it thus suffices to show that $\psi$ is a homotopy equivalence.

However, we have already seen that $Q$ is homotopy equivalent to a sum of complexes of the form ${}_iB\rightarrow{}_{i-1}B$ where ${}_iB=\Hom_A(Q_i,Q)$, and an analogous argument shows that $Q^*$ is homotopy equivalent to a sum of complexes of the form $B_{i-1}\rightarrow B_i$ where $B_i:=\Hom_A(Q,Q_i)$, and $B$ is homotopy equivalent to a sum of complexes of the form ${}_iB_{j-1}\rightarrow({}_{i-1}B_{j-1}\oplus{}_iB_j)\rightarrow {}_{i-1}B_j$ where ${}_iB_j:=\Hom_A(Q_i,Q_j)$. Moreover, one can check that under these various homotopy equivalences, the map $\psi$ corresponds to the canonical map from $({}_iB\rightarrow{}_{i-1}B)\otimes_B(B_{j-1}\rightarrow B_j)$ to ${}_iB_{j-1}\rightarrow({}_{i-1}B_{j-1}\oplus{}_iB_j)\rightarrow {}_{i-1}B_j$, and now the fact that $\psi$ is a homotopy equivalence follows from the identities \[{}_iB\otimes_B B_j=({}_i\Id_i)B\otimes_BB({}_j\Id_j)=({}_i\Id_i)B({}_j\Id_j)={}_iB_j.\]

\end{proof}





To understand the braid group action on $D_\infty(B^{Kh})$, recall (see \cite[Sec. 2d]{MR1862802}) that Khovanov-Seidel associate 
\begin{itemize} 
	\item to the elementary Artin generator $\sigma_k^+$ the dg $A$--bimodule \[\cM_{\sigma_k^+} := \xymatrix{0  \ar[r] &  P_k \otimes {}_k{P}  \ar[r]^-{\beta_k} & A \ar[r] & 0},\] where $\beta_k$ is the $A$--bimodule map specified by $\beta_k((k)\otimes(k)) = (k)$, and 
	\item to the elementary Artin generator $\sigma_k^{-}$ the dg $A$--bimodule \[\cM_{\sigma_k^-} := \xymatrix{0  \ar[r] &  A\ar[r]^-{\gamma_k} & P_k \otimes {}_k{P}\{-1\} \ar[r] & 0},\] where \[\gamma_k(1) = (k-1|k) \otimes (k|k-1) + (k+1|k)\otimes (k|k+1) + (k) \otimes (k|k-1|k) + (k|k-1|k) \otimes (k).\]  Here, ``$1$" denotes the identity element $1 = \sum_{i=0}^m (i)$.
\end{itemize}

We can therefore understand the induced braid group action on $D_\infty(B^{Kh})$ by understanding the images of $\cM_{\sigma_k^\pm}$ under the derived equivalence $D_\infty(A) \rightarrow D_\infty(B) \rightarrow D_\infty(B^{Kh})$.

Accordingly, we denote by $\widetilde{\cM}_{\sigma_k^+}$ (resp., $\widetilde{\cM}_{\sigma_k^-}$) the mapping cone
\[\xymatrix{0 \ar[r] & \Hom_A\left(\bigoplus_{i = 0}^m Q_i, P_k\right) \otimes \Hom_A\left(P_k,\bigoplus_{j=0}^m Q_j\right) \ar[r]^{\hskip 90pt\widetilde{\beta}_k} & B \ar[r] & 0}\]
(resp.,
\[\xymatrix{0 \ar[r] & B \ar[r]^{\hskip -100pt \widetilde{\gamma}_k} & \Hom_A\left(\bigoplus_{i = 0}^m Q_i, P_k\right) \otimes \Hom_A\left(P_k,\bigoplus_{j=0}^m Q_j\right) \{-1\}\ar[r] & 0}),\] considered as a $B^{Kh}$-$B^{Kh}$ dg bimodule.

After an application of  
Lemma \ref{lem:Homcpx}:
\[\Hom_A\left(\bigoplus_{i = 0}^m Q_i, P_k\right) \otimes \Hom_A\left(P_k,\bigoplus_{j=0}^m Q_j\right) = \left(\bigoplus_{i=0}^m {}_iQ \right) \otimes_A P_k \otimes \,{}_kP \otimes_A \left(\bigoplus_{i=0}^m Q_j\right),\] the induced maps $\widetilde{\beta}_k$, $\widetilde{\gamma}_k$ can be described as $\widetilde{\beta}_k = \mbox{Id} \otimes \beta_k \otimes \mbox{Id}$ and $\widetilde{\gamma}_k = \mbox{Id} \otimes \gamma_k \otimes \mbox{Id}$.

To further streamline notation, we set \[\widetilde{P}_k := \Hom_A\left(\bigoplus_{i = 0}^m Q_i, P_k\right)\] and
\[{}_k\widetilde{P} := \Hom_A\left(P_k,\bigoplus_{j=0}^m Q_j\right).\]

We will also find it convenient to replace the mapping cones $\widetilde{\cM}_{\sigma_k^\pm}$ with simpler, quasi-isomorphic, mapping cones.  We do this by replacing each bimodule $B$ and $\widetilde{P}_k \otimes {}_k\widetilde{P}$ by its homology and the maps $\widetilde{\beta}_k, \widetilde{\gamma}_k$ by the induced maps on homology.  

We already understand the structure of $B^{Kh} = H_*(B)$ (Lemma \ref{lem:Bformal}).  The homology of $\widetilde{P}_k$ (resp., ${}_k\widetilde{P}$) is described by:

\begin{lemma} \label{lem:PkHom} $\widetilde{P}_k$ (resp., ${}_k\widetilde{P}$) is formal as a left (resp., right) $B^{Kh}$ module.  

Furthermore, $P_k^{Kh} := H_*\left(\widetilde{P}_k\right)$ and ${}_kP^{Kh} := H_*\left({}_k\widetilde{P}\right)$ have the following explicit descriptions. \[P_k^{Kh} = \mbox{Span}_\F\langle {\bf u}^*, {\bf v}^* \rangle, \hskip 20pt {}_kP^{Kh} = \mbox{Span}_\F\langle {\bf u}, {\bf v} \rangle\]
where the bigradings of ${\bf u}^*, {\bf v}^*, {\bf u}, {\bf v}$ are given by:
\[\mbox{gr}({\bf u}^*) = (0,1), \hskip 20pt \mbox{gr}({\bf v}^*) = (1,0), \hskip 20pt \mbox{gr}({\bf u}) = (0,0), \hskip 20pt \mbox{gr}({\bf v}) = (-1,1),\]
and left multiplication by a generator $\theta \in B^{Kh}$ on $P_k^{Kh}$ is given by:
\[\theta \cdot {\bf u}^* = \left\{\begin{array}{cl}
					{\bf u}^* & \mbox{if $\theta = {}_k\Id_k$,}\\
					0 & \mbox{otherwise.}\end{array}\right.
\hskip 20pt \theta \cdot {\bf v}^* = \left\{\begin{array}{cl}
						{\bf v}^* & \mbox{if $\theta = {}_{k-1}\Id_{k-1}$}\\
						{\bf u}^* & \mbox{if $\theta = {}_k{\bf x}_{k-1}$,}\\
						0 & \mbox{otherwise.}\end{array}\right.\]
and right multiplication by a generator $\theta \in B^{Kh}$ on ${}_kP^{Kh}$ is given by:
\[{\bf u}\cdot \theta = \left\{\begin{array}{cl}
						{\bf u} & \mbox{if $\theta = {}_k\Id_k$}\\
						{\bf v} & \mbox{if $\theta = {}_k{\bf x}_{k-1}$}\\
						0 & \mbox{otherwise.}\end{array}\right.
\hskip 20pt {\bf v}\cdot \theta = \left\{\begin{array}{cl}
						{\bf v} & \mbox{if $\theta = {}_{k-1}\Id_{k-1}$}\\
						0 & \mbox{otherwise.}\end{array}\right.\]

\end{lemma}

\begin{proof} 

By Lemma \ref{lem:Homcpx}, $\Hom_A(Q_i, P_k)$ is given by the complex ${}_iQ \otimes_A P_k$ and $\Hom_A(P_k,Q_i)$ by ${}_kP \otimes_A Q_i$, where
\[{}_iQ := \xymatrix{{}_0P & \ar[l]_{(0|1)\cdot} {}_1P & \ar[l]_{(1|2)\cdot} \ldots & \ar[l]_{(i-1|i)\cdot}  {}_iP}\]

This implies that $\widetilde{P}_k, {}_k\widetilde{P}$ are given by: 
\begin{eqnarray*}
\widetilde{P}_k &:=& \bigoplus_{i=0}^m \xymatrix{{}_0P_k & \ar[l]_{(0|1)\cdot} {}_1P_k & \ar[l]_{(1|2)\cdot} \ldots & \ar[l]_{(i-1|i)\cdot} {}_iP_k}\\
{}_k\widetilde{P} &:=& \bigoplus_{i=0}^m \xymatrix{{}_kP_0 \ar[r]^{\cdot (0|1)}& {}_kP_1 \ar[r]^{\cdot (1|2)}& \ldots \ar[r]^{\cdot (i-1|i)} & {}_kP_i}
\end{eqnarray*}

We see from above that ${}_iQ \otimes_A P_k$ is:

\begin{itemize}
	\item $0$ when $i < k-1$,
	\item rank one, generated by $(k-1|k) \in {}_{k-1}P_k$, with $0$ differential, when $i = k-1$,
	\item a direct sum of Span$\langle(k|k-1|k)\rangle \subset {}_kP_k$ and the acyclic subcomplex \[(k-1|k) \leftarrow (k)\subset \left\{{}_{k-1}P_k \leftarrow {}_kP_k\right\}\] when $i=k$, and
	\item a direct sum of the two acyclic subcomplexes 
		\[(k-1|k) \leftarrow (k)\subset \left\{{}_{k-1}P_k \leftarrow {}_kP_k\right\} \mbox{ and } (k|k-1|k) \leftarrow (k+1|k)\subset \left\{{}_{k}P_k \leftarrow {}_{k+1}P_k\right\}\] when $i > k$.
\end{itemize}

To show formality of $\widetilde{P}_k$, we use Lemma \ref{lem:FormalCond2} to show that all induced multiplications \[m_{(n-1|1|0)}: \left(B^{Kh}\right)^{\otimes n-1} \otimes H_*(\Hom_A(Q_i,P_k)) \rightarrow H_*(\Hom_A(Q_j,P_k))\] vanish for $n>2$.

When $i \leq k-1$, $\Hom_A(Q_i,P_k)$ has trivial differential, so the maps $\iota_i,p_i,h_i$ are clear.
In the case $i \geq k$, we define:
\begin{eqnarray*}
	\iota_i&:& H_*\left(\Hom_A(Q_i,P_k)\right) \rightarrow \Hom_A(Q_i,P_k),\\ 
	p_i&:& \Hom_A(Q_i,P_k) \rightarrow H_*\left(\Hom_A(Q_i,P_k)\right), \mbox{and}\\
	h_i&:& \Hom_A(Q_i,P_k) \rightarrow \Hom_A(Q_i,P_k)[-1]
\end{eqnarray*} 
as follows.

Let $\theta$ denote any generator of $\Hom_A(Q_i, P_k)$, let ${\bf{u}^*}$ denote the lone generator of $H_*(\Hom_A(Q_k,P_k))$, and let $\partial$ denote the differential on the complex $\Hom_A(Q_i,P_k)$.  Note that $H_*(\Hom_A(Q_i,P_k)) = 0$ for $i>k$.  Then we define $\iota_i, p_i, h_i$ to be the \F--linear extensions of:
\begin{eqnarray*}
	\iota_k({\bf u}^*) &:=& (k|k-1|k)\\
	\iota_{i > k} &:=& 0
\end{eqnarray*}
\begin{eqnarray*}
	p_i(\theta) &:=& \left\{\begin{array}{ll}
					 \bf{u}^* & \mbox{if $i=k$ and ${\theta}= (k|k-1|k)$}\\
					 0 & \mbox{otherwise}\end{array}\right.
\end{eqnarray*}
and
\[h_i({\bf \theta}) := \left\{\begin{array}{ll}
					\partial^{-1}({\theta}) & \mbox{if } {\theta} \in \mbox{Im}(\partial),\\
					0 & \mbox{otherwise.}
		 			\end{array}\right.\]

In the above, $\partial^{-1}(\theta)$ is defined to be the (unique) basis element $\theta'$ satisfying $\partial(\theta') = \theta$.


It is now straightforward to verify that 
\begin{enumerate}
	\item $p_i h_i = 0$ for all $i$, and
	\item Im($h_i$) and Im($\partial$) are left $B^{Kh}$-submodules.
\end{enumerate}

Therefore $\widetilde{P}_k$ is formal by Lemma \ref{lem:FormalCond2}.  

To see that ${}_k\widetilde{P}$ is also formal, we perform a very similar computation, observing that ${}_k\widetilde{P}$ satisfies the assumptions of Lemma \ref{lem:FormalCond} as a {\em right} $B^{Kh}$--module, hence is formal.

Now, we simply note that $H_*(\widetilde{P}_k)$ is rank $2$, generated by 
\begin{itemize}
	\item ${\bf u}^* := p_k(k|k-1|k) \in {}_kP_k \subset \Hom_A\left(Q_k,P_k\right)$ and
	\item ${\bf v}^* := p_{k-1}(k-1|k) \in {}_{k-1}P_k \subset \Hom_A\left(Q_{k-1},P_k\right)$,
\end{itemize}
as is $H_*({}_k\widetilde{P})$, generated by
\begin{itemize}
	\item ${\bf u} := p_k(k) \in {}_kP_k \subset \Hom_A\left(P_k, Q_k\right)$ and
	\item ${\bf v} := p_{k-1}(k|k-1) \in {}_kP_{k-1} \subset \Hom_A\left(P_k,Q_{k-1}\right)$.
\end{itemize}

Recalling (see the proof of Lemma \ref{lem:Bformal}) that the generators ${}_i\Id_j$ (for $i\geq j$) and ${}_i{\bf x}_j$ (for $i>j$) of $B^{Kh}$ are represented by $(0) + \ldots + (j)$ and $(1|0) + \ldots + (j+1|j)$, we see that the multiplication is also as claimed.
\end{proof}

We now have the proposed model 
\[MC\left(\xymatrix{P_k^{Kh} \otimes {}_kP^{Kh} \ar[r]^-{{\beta}^{Kh}_k} & B^{Kh}}\right)\] for $\cM^{Kh}_{\sigma_k^+}$ and the model
 \[MC\left(\xymatrix{B^{Kh} \ar[r]^-{{\gamma}^{Kh}_k} & P_k^{Kh} \otimes {}_kP^{Kh} \{-1\}}\right)\] for $\cM^{Kh}_{\sigma_k^-}$, where $\beta_k^{Kh}$ and $\gamma_k^{Kh}$ are the $A_\infty$ morphisms on homology induced by $\widetilde{\beta}_k$ and $\widetilde{\gamma}_k$.
 
To understand the induced maps on homology, we must explicitly understand the quasi-isomorphisms $B \leftrightarrow B^{Kh}$ and $\widetilde{P}_k \otimes {}_k\widetilde{P} \leftrightarrow P_k^{Kh} \otimes {}_kP^{Kh}$.

Explicitly, if $\iota_{P} \otimes \iota_P':  P_k^{Kh} \otimes {}_kP^{Kh} \rightarrow \widetilde{P}_k \otimes {}_k\widetilde{P}$ and $p_B: B \rightarrow B^{Kh}$ are $A_\infty$ quasi-isomorphisms, then the induced $A_\infty$ morphism on homology is given by:
 \[{\beta}^{Kh}_k=p_B\circ \widetilde{\beta}_k \circ (\iota_{P} \otimes \iota_P'): P_k^{Kh} \otimes {}_kP^{Kh} \rightarrow B^{Kh}.\]

Furthermore, (cf. \cite[Cor. 3.16]{SeidelBook}), the mapping cones satisfy:
 \begin{eqnarray*}
 	\left(\xymatrix{0 \ar[r] & P_k^{Kh} \otimes {}_kP^{Kh} \ar[rrr]^-{{\beta}^{Kh}_k=p_B\circ \widetilde{\beta}_k \circ (\iota_{P} \otimes \iota_P')} & & & B^{Kh} \ar[r] & 0}\right) = \\
	\left(\xymatrix{0 \ar[r] & \widetilde{P}_k \otimes {}_k\widetilde{P} \ar[r]^-{\widetilde{\beta}_k} & B \ar[r] & 0}\right)
\end{eqnarray*} 
as elements of $D_\infty(B^{Kh})$.

Similarly, if $\iota_B: B^{Kh} \rightarrow B$ and $p_{P}:  \widetilde{P}_k \otimes {}_k\widetilde{P} \rightarrow P_k^{Kh} \otimes {}_kP^{Kh}$ are $A_\infty$ quasi-isomorphisms, then:

\begin{eqnarray*}
 	\left(\xymatrix{0 \ar[r] & B^{Kh} \ar[rrr]^-{\gamma_k^{Kh} = (p_P \otimes p_P') \circ \widetilde{\beta}_k \circ \iota_{B}} & & & P_k^{Kh} \otimes {}_kP^{Kh}\{-1\}  \ar[r] & 0}\right) = \\
	\left(\xymatrix{0 \ar[r] &  B \ar[r]^-{\widetilde{\gamma}_k} &\widetilde{P}_k \otimes {}_k\widetilde{P} \{-1\} \ar[r] & 0}\right)
\end{eqnarray*} 
as elements of $D_\infty(B^{Kh})$.

\begin{proposition} \label{prop:MCgammak} The image of $\mathcal{M}_{\sigma_k^-} \in D_\infty(A)$ under the derived equivalence $D_\infty(A) \rightarrow D_\infty\left(B^{Kh}\right)$ is $MC\left(\gamma_k^{Kh}\right)$, where
\[\gamma_k^{Kh}: B^{Kh} \rightarrow P_k^{Kh} \otimes {}_kP^{Kh}\{-1\}\] is the \F--linear $B^{Kh}$--bimodule map (i.e., strict $A_\infty$ morphism) determined by
\[{}_i\Id_i \mapsto \left\{\begin{array}{cl}
				{\bf v}^* \otimes {\bf v} & \mbox{when } i = k-1,\\
				{\bf u}^* \otimes {\bf u} & \mbox{when } i=k, \mbox{ and}\\
				0 & \mbox{otherwise.}\end{array}\right.\]
Accordingly, we define $\mathcal{M}^{Kh}_{\sigma_k^-} := MC\left(\gamma_k^{Kh}\right)$
\end{proposition}

\begin{proof}
We must compute the terms of the induced $A_\infty$ morphism $\gamma_k^{Kh} := (p_P \otimes p_P') \circ \widetilde{\beta}_k \circ \iota_B$, as described above.

We begin by noting that the $(n_1|1|n_2)$ map of the $A_\infty$ morphism $\gamma_k^{Kh}$, i.e., the map
\[\left(\gamma_k^{Kh}\right)_{(n_1|1|n_2)}: \left(B^{Kh}\right)^{\otimes n_1} \otimes B^{Kh} \otimes \left(B^{Kh}\right)^{\otimes n_2} \rightarrow \left(P_k^{Kh} \otimes {}_kP^{Kh}\right)\{-1\}\] is degree $(-(n_1+n_2),0)$ with respect to the bigrading.  This follows from the $A_\infty$ relations for morphisms, combined with Lemma \ref{lem:degree}.

An examination of the bigradings of elements of $B^{Kh}$ and $P_k^{Kh} \otimes {}_kP^{Kh}$ then immediately implies that $\left(\gamma_k^{Kh}\right)_{(n_1|1|n_2)} = 0$ unless $n_1=n_2=0$, so $\gamma_k^{Kh}$ is a strict $A_\infty$ isomorphism, as desired.  A quick way to see this is to notice that the sum of the two gradings associated to each element in $B^{Kh}$ and $\left(P_k^{Kh} \otimes {}_kP^{Kh}\right)\{-1\}$ is $0$, and $\left(\gamma_k\right)_{(n_1|1|n_2)}$ is degree $-(n_1+n_2)$ on this sum.

It is now easy to verify that \[\left(\gamma_k^{Kh}\right) := \left(\gamma_k^{Kh}\right)_{(0|1|0)} = \left(p_P \otimes p_P'\right)_{(0|1|0)} \circ \widetilde{\gamma}_k \circ \left(\iota_B\right)_{(0|1|0)}\] is as described.  In particular, $\gamma_k^{Kh}$ is determined by its behavior on the $(m+1)$ idempotents ${}_i\Id_i \in B^{Kh}$, since it is a $B^{Kh}$--bimodule map.

For example:
\begin{eqnarray*}
	\gamma_k^{Kh}\left({}_{k-1}\Id_{k-1}\right) &:=& \left(p_P\right)_{(0|1|0)} \circ \widetilde{\gamma}_k \circ \left(\iota_B\right)_{(0|1|0)}\left({}_{k-1}\Id_{k-1}\right)\\
	&=& \left(p_P\right)_{(0|1|0)}\circ\widetilde{\gamma}_k\left[(0) + \ldots + (k-1)\right]\\
	&=& \left(p_P\right)_{(0|1|0)}\left[(k-1|k) \otimes (k|k-1)\right]\\
	&=& {\bf v}^* \otimes {\bf v}
\end{eqnarray*}

We leave the remaining similarly straightforward computations to the reader.
\end{proof}

\begin{proposition} \label{prop:MCbetak} The image of $\mathcal{M}_{\sigma_k^+} \in D_\infty(A)$ under the derived equivalence $D_\infty(A) \rightarrow D_\infty\left(B^{Kh}\right)$ is $MC\left(\beta_k^{Kh}\right)$, where the terms of the $A_\infty$ morphism $\beta_k^{Kh}$ are given as follows.

\[\left(\beta_k^{Kh}\right)_{(n_1|1|n_2)}: \left(B^{Kh}\right)^{\otimes n_1} \otimes \left(P_k^{Kh} \otimes {}_kP^{Kh}\right) \otimes \left(B^{Kh}\right)^{\otimes n_2} \rightarrow B^{Kh}\] is identically zero unless $n_1+n_2 =1$.

When $n_1=1, n_2=0$: \[\betatwoleft: B^{Kh} \otimes \left(P_k^{Kh} \otimes {}_kP^{Kh}\right) \rightarrow B^{Kh}\] is the trilinear map satisfying:
\[\betatwoleft: \left\{\begin{array}{cc}
	\left[{}_{i}\Id_k \otimes \left({\bf u}^* \otimes {\bf u}\right)\right] \mapsto {}_{i} {\bf x}_k & (i \geq k+1)\\
	\left[{}_i\Id_{k-1} \otimes \left({\bf v}^* \otimes {\bf u}\right)\right] \mapsto {}_i\Id_k & (i \geq k)\\
	\left[{}_i{\bf x}_{k-1} \otimes \left({\bf v}^* \otimes {\bf u} \right)\right] \mapsto {}_i{\bf x}_k & (i \geq k+1)\\
	\left[{}_i\Id_{k-1} \otimes \left({\bf v}^* \otimes {\bf v}\right)\right] \mapsto {}_i{\bf x}_{k-1} & (i \geq k)
\end{array}\right.\] and $\betatwoleft(b \otimes \theta) = 0$ for all other basis elements $b \in B^{Kh}, \theta \in \left(P_k^{Kh} \otimes {}_kP^{Kh}\right)$.

When $n_1=0, n_2=1$: \[\betatworight: \left(P_k^{Kh} \otimes {}_kP^{Kh}\right) \otimes B^{Kh} \rightarrow B^{Kh}\] is the trilinear map satisfying:
\[\betatworight: \left\{\begin{array}{cc}
	\left[\left({\bf u}^* \otimes {\bf u}\right) \otimes {}_k\Id_j\right] \mapsto {}_k{\bf x}_j & (j \leq k-1)\\
	\left[\left({\bf v}^* \otimes {\bf u}\right) \otimes {}_k\Id_{j}\right] \mapsto {}_{k-1}\Id_{j} & (j \leq k-1)\\
	\left[\left({\bf v}^* \otimes {\bf u}\right) \otimes {}_{k}{\bf x}_j\right] \mapsto {}_{k-1}{\bf x}_j & (j \leq k-2)\\
	\left[\left({\bf v}^* \otimes {\bf v}\right) \otimes {}_{k-1}\Id_{j}\right] \mapsto {}_{k-1}{\bf x}_{j} & (j \leq k-2)
\end{array}\right.\]
and $\betatworight(\theta \otimes b) = 0$ for all other basis elements $b \in B^{Kh}, \theta \in \left(P_k^{Kh} \otimes {}_kP^{Kh}\right)$.

Accordingly, we define $\mathcal{M}^{Kh}_{\sigma_k^+} := MC\left(\beta_k^{Kh}\right)$
\end{proposition}

\begin{proof}
As in the proof of Proposition \ref{prop:MCgammak}, the $(n_1|1|n_2)$ map of the $A_\infty$ morphism $\beta_k^{Kh}$ is degree $(-(n_1+n_2),0)$ with respect to the bigrading.

In this case, however, we see that the sum of the two gradings for each element in $P_k^{Kh} \otimes {}_kP^{Kh}$ is $1$, while the sum of the two gradings associated to each element in $B^{Kh}$ is $0$.  Since $\left(\beta_k^{Kh}\right)_{(n_1|1|n_2)}$ is degree $-(n_1+n_2)$ on this sum, we conclude that $\left(\beta_k^{Kh}\right)_{(n_1|1|n_2)} = 0$ unless $-(n_1+n_2) = -1$, as claimed.

To calculate $\left(\beta_k^{Kh}\right)_{(n_1|1|n_2)}$ in the relevant cases ($n_1=1, n_2=0$) and ($n_1=0, n_2=1$), we recall that $\beta_k^{Kh}:P_k^{Kh} \otimes {}_kP^{Kh} \rightarrow B^{Kh}$ is given by the composition
\[\xymatrix{P_k^{Kh} \otimes {}_kP^{Kh} \ar[r]^-{\iota \otimes \iota'} & \widetilde{P}_k \otimes {}_k\widetilde{P} \ar[r]^-{\widetilde{\beta}_k} & B \ar[r]^-{p} & B^{Kh}}.\]  

{\flushleft {\bf Calculation of \betatwoleft:}}

Since $\widetilde{\beta}_k$ is, by definition, a strict $A_\infty$ morphism, we see that
\[\betatwoleft := p_{(1|1|0)}\circ\widetilde{\beta}_k\circ \left(\iota_{(0|1|0)} \otimes \iota'_{(0|1|0)}\right) + p_{(0|1|0)} \circ \widetilde{\beta}_k\circ\left(\iota_{(1|1|0)}\otimes\iota'_{(0|1|0)}\right).\]

Furthermore, we showed during the proof of Lemma \ref{lem:iotapB} that $p_{(1|1|0)}:= 0$, so the first term above also vanishes, leaving:

\[\betatwoleft := p_{(0|1|0)} \circ \widetilde{\beta}_k\circ\left(\iota_{(1|1|0)}\otimes\iota'_{(0|1|0)}\right).\]

Another application of the Transfer Theorem \cite[Thm. 2.1]{Berglund} tells us that on basis elements $b \in B^{Kh}$ and $\theta \in P_k^{Kh}$, we have
\[\iota_{(1|1|0)}\left[b \otimes \theta \right] = \left\{\begin{array}{cl}
				(k+1|k) \in \mbox{Hom}_A(Q_i,P_k) & \mbox{when $b = {}_i\Id_k$, $\theta = {\bf u}^*$, and $i\geq k+1$,}\\
				(k) \in \mbox{Hom}_A(Q_i,P_k) & \mbox{when $b={}_i\Id_k$, $\theta = {\bf v}^*$, and $i\geq k$,}\\
				(k+1|k) \in \mbox{Hom}_A(Q_i,P_k) & \mbox{when $b = {}_i{\bf x}_k$, $\theta = {\bf v}^*$, and $i\geq k+1$, and}\\
				0 & \mbox{otherwise.}
				\end{array}\right.\]
				
Composing the above with $p_{(0|1|0)}\circ \widetilde{\beta}_k$ yields the desired result.  We perform this computation in one case, leaving the small number of remaining (similarly straightforward) computations to the reader.  Assume $i\geq k+1$.  Then:
\begin{eqnarray*}
	\betatwoleft({}_i\Id_k \otimes \left({\bf u}^* \otimes {\bf u}\right)) &:=& p_{(0|1|0)} \circ \widetilde{\beta}_k \circ \left[\iota_{(1|1|0)}({}_i\Id_k \otimes {\bf u}^*) \otimes \iota'_{(0|1|0)}({\bf u})\right]\\
	&=& p_{(0|1|0)} \circ \widetilde{\beta}_k \left[(k+1|k) \otimes (k)\right]\\
	&=& p_{(0|1|0)} [(k+1|k)]\\
	&=& {}_i{\bf x}_k\\
\end{eqnarray*}			

{\flushleft {\bf Calculation of \betatworight:}}
Similarly, we have:
\[\betatworight := p_{(0|1|1)}\circ\widetilde{\beta}_k\circ \left(\iota_{(0|1|0)} \otimes \iota'_{(0|1|0)}\right) + p_{(0|1|0)} \circ \widetilde{\beta}_k\circ\left(\iota_{(0|1|0)}\otimes\iota'_{(0|1|1)}\right),\] 
and an application of Lemma \ref{lem:FormalCond} (see the proof of Lemma \ref{lem:PkHom}) implies that $\iota'_{(0|1|1)} := 0$, leaving:
\[\betatworight := p_{(0|1|1)}\circ\widetilde{\beta}_k\circ \left(\iota_{(0|1|0)} \otimes \iota'_{(0|1|0)}\right).\]

Referring to Lemma \ref{lem:iotapB}, we again perform a sample computation, leaving the remaining computations to the reader.  Assume $j \leq k-1$.  Then:
\begin{eqnarray*}
	\betatworight\left(({\bf v}^* \otimes {\bf u}) \otimes {}_k\Id_j\right) &:=& p_{(0|1|1)} \circ \left(\widetilde{\beta}_k[(k-1|k)\otimes (k)] \otimes {}_k\Id_j\right)\\
	&=& p_{(0|1|1)} \left[(k-1|k) \otimes {}_k\Id_j\right]\\
	&=& {}_{k-1}\Id_j
\end{eqnarray*}

\end{proof}

Now, if we have a general braid group element $\sigma \in B_{m+1}$ that decomposes as $\sigma = \sigma_{k_1}^\pm \cdots \sigma_{k_n}^\pm$, \cite{MR1862802} associates to $\sigma \in B_{m+1}$ the dg bimodule: \[\cM_{\sigma} := \cM_{\sigma_{k_1}^\pm} \otimes_A \ldots \otimes_A \cM_{\sigma_{k_n}^\pm}\] over the algebra $A$ (or, rather, its equivalence class in $D^b(A)$). 

Considered as an element of $D_\infty(A)$, we can alternatively describe $\cM_\sigma$ in terms of an $A_\infty$ tensor product, by the following.

\begin{definition} \cite[Defn. 1]{Tangles} Given rings ${\bf A}, {\bf B}$, an ${\bf A}$-${\bf B}$ bimodule ${\bf M}$ is called {\em sweet} if it is finitely-generated and projective as a left ${\bf A}$ module and as a right ${\bf B}$ module.
\end{definition}

\begin{remark} The tensor product ${\bf N} \otimes_{\bf A} {\bf M}$ of an ${\bf A}'$-${\bf A}$ bimodule ${\bf N}$ with an ${\bf A}$-${\bf B}$ bimodule is a sweet ${\bf A}'$-${\bf B}$ bimodule.
\end{remark}

Since each $\cM_{\sigma_k^\pm}$ is a bounded complex of {\em sweet} bimodules over $A$ whose higher multiplications are all trivial, the ordinary tensor product above agrees with the $A_\infty$ tensor product in $D_\infty(A)$.  In other words, \[\cM_{\sigma} := \cM_{\sigma_{k_1}^\pm} \widetilde{\otimes}_A \ldots \widetilde{\otimes}_A \cM_{\sigma_{k_n}^\pm}.\]  Since $A_\infty$ tensor products are sent to $A_\infty$ tensor products under the derived equivalence $D_\infty(A) \leftrightarrow D_\infty(B) \leftrightarrow D_\infty(B^{Kh})$, we see that the element of $D_\infty(B^{Kh})$ associated to a general braid $\sigma = \sigma_{k_1}^\pm \cdots \sigma_{k_n}^\pm \in B_{m+1}$ is given by:
\[\cM^{Kh}_{\sigma} := \cM^{Kh}_{\sigma_{k_1}^\pm} \widetilde{\otimes}_{B^{Kh}} \,\, \ldots \widetilde{\otimes}_{B^{Kh}} \,\, \cM^{Kh}_{\sigma_{k_n}^\pm}.\]

\begin{remark}\label{rmk:gradings}
The $B^{Kh}$ modules described here (and, more generally, any $A_\infty$ module over the Hom algebra of a basis of curves) are equipped with three gradings:
\begin{enumerate}
	\item a (co)homological grading,
	\item an internal grading counting steps to the left in the path algebra, $A_m$, which corresponds to the power of $t$ under the identification of the Khovanov-Seidel construction with a categorification of the Burau representation (see \cite[Sec. 2e]{MR1862802}),
	\item a grading by path length in the path algebra, $A_m$, which corresponds to Khovanov's $j$ (quantum) grading if one identifies the Khovanov-Seidel quiver algebra $A_m$ with the algebra $A^{1,m}$ appearing in \cite{QA0610054, StroppelAlgebra}.
\end{enumerate}

The first two of these gradings constitute the bigrading described in \cite[Sec. 3d]{MR1862802} and discussed throughout this section.

For the benefit of those readers interested in the {\em trigradings} of generators of $B^{Kh}$, $P_k^{Kh}$, and ${}_kP^{Kh}$, we record them here:

\begin{itemize}
	\item $gr({}_i\Id_j) = (0,0,0)$ for ${}_i\Id_j \in {}_iB_j$ for all $i,j \in \{0, \ldots, m\}$,
	\item $gr({}_i{\bf x}_j) = (-1,1,1)$ for ${}_i{\bf x}_j \in {}_iB_j$ for all $i>j \in \{0, \ldots, m\}$,
	\item $gr({\bf v}^*) = (1,0,1)$ and $gr({\bf u}^*) = (0,1,2)$ for ${\bf v}^*, {\bf u}^* \in P_k^{Kh}$ for all $k \in \{1, \ldots, m\}$, and
	\item $gr({\bf v}) = (-1,1,1)$ and $gr({\bf u}) = (0,0,0)$ for ${\bf v}, {\bf u} \in {}_kP^{Kh}$ for all $k \in \{1, \ldots, m\}$.
\end{itemize}

\end{remark}

\subsection{$B^{Kh}$ and Fukaya categories} \label{sec:BKhFukaya}
For completeness, and to motivate the constructions in the next section, we
briefly outline a geometric interpretation of the algebra $B^{Kh}$ and the
bimodules $\mathcal{M}^{Kh}_{\sigma_i^\pm}$, in terms of the Fukaya category
of a suitable Lefschetz fibration \cite{SeidelVCM,SeidelBook,SeidelLF1}.
(Since the construction in \cite{SeidelBook} 
does not work over \Ztwo, the setup of \cite{SeidelLF1} is the most
appropriate one here.)

Namely, denote by $p$ a polynomial of degree $m+1$ whose roots are
exactly the points of $\Delta$, and consider the complex surface
$S=\{(x,y,z)\in\mathbb{C}^3\,|\,x^2+y^2=p(z)\}$. The projection to the
$z$ coordinate defines a Lefschetz fibration $\pi_S:S\to\mathbb{C}$, whose generic
fiber is an affine conic, and whose $m+1$ vanishing cycles are all isotopic
to each other. The basis of arcs $\mathcal{Q}=\{q_0,\dots,q_m\}$ of Figure 
\ref{fig:Qbasis} then determines a collection of Lefschetz thimbles
$Q_0^S,\dots,Q_m^S$ (i.e., Lagrangian disks in $S$ whose boundaries are
the vanishing cycles in the fiber $\pi_S^{-1}(-1)$). These form an exceptional
collection which generates the (directed) Fukaya category $\mathcal{F}(\pi_S)$ of the Lefschetz 
fibration $\pi_S$ \cite{SeidelBook,SeidelLF1}.

Perturbing the symplectic structure slightly, we can ensure that the
vanishing cycles (which are Hamiltonian isotopic loops in
$\pi_S^{-1}(-1)\simeq \mathbb{C}^*$) are mutually transverse and intersect
in a suitable manner (i.e., they pairwise intersect in exactly two points, 
and the intersection points are arranged in a configuration which forces the 
vanishing of higher products on Floer complexes within the ordered collection).

The Floer complexes which determine morphisms from $Q_i^S$ to $Q_j^S$ in
the directed Fukaya category
then have rank 2 whenever $i>j$, while by definition these morphism spaces
have rank 1 for $i=j$ and 0 for $i<j$~\cite{SeidelVCM}. (Note: our ordering 
convention for bases of arcs is the opposite of Seidel's.) Moreover,
an easy calculation in Floer homology then shows that $$B^S:=\bigoplus\limits_{i,j=0}^m
\mathrm{Hom}_{\mathcal{F}(\pi_S)}(Q_i^S,Q_j^S)$$ is isomorphic to $B^{Kh}$
(viewing both as $A_\infty$-algebras, in which $m_n$ happens to vanish for
$n\neq 2$). The categories of modules over $\mathcal{F}(\pi_S)$ and
$B^{Kh}$ are therefore equivalent. 

In fact, the $B^{Kh}$-module $P_k^{Kh}$ 
has a geometric counterpart via this equivalence, namely a Lagrangian sphere 
$P_k^S$ in $S$ which projects under $\pi_S$ to a line segment connecting two 
consecutive points of $\Delta$. Indeed, $P_k^S$ intersects $Q_{k-1}^S$ and 
$Q_k^S$ in one point each, and is disjoint from the other $Q_i^S$; it is then
not hard to check that $\bigoplus_i \mathrm{Hom}_{\mathcal{F}(\pi_S)}(Q_i^S,
P_k^S)\simeq P_k^{Kh}$ as an $A_\infty$-module over $B^S\simeq B^{Kh}$).
See Chapter 20 of \cite{SeidelBook} for more about the
symplectic geometry of~$S$.

Elements of the braid group $B_{m+1}$ acting on $(D_m,\Delta)$ lift to
symplectic automorphisms of $S$ preserving the fiber $\pi_S^{-1}(-1)$;
specifically, the Artin generator $\sigma_k$ lifts to the Dehn twist about
the Lagrangian sphere $P_k^S$. Denoting again by $\sigma$ the symplectic
automorphism of $S$ which corresponds to a braid $\sigma\in B_{m+1}$, we 
associate to it the $A_\infty$-bimodule 
$$\mathcal{M}^S_\sigma=
\bigoplus_{i,j=0}^m \mathrm{Hom}_{\mathcal{F}(\pi_S)}(Q_i^S,\sigma(Q_j^S))$$
over $B^S\simeq B^{Kh}$. It then follows from Seidel's long exact sequence
for Dehn twists \cite{MR1978046} that the bimodules 
$\mathcal{M}^S_{\sigma_k^{\pm}}$ and $\mathcal{M}^{Kh}_{\sigma_k^\pm}$
associated to Artin generators (or their inverses) are quasi-isomorphic.

\section{Bordered Floer algebras and bimodules} \label{sec:HFalgmod}
We now consider the analogues in bordered Floer homology of the Khovanov-Seidel bimodules described in Section \ref{sec:Khalgmod}.  We follow Lipshitz-Ozsv{\'a}th-Thurston in \cite{GT08100687, GT10030598, GT10051248} and  Zarev in \cite{GT09081106}, using a symplectic reinterpretation of their work due to the first author \cite{GT10014323}.  

\subsection{The bordered Floer algebra}
Denote by $\Sigma$ the double cover of $D_m$ branched at the
$m+1$ points of $\Delta$ (with covering map $\pi_\Sigma:\Sigma\to D_m$).
We make $\Sigma$ a {\em parametrized surface} by equipping it with two 
marked points $z_\pm$ on its boundary (the two preimages by $\pi_\Sigma$ of a point in 
$\partial D_m$) and the collection of arcs
$\mathcal{Q}_\Sigma=\{Q_0^\Sigma,\dots,Q_m^\Sigma\}$, where
$Q_k^\Sigma=\pi_\Sigma^{-1}(q_k)$. 

In the language introduced by Lipshitz, Ozsv\'ath and Thurston
\cite{GT08100687}, the parametrized surface
$(\Sigma,z_\pm,\mathcal{Q}_\Sigma)$ is described combinatorially by
a {\it (twice) pointed matched circle} (or pair of circles when $m$ is odd),
$\mathcal{Z}_\mathcal{Q}$. This consists of
 a pair of oriented intervals (the two components of
$\partial\Sigma\setminus \{z_\pm\}$), each carrying $m+1$ distinguished
points (the end points of disjoint pushoffs of the $Q_k^\Sigma$), labeled
successively in decreasing order $m,\dots,1,0$ along each interval
(according to the manner in which the end points of the 1-handles
$Q_k^\Sigma$ match up).

Recall that the $1$--moving strands algebra
$\mathcal{A}\left(\mathcal{Z}_\mathcal{Q},1\right)$,\footnote{Here we use the notation convention from \cite{GT09081106}, which differs by a shift from the one in \cite{GT08100687}.  See the note in \cite[Sec. 2.2]{GT09081106}.} which we
also denote by $B^{HF}$ for consistency with the preceding sections,
can be described as: 
\[\mathcal{A}\left(\mathcal{Z}_{\mathcal{Q}},1\right) = \bigoplus_{i,j=0}^m {}_iB^{HF}_{j},\] 
where \[{}_iB^{HF}_{j} = \mbox{Span}_{\F}\left\{\begin{array}{ll}
										0 				 & \mbox{if } i<j,\\
										{}_i\Id_i 			 & \mbox{if } i=j,\\
										{}_i\rho_j, {}_i\sigma_{j} & \mbox{if } i>j
									\end{array}\right\},\]
and the multiplication $m_2^{HF}:{}_iB^{HF}_{j} \otimes {}_jB^{HF}_{k} \rightarrow {}_iB^{HF}_{k}$
is defined by
\begin{itemize}
	\item $m_2^{HF}({}_i\Id_i \otimes a) = m_2^{HF}(a \otimes {}_j\Id_j) = a$ for all $a \in {}_iB^{HF}_{j}$, and
	\item $m_2^{HF}({}_i\rho_{j} \otimes {}_j\rho_{k}) = {}_i\rho_{k}$ and $m_2^{HF}({}_i\sigma_{j} \otimes {}_j\sigma_{k}) = {}_i\sigma_{k}$, but \\ $m_2^{HF}({}_i\rho_{j} \otimes {}_j\sigma_{k}) = m_2^{HF}({}_i\sigma_{j} \otimes {}_j\rho_{k}) = 0$.
\end{itemize}
As usual, the multiplication map $m_2^{HF}:{}_iB^{HF}_{j}\otimes {}_kB^{HF}_{\ell} \rightarrow {}_iB^{HF}_{\ell}$ is
zero unless $j=k$. We also set $m_n^{HF}=0$ for $n\neq 2$.

\begin{remark} \label{rmk:BHFMatrix} Let $\mathbb{F}\rho\oplus\mathbb{F}\sigma$ denote the $\mathbb{F}$-algebra generated by two orthogonal idempotents $\rho$ and $\sigma$, and let $1:=\rho+\sigma$ be its identity element.  As we did in the previous section for $B^{Kh}$ (Remark \ref{rmk:BKhMatrix}), we can interpret $B^{HF}$ as the algebra of all lower triangular $(m+1)\times (m+1)$ matrices over $\mathbb{F}\rho\oplus\mathbb{F}\sigma$ which have only $0$'s and $1$'s on the main diagonal:
\[
B^{HF}\cong
\left\{
\left.
\left(\begin{matrix}
d_0        & 0      & \ldots         & 0     \\
\phi_{1,0} & d_1    & \ddots         & \vdots\\
\vdots     & \ddots & \ddots         & 0     \\
\phi_{m,0} & \ldots & \phi_{m,m-1} & d_m
\end{matrix}\right)
\right |\,d_i\in\{0,1\}
\right\}\subset M_{m+1}(\mathbb{F}\rho\oplus\mathbb{F}\sigma)
\]
We identify the generator ${}_i\rho_j\in{}_iB^{HF}_j$ (resp., ${}_i\sigma_j\in{}_iB^{HF}_j$) with the $(m+1)\times(m+1)$ matrix whose only nonzero matrix entry is a $\rho$ (resp., a $\sigma$), located in row number $i$ and column number $j$; and we identify the generator ${}_i\mathbb{1}_i\in {}_iB^{HF}_i$ with the $(m+1)\times (m+1)$ matrix whose only nonzero entry is a $1$, located on the diagonal in row number $i$. (Here we assume that rows and columns are numbered from $0$ to $m$).
\end{remark}

The 1-moving strands algebra has a more 
geometric interpretation in terms of 
the arcs $Q_0^\Sigma,\dots,Q_m^\Sigma$ on the surface $\Sigma$. Namely, these
arcs (or small isotopic deformations of them) are objects of (and in fact generate) 
the ``partially wrapped'' Fukaya category
of $\Sigma$ relatively to the two marked points $z_\pm$ (see \cite{GT10032962,GT10014323}).
In this category, the morphism spaces $\hom(Q_i^\Sigma,Q_j^\Sigma)$ are
the Floer complexes generated by intersections between suitably perturbed copies 
of the arcs (namely, using the flow of a suitable Hamiltonian to ensure
transversality and push the end points so that they lie
in a specific position along the components of $\partial \Sigma\setminus \{z_\pm\}$).
In our case, $\{z_\pm\}$ is a fiber of the covering map $\pi_\Sigma$, which is in
fact a Lefschetz fibration. The partially 
wrapped Fukaya category is then equivalent to $\mathcal{F}(\pi_\Sigma)$, Seidel's
Fukaya category of the Lefschetz fibration $\pi_\Sigma$ 
(see the Remark in section 4 of \cite{GT10032962}), and the $Q_i^\Sigma$ are nothing
but the Lefschetz thimbles associated to the basis of arcs 
$\mathcal{Q}$ of Figure \ref{fig:Qbasis}. 

Note that the technical setup in \cite{GT10014323} is somewhat
different from those in \cite{SeidelVCM} and \cite{SeidelLF1}, even though
the resulting categories are equivalent and, in the case at hand, all
calculations for the thimbles $Q_i^\Sigma$ give exactly the same answer on the nose.
We use the notation $\mathcal{F}(\pi_\Sigma)$ for familiarity;
however the comparison with bordered Floer homology is simpler in the setup 
of \cite{GT10014323}, see Remark \ref{rmk:hamperturb} below.

The Floer complexes which determine morphisms from $Q_i^\Sigma$ to $Q_j^\Sigma$
have rank 2 whenever $i>j$, while these morphism spaces
have rank 1 for $i=j$ and 0 for $i<j$. In the setting of \cite{GT10014323},
this is because the image of $Q_i^\Sigma$ under the appropriate Hamiltonian 
\cite[\S 4.2]{GT10014323} intersects $Q_j^\Sigma$ transversely in 0, 1 or 2
points depending on cases; while in the directed Fukaya category of
\cite{SeidelVCM}, this is because the vanishing cycles consist of the same 
two points in the case $i>j$, and by definition in the other cases.
(As before, our ordering convention for bases of arcs is the opposite of Seidel's.) 
An easy calculation in Floer homology then shows that 
$$B^\Sigma:=\bigoplus\limits_{i,j=0}^m
\mathrm{Hom}_{\mathcal{F}(\pi_\Sigma)}(Q_i^\Sigma,Q_j^\Sigma)$$ is isomorphic to $B^{HF}$,
viewing both as $A_\infty$-algebras in which $m_n$ happens to vanish for
$n\neq 2$ (cf.\ \cite{GT10032962,GT10014323}). The categories of modules over 
$\mathcal{F}(\pi_\Sigma)$ (in any of its incarnations) and
$B^{HF}$ are therefore equivalent.

\subsection{Bordered Floer bimodules}
\begin{figure}
\begin{center}
\resizebox{4.5in}{!}{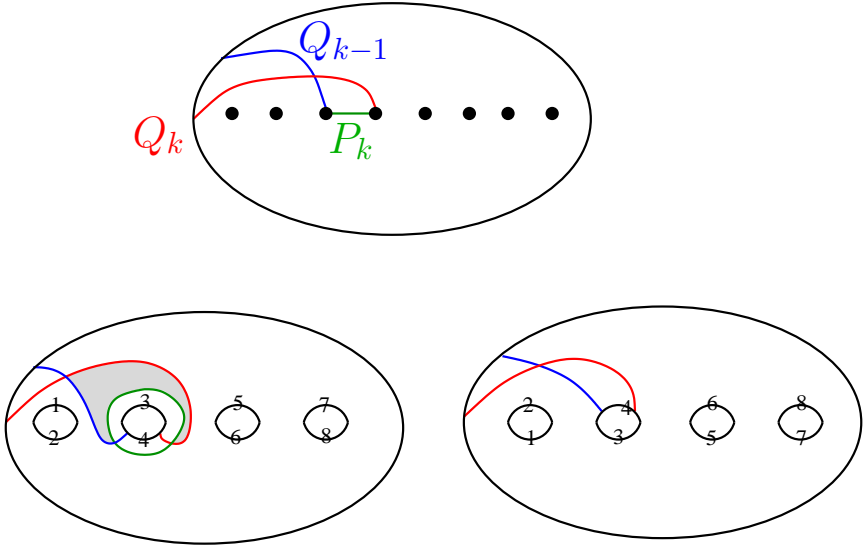}
\end{center}
\caption{The top row above shows curves $p_k, q_{k-1},$ and $q_k$ in the disk $D_m$, 
while the bottom row shows their lifts to Lagrangians in the double branched cover 
$\Sigma$ (the figures on the left and right are identified according to the numbers).  
The shaded triangle gives rise to a non-trivial multiplication map $m_{(1|1|0)}: 
\mathrm{Hom}(Q_{k}^\Sigma,Q_{k-1}^\Sigma) \otimes \mathrm{Hom}(Q_{k-1}^\Sigma, P_k^\Sigma) \rightarrow 
\mathrm{Hom}(Q_k^\Sigma,P_k^\Sigma).$}
\label{fig:ngonsample}
\end{figure}

Elements of the braid group $B_{m+1}$ acting on $(D_m,\Delta)$ lift to elements
of the mapping class group of the double cover $\Sigma$; specifically, the Artin
generator $\sigma_k$ lifts to the Dehn twist about the simple closed curve
$P_k^\Sigma=\pi_\Sigma^{-1}(p_k)$, where $p_k$ is the line segment in $D_m$ 
joining the two points labeled $k-1$ and $k$ (see Figure \ref{fig:ngonsample}).
We denote by $\hat{\sigma}$ the mapping class group element which lifts a braid
$\sigma\in B_{m+1}$. With this understood, there are two natural ways of 
associating an $A_\infty$-bimodule over $B^{HF}$ to a braid~$\sigma$. 

On one hand, Lipshitz, Ozsv\'ath and Thurston
\cite{GT10030598} associate to the element $\hat\sigma$ of the mapping class group
a bimodule $\widehat{CFDA}(\hat\sigma)$ over the strands algebra, defined in terms
of a suitable Heegaard diagram for the ``mapping cylinder'' of $\hat\sigma$, i.e.\
the 3-manifold $\Sigma\times [0,1]$ equipped with parametrizations of the two
boundary components which differ by the action of $\hat\sigma$ (see 
\cite{GT08100687, GT10030598} for details). We denote by
$\cM^{HF}_\sigma$ the 1-moving strand part of $\widehat{CFDA}(\hat\sigma)$;
this is an $A_\infty$-bimodule over $B^{HF}$ (in fact a ``type DA'' bimodule,
which has nicer algebraic properties).

On the other hand, $\hat\sigma$ acts on the Fukaya category of $\pi_\Sigma$,
and the $A_\infty$-functor induced by $\hat\sigma$ naturally yields a bimodule
over $\mathcal{F}(\pi_\Sigma)$, hence over $B^{\Sigma}$. More concretely,
following \cite{GT10032962} (see also \cite{GT10051248}) we set
\[\cM_{\sigma}^{\Sigma}:=\bigoplus_{i,j=0}^m \mathrm{Hom}_{\mathcal{F}(\pi_\Sigma)}
\left(Q_i^\Sigma, \hat\sigma(Q_j^\Sigma)\right),\] 
which is naturally an $A_\infty$-bimodule over $B^\Sigma\simeq B^{HF}$.

\begin{figure}
\begin{center}
\resizebox{4.5in}{!}{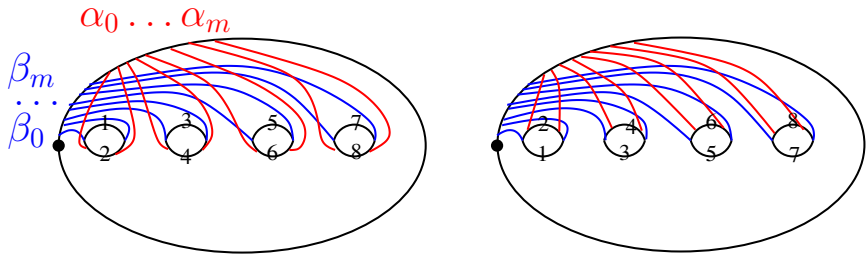}
\end{center}
\caption{A Heegaard diagram for the identity mapping class on $\Sigma$ (the left and right hand side pictures are 
glued according to the numbers). Note that the $\alpha$ and $\beta$ arcs are perturbed
copies of the arcs $Q_k^\Sigma$.}
\label{fig:heegaardAZ}
\end{figure}

\begin{lemma}
The $A_\infty$-bimodules $\cM_\sigma^\Sigma$ and $\cM_\sigma^{HF}$ are quasi-isomorphic.
\end{lemma}

\proof
It is known \cite{GT10030598} that the bordered bimodule 
$\widehat{CFDA}(\mathrm{id})$ is quasi-isomorphic to the strands 
algebra viewed as a bimodule over itself; therefore
$\cM_{\mathrm{id}}^{HF}\simeq B^{HF}\simeq B^{\Sigma}\simeq \cM_{\mathrm{id}}^{\Sigma}$
(as bimodules). We now give a more geometric interpretation, still in the
case $\sigma=\mathrm{id}$.

Following the terminology in \cite{GT10051248}, denote by $AZ$ the bordered Heegaard
diagram depicted in Figure \ref{fig:heegaardAZ}, in which the $\alpha$-arcs and the
$\beta$-arcs are obtained from $Q_k^{\Sigma}$ by pushing the end points along the
boundary of $\Sigma$, in such a manner that the end points of the $\alpha$-arcs
all lie {\it before} those of the $\beta$-arcs along the oriented intervals
$\partial\Sigma\setminus \{z_\pm\}$. Then the 1-moving strand part of the 
$A_\infty$-bimodule $\widehat{CFAA}(AZ)$ is quasi-isomorphic to $\cM^{HF}_{\mathrm{id}}\simeq
B^{HF}$; in fact, $\widehat{CFAA}(AZ)\simeq 
\widehat{CFDA}(\mathrm{id})\simeq \mathcal{A}(\mathcal{Z}_{\mathcal{Q}})$ 
\cite{GT10014323,GT09081106,GT10051248}.
Thus it is enough to show that the 1-moving strand part of
$\widehat{CFAA}(AZ)$ is quasi-isomorphic to $\cM^\Sigma_{\mathrm{id}}=B^{\Sigma}$.

To understand this, recall that morphisms in $\mathcal{F}(\pi_\Sigma)$ are computed by perturbing
the arcs to the same positions used in the Heegaard diagram $AZ$. Hence,
the generators of $\mathrm{Hom}(Q_i^\Sigma,Q_j^\Sigma)$ are precisely the
intersection points between $\beta_i$ and $\alpha_j$, i.e.\ the generators
of the 1-moving strand type AA bimodule. Moreover, the structure maps 
$m_{(k|1|\ell)}$ count:

\begin{itemize}
\item in the case of the type AA bordered Floer bimodule $\widehat{CFAA}(AZ)$,
holomorphic strips in $\Sigma$
connecting two generators of the Heegaard-Floer complex, and
with $k$ (resp.\ $\ell$) additional strip-like ends corresponding to
chords between $\beta$ (resp.\ $\alpha$) arcs;
\item in the case of $\cM^{\Sigma}_{\mathrm{id}}$ (bimodule over the Fukaya
category), rigid holomorphic polygons
bounded by $k+1$ successively perturbed copies of the $\beta$-arcs and $\ell+1$
successively perturbed copies of the $\alpha$-arcs (in the setting
of \cite{GT10014323}; with other definitions of $\mathcal{F}(\pi_\Sigma)$
the interpretation is slightly different).
\end{itemize}

\noindent However, there is a natural correspondence between these two types of objects;
see Proposition 6.5 of \cite{GT10014323} and its proof for details.

In the case of an arbitrary braid $\sigma$, denote by $\hat\sigma(AZ)$ the
bordered Heegaard diagram obtained from $AZ$ by having $\hat\sigma$ act on
the $\alpha$-arcs (leaving the $\beta$-arcs unchanged). From the perspective of
Heegaard-Floer theory, the bordered 3-manifold represented by $\hat\sigma(AZ)$
differs from that corresponding to $AZ$ by a reparametrization of its $\alpha$-boundary
via the action of $\hat\sigma$, or equivalently, by attaching the mapping cylinder
of $\hat\sigma$. Thus $$\widehat{CFAA}(\hat\sigma(AZ))\simeq
\widehat{CFAA}(AZ)\,\tildeotimes\,\widehat{CFDA}(\hat\sigma)\simeq
\widehat{CFDA}(\hat\sigma).$$ Hence $\cM^{HF}_\sigma$ is quasi-isomorphic to
the 1-moving strands part of $\widehat{CFAA}(\hat\sigma(AZ))$. On the other
hand, by the same argument as above this latter bimodule is quasi-isomorphic
to $\cM^\Sigma_\sigma=\bigoplus_{i,j} \mathrm{Hom}_{\mathcal{F}(\pi_\Sigma)}(Q_i^\Sigma,\hat\sigma(Q_j^
\Sigma))$.
\endproof

\begin{remark}\label{rmk:hamperturb}
The comparison between the higher structure maps of the bimodules defined
from $\mathcal{F}(\pi_\Sigma)$ and from bordered Floer homology is easiest
in the setup of \cite{GT10014323}, where a specific Hamiltonian flow is
used to perturb the Lagrangians and ensure transversality, and the structure
maps count honest holomorphic curves bounded by successively perturbed 
copies of the Lagrangians (see Lemma 4.7 of \cite{GT10014323}: the 
definition of the partially wrapped category is much more cumbersome, 
but in the case at hand it simplifies vastly). 

The reader who wishes to reproduce this argument using Seidel's
definition of $\mathcal{F}(\pi_\Sigma)$ instead is referred to
\cite{SeidelLF1}, where the directed Fukaya category
is recast in terms of the symplectic geometry 
of the thimbles and solutions to Floer's equation with Hamiltonian
perturbations. The relevant Hamiltonians behave essentially in the same
manner as that of \cite{GT10014323}, and the main remaining difference is 
that one counts solutions to a perturbed holomorphic curve equation with
boundary on the original Lagrangians, rather than (cascades of)
honest holomorphic curves with boundary on perturbed copies of the
Lagrangians. The two counts can be compared by a fairly standard argument,
or alternatively the proof of \cite[Proposition 6.5]{GT10014323} can be
adapted to that setting.
\end{remark}

If a braid $\sigma$ can be expressed in terms of the Artin generators as
$\sigma=\sigma_{k_1}^\pm\dots\sigma_{k_n}^\pm$, then its lift can be written as
$\hat\sigma={\hat\sigma_{k_1}}^\pm\dots{\hat\sigma_{k_n}}^\pm$, and the
pairing theorem for CFDA bimodules \cite{GT08100687, GT10030598} implies that
$$\cM^{HF}_\sigma\simeq \cM^{HF}_{\sigma_{k_1}^\pm}\,\tildeotimes_{B^{HF}}
\dots\,\tildeotimes_{B^{HF}}\cM^{HF}_{\sigma_{k_n}^\pm}.$$
Thus it is enough to understand the bimodules $\cM^{HF}_{\sigma_k^\pm}\simeq
\cM^{\Sigma}_{\sigma_k^\pm}$ associated to the Artin generators
and their inverses.  We do this working in the category $\mathcal{F}(\pi_\Sigma)$.
Recall that morphism spaces in that category are defined by Lagrangian Floer
theory after a suitable perturbation (so the end points of arcs lie in the correct
order along the boundary of $\Sigma$); in particular they are generated by 
intersection points.


Focusing first on $\cM_{\sigma_k^+}^{HF}$, and recalling that $\hat\sigma_k^+$ is
the positive Dehn twist about $P_k^\Sigma$, Seidel's exact triangle
for Lagrangian Floer homology \cite{MR1978046} tells us that, for each $i,j \in \{0, \ldots, m\}$, 
$\mathrm{Hom}\left(Q_i^\Sigma,\hat\sigma_k^+(Q_j^\Sigma)\right)$ is quasi-isomorphic
to the complex
\[\xymatrix{0 \ar[r] & \mathrm{Hom}\left(Q_i^\Sigma,P_k^\Sigma\right) \otimes 
\mathrm{Hom}\left(P_k^\Sigma,Q_j^\Sigma\right) \ar[r]^{\hskip 40pt \beta_k^{HF} } & 
\mathrm{Hom}\left(Q_i^\Sigma,Q_j^\Sigma\right) \ar[r] & 0},\] where $\beta_k^{HF}$ 
is the Floer product map (cf. \cite{MR1978046}) induced by counting holomorphic triangles 
in $\Sigma$ whose sides lie on (suitable perturbations of) $Q_i^\Sigma,P_k^\Sigma, Q_j^\Sigma$, appearing in counterclockwise order around the boundary.
Moreover, these quasi-isomorphisms are compatible with Floer products, in the sense
that in $D_\infty(B^{HF})$ the bimodule $\cM^{HF}_{\sigma_k^+}$ is equivalent to
the complex of bimodules obtained by taking the direct sum of the above complexes 
over all $i,j$.

In analogy to the previous section, we introduce the $A_\infty$-modules
\[P_k^{HF} := \bigoplus_{i=0}^m \mathrm{Hom}_{\mathcal{F}(\pi_\Sigma)}(Q_i^\Sigma,P_k^\Sigma)\hskip20pt \mbox{and} 
\hskip20pt {}_kP^{HF} := \bigoplus_{j=0}^m \mathrm{Hom}_{\mathcal{F}(\pi_\Sigma)}(P_k^\Sigma, Q_j),\]
which allows us to write
\[\cM_{\sigma_k^+}^{HF}\ \simeq \Bigl\{
\xymatrix{0 \ar[r] & P_k^{HF} \otimes {}_kP^{HF} \ar[r]^{\hskip 20pt\beta_k^{HF}} & \ar[r] B^{HF} \ar[r] &0}
\Bigr\}\] 
Like the linear term described above, the higher terms 
\begin{multline*}
(\beta_k^{HF})_{(n_1|1|n_2)}:\\\bigoplus_{\substack{i_0,\dots,i_{n_1}\\
j_0,\dots,j_{n_2}}} \mathrm{Hom}(Q_{i_{n_1}}^\Sigma,Q_{i_{n_1-1}}^\Sigma)
\otimes\dots\otimes \mathrm{Hom}(Q_{i_1}^\Sigma,Q_{i_0}^\Sigma)
\otimes \mathrm{Hom}(Q_{i_0}^\Sigma,P_k^\Sigma)\otimes\hfill\\[-12pt]
\otimes \mathrm{Hom}(P_k^\Sigma,Q_{j_0}^\Sigma)\otimes
\dots\otimes \mathrm{Hom}(Q_{j_{n_2-1}}^\Sigma,Q_{j_{n_2}}^\Sigma)\longrightarrow
\bigoplus_{i_{n_1},j_{n_2}}\mathrm{Hom}(Q_{i_{n_1}}^\Sigma,Q_{j_{n_2}}^\Sigma)$$
\end{multline*}
of the $A_\infty$-bimodule 
homomorphism $\beta_k^{HF}$ count rigid holomorphic polygons in $\Sigma$ whose
sides lie on (suitable perturbations of) $Q_{i_{n_1}}^\Sigma,\dots,Q_{i_0}^\Sigma,
P_k^\Sigma,Q_{j_0}^\Sigma,\dots,Q_{j_{n_2}}^\Sigma$ in that order.

Similarly, $\cM_{\sigma_k^-}^{HF}$ is equivalent in $D_\infty\left(B^{HF}\right)$ to 
the direct sum of the complexes
\[\xymatrix{0 \ar[r] & \mathrm{Hom}\left(Q_i^\Sigma,Q_j^\Sigma\right) 
\ar[r]^{\hskip -30pt \gamma_k^{HF} } & \mathrm{Hom}\left(Q_i^\Sigma,P_k^\Sigma\right) 
\otimes \mathrm{Hom}\left(P_k^\Sigma,Q_j^\Sigma\right)\ar[r] & 0},\]
where $\gamma_k^{HF}$ is induced by counting holomorphic triangles in $\Sigma$ 
whose sides lie on (suitable perturbations of) $P_k^\Sigma,Q_i^\Sigma,Q_j^\Sigma$, 
appearing in counterclockwise order around the boundary. Thus,
in $D_\infty(B^{HF})$ we have
\[\cM_{\sigma_k^-}^{HF}\ \simeq \Bigl\{
\xymatrix{0 \ar[r] & B^{HF} \ar[r]^{\hskip -20pt\gamma_k^{HF}} & \ar[r] P_k^{HF} \otimes {}_kP^{HF} \ar[r] &0}
\Bigr\}\]
where the higher terms of the $A_\infty$-bimodule homomorphism $\gamma_k^{HF}$ again count
rigid holomorphic polygons in $\Sigma$.

We remark that, in our very simple setting, these counts are equivalent (by the Riemann mapping theorem) to counts of topological immersed
triangles in $\Sigma$ with the stated boundary conditions, and satisfying a local
convexity condition at their corners.

\subsection{Explicit calculations}
\begin{figure}
\begin{center}
\resizebox{4.5in}{!}{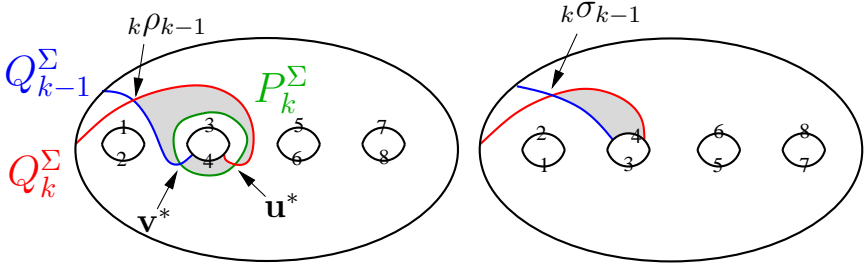}
\end{center}
\caption{The holomorphic triangles giving rise to the nontrivial multiplication maps
\hfil $m_{(1|1|0)}: 
\mathrm{Hom}(Q_{k}^\Sigma,Q_{k-1}^\Sigma) \otimes \mathrm{Hom}(Q_{k-1}^\Sigma, P_k^\Sigma) \rightarrow 
\mathrm{Hom}(Q_k^\Sigma,P_k^\Sigma).$
The other nontrivial multiplication maps can be seen in a similar manner.}
\label{fig:nontrivialmult}
\end{figure}

We now make the above story more explicit, by determining the left (resp., right)
$A_\infty$-modules $P_k^{HF}$ (resp., ${}_kP^{HF}$) and the maps $\beta_k^{HF}$ and
$\gamma_k^{HF}$. Since $P_k^\Sigma$ intersects
$Q_{k-1}^\Sigma$ and $Q_k^\Sigma$ transversely once each and is disjoint from all
the other $Q_j^\Sigma$, the vector spaces underlying these modules have rank 2.
The multiplication maps
\[m_{(n|1|0)}: \left(B^{HF}\right)^{\otimes n} \otimes P_k^{HF} \rightarrow P_k^{HF}
\qquad \text{and}\qquad
m_{(0|1|n)}: {}_kP^{HF} \otimes \left(B^{HF}\right)^{\otimes n} \rightarrow {}_kP^{HF}\]  
are given by counting holomorphic ($n+2$)--gons in $\Sigma$ as in Figure \ref{fig:ngonsample}.  Again letting the two generators of $P_k^{HF}$ (resp., of ${}_kP^{HF}$) be denoted by ${\bf u}^*, {\bf v}^*$ (resp., by {\bf u}, {\bf v}) and letting $\theta$ represent an element of $B^{HF}$, it is easily verified (see Figure \ref{fig:nontrivialmult}) that the $m_{(1|1|0)}$ (resp., $m_{(0|1|1)}$) multiplication is given by:
\[\theta \cdot {\bf u}^* = \left\{\begin{array}{cl}
  					{\bf u}^* & \mbox{if $\theta = {}_k\Id_k$},\\
					0 & \mbox{otherwise.}\end{array}\right.
\hskip 20pt\theta \cdot {\bf v}^* = \left\{\begin{array}{cl}
					{\bf v}^* & \mbox{if $\theta = {}_{k-1}\Id_{k-1}$},\\
					{\bf u}^* & \mbox{if $\theta = {}_k{\rho}_{k-1}$ or ${}_k{\sigma}_{k-1}$},\\
					0 & \mbox{otherwise.}\end{array}\right.\]
(resp., given by:
\[{\bf u} \cdot \theta = \left\{\begin{array}{cl}
					{\bf u} & \mbox{if $\theta = {}_k\Id_k$},\\
					{\bf v} & \mbox{if $\theta = {}_k\rho_{k-1}$ or ${}_k\sigma_{k-1}$},\\
					0 & \mbox{otherwise.}\end{array}\right.
\hskip 20pt  {\bf v} \cdot \theta = \left\{\begin{array}{cl}
  					{\bf v} & \mbox{if $\theta = {}_{k-1}\Id_{k-1}$},\\
					0 & \mbox{otherwise.}\end{array}\right)\]

The multiplications $m_{(1|1|0)}$ and $m_{(0|1|1)}$ are associative. Moreover,
the higher multiplications are all identically zero. One way to see the vanishing
of $m_{(n|1|0)}$ is to
observe that, for any sequence $i_n\ge\dots\ge i_1\ge i_0$ ($n\ge 2$), and 
perturbing $Q_{i_0}^\Sigma,\dots,Q_{i_n}^\Sigma$ so that their end points are
in counterclockwise order along the boundary of $\Sigma$ (but preserving minimal
intersection otherwise), there are no convex $(n+2)$-gons with edges lying successively
on $Q_{i_n}^\Sigma,\dots,Q_{i_0}^\Sigma,P_k^\Sigma$ (and similarly for the vanishing
of $m_{(0|1|n)}$). 

A more conceptual
explanation is that it is possible to find a trivialization of the tangent
bundle of $\Sigma$ and graded lifts \cite{SeidelBook} of the Lagrangians 
$P_k^\Sigma,Q_0^\Sigma,
\dots,Q_m^\Sigma$, and hence a $\mathbb{Z}$-grading by Maslov index on $B^{HF}$ and the modules
$P_k^{HF}$, ${}_kP^{HF}$, with the following properties:
\begin{itemize} 
\item all the generators of $B^{HF}$ have degree 0;
\item the generators $\bf u^*$, $\bf v^*$ of $P_k^{HF}$ have the same degree.
\item the generators $\bf u$, $\bf v$ of ${}_kP^{HF}$ have the same degree.
\end{itemize}
Not all degrees can be taken to be zero: in fact $\deg {\bf u}+\deg {\bf u}^*=
\deg {\bf v}+\deg {\bf v}^*=1$.

Since the maps $m_{(n|1|0)}$ and $m_{(0|1|n)}$ are compatible with the grading and
have degree $1-n$, this forces their vanishing unless $n=1$.
\medskip

We now turn to the $A_\infty$ morphisms $\beta_k^{HF}$ and $\gamma_k^{HF}$.
The calculations are simplified by constraints arising from the
Maslov $\Z$-grading. 

First, we observe that $\beta_k^{HF}$ is a degree-preserving $A_\infty$-homomorphism of
bimodules. Namely, since $(\beta_k^{HF})_{(n_1|1|n_2)}$ corresponds to
a Floer product of order $(n_1+n_2+2)$ in $\mathcal{F}(\pi_\Sigma)$, it has degree 
$-(n_1+n_2)$. However, 
$P_k^{HF}\otimes {}_k P^{HF}$ is concentrated in degree 1, while all the
generators of $B^{HF}$ have degree 0. Therefore, the only non-trivial terms
in $\beta_k^{HF}$ are those of degree $-1$, namely
$(\beta_k^{HF})_{(1|1|0)}$ and $(\beta_k^{HF})_{(0|1|1)}$. In particular
the linear term
$\beta_k^{HF}: \mathrm{Hom}\left(Q_i^\Sigma,P_k^\Sigma\right) \otimes \mathrm{Hom}\left(P_k^\Sigma,Q_j^\Sigma\right) \rightarrow \mathrm{Hom}\left(Q_i^\Sigma,Q_j^\Sigma\right)$
vanishes identically.

Similarly, $\gamma_k^{HF}$, which is an 
$A_\infty$-refinement of the pair of pants coproduct in Floer homology,
has degree $\dim_{\mathbb{C}}(\Sigma)=1$ with respect to the Maslov
$\Z$-grading. Hence, the map $(\gamma_k^{HF})_{(n_1|1|n_2)}$ has degree
$1-(n_1+n_2)$ and, for degree reasons, it must vanish identically unless
$n_1+n_2=0$. Thus, the only nontrivial term of $\gamma_k^{HF}$ is the linear one.

The calculations are further simplified by recalling that
\begin{itemize}
	\item $\mathrm{Hom}\left(Q_i^\Sigma,P_k^\Sigma\right) = \mathrm{Hom}\left(P_k^\Sigma,Q_i^\Sigma\right) = 0$ whenever $i \neq k, k-1$ and
	\item $\mathrm{Hom}\left(Q_i^\Sigma,Q_j^\Sigma\right) = 0$ whenever $i<j$.
\end{itemize}

\begin{figure}
\begin{center}
\resizebox{4.5in}{!}{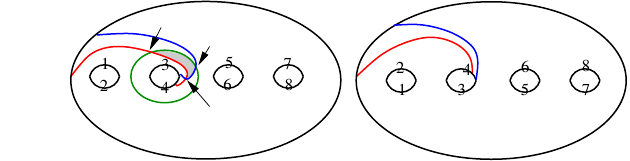}
\end{center}
\caption{The above diagram verifies both that the linear part of \smallskip\newline
\centerline{$\beta_k^{HF}: \mathrm{Hom}\left(Q_k^{\Sigma},P_k^{\Sigma}\right) \otimes 
\mathrm{Hom}\left(P_k^\Sigma,Q_k^\Sigma\right) \rightarrow 
\mathrm{Hom}\left(Q_k^\Sigma,Q_k^\Sigma\right)$}\smallskip\newline
is zero, and that the map \smallskip\newline
\centerline{$\gamma_k^{HF}: \mathrm{Hom}\left(Q_k^\Sigma,Q_k^\Sigma\right) 
\rightarrow \mathrm{Hom}\left(Q_k^\Sigma,P_k^{\Sigma}\right) \otimes 
\mathrm{Hom}\left(P_k^\Sigma,Q_k^\Sigma\right)$}\smallskip\newline 
sends ${}_k\Id_k \in 
\mathrm{Hom}\left(Q_k^\Sigma,Q_k^\Sigma\right)$ to 
${\bf u}^* \otimes {\bf u} \in \mathrm{Hom}\left(Q_k^{\Sigma},P_k^{\Sigma}\right) 
\otimes \mathrm{Hom}\left(P_k^\Sigma,Q_k^\Sigma\right)$.
\smallskip\newline
By definition, these maps count holomorphic triangles with boundary on $P_k^\Sigma$
and on two perturbed copies of $Q_k^\Sigma$, denoted by ${(Q_k^\Sigma)}_1$ and 
${(Q_k^\Sigma)}_2$ in the picture; in counterclockwise order, the successive edges 
must lie on ${(Q_k^\Sigma)}_1,P_k^\Sigma,{(Q_k^\Sigma)}_2$ for $\beta_k^{HF}$, and on
$P_k^\Sigma,{(Q_k^\Sigma)}_1,{(Q_k^\Sigma)}_2$ for $\gamma_k^{HF}$. Hence,
the shaded topological triangle does not contribute to $\beta_k^{HF}$, because
its boundary has the incorrect orientation, hence it does not admit a holomorphic 
representative.  However, it does contribute to the map $\gamma_k^{HF}$.
Computations for the pairs $(i,j) = (k,k-1), (k-1,k-1)$ are similarly straightforward.}
\label{fig:BetakSample}
\end{figure}

\begin{lemma} \label{lem:MCgammakHF}
$\gamma_k^{HF}: B^{HF} \rightarrow P_k^{HF} \otimes {}_kP^{HF}$ is the bimodule map determined by
\[{}_i\Id_i \mapsto \left\{\begin{array}{cl}
				{\bf v}^* \otimes {\bf v} & \mbox{when } i = k-1,\\
				{\bf u}^* \otimes {\bf u} & \mbox{when } i=k, \mbox{ and}\\
				0 & \mbox{otherwise}\end{array}\right.\]
and by associativity with respect to the multiplication. Moreover,
the higher order maps $(\gamma_k^{HF})_{(n_1|1|n_2)}$ vanish identically
for $(n_1,n_2)\neq (0,0)$.
\end{lemma}

\begin{proof}
The map $\gamma_k^{HF}: \mathrm{Hom}\left(Q_i^\Sigma,Q_j^\Sigma\right) \rightarrow \mathrm{Hom}\left(Q_i^\Sigma,P_k^\Sigma\right) \otimes \mathrm{Hom}\left(P_k^\Sigma,Q_j^\Sigma\right)$ is $0$ unless $(i,j) = (k,k), (k-1,k-1),$ or
$(k,k-1)$, since in all other cases either the domain or the target is zero.
The nontrivial cases are then determined by counting immersed triangles in
$\Sigma$; the case $(i,j)=(k,k)$ is shown in Figure \ref{fig:BetakSample}. 
By inspection, we see that  $\gamma_k^{HF}$ is given by:
\begin{itemize}
	\item When $(i,j) = (k,k)$ or $(k-1,k-1)$, $\gamma_k^{HF}$ sends the unique generator of $\mathrm{Hom}\left(Q_i^\Sigma,Q_j^\Sigma\right)$ to the unique generator of $\mathrm{Hom}\left(Q_i^\Sigma,P_k^\Sigma\right) \otimes \mathrm{Hom}\left(P_k^\Sigma,Q_j^\Sigma\right)$, and
	\item When $(i,j) = (k,k-1)$, $\gamma_k^{HF}$ sends both ${}_k\rho_{k-1}$ and ${}_k\sigma_{k-1} \in \mathrm{Hom}\left(Q_i^\Sigma,Q_j^\Sigma\right)$ to the unique generator of $\mathrm{Hom}\left(Q_i^\Sigma,P_k^\Sigma\right) \otimes \mathrm{Hom}\left(P_k^\Sigma,Q_j^\Sigma\right)$.
\end{itemize}
The vanishing of the higher maps follows from the degree argument explained
above.
\end{proof}

The story for $\beta_k^{HF}$ is slightly more complicated, because the
maps 
$$(\beta_k^{HF})_{(1|1|0)}:\mathrm{Hom}(Q_{i_1}^\Sigma,Q_{i_0}^\Sigma)\otimes
\mathrm{Hom}(Q_{i_0}^\Sigma,P_k^\Sigma)\otimes 
\mathrm{Hom}(P_k^\Sigma,Q_j^\Sigma)\longrightarrow 
\mathrm{Hom}(Q_{i_1}^\Sigma,Q_j^\Sigma)$$ and 
$$(\beta_k^{HF})_{(0|1|1)}:\mathrm{Hom}(Q_{i}^\Sigma,P_k^\Sigma)\otimes
\mathrm{Hom}(P_k^\Sigma,Q_{j_0}^\Sigma)\otimes
\mathrm{Hom}(Q_{j_0}^\Sigma,Q_{j_1}^\Sigma)\longrightarrow \mathrm{Hom}(Q_{i_1}^\Sigma,
Q_j^\Sigma),$$
which count holomorphic 4-gons in $\Sigma$, depend on the choice of
Hamiltonian perturbations used to resolve triple intersections at the branch
points of $\pi_\Sigma$. (Of course, the behavior of Lagrangian Floer
homology under Hamiltonian isotopies guarantees that the maps obtained from
different choices are homotopic.) To fix a convention, we perturb
$P_k^\Sigma$ away from the branch points of $\pi_\Sigma$ in such a way 
that its intersections with $Q_k^\Sigma$ and $Q_{k-1}^\Sigma$ occur on the 
sheet of the double cover that contains the generators ${}_i\rho_j$. With
this understood, we have:

\begin{lemma} \label{lem:MCbetakHF}
The only nontrivial terms of $\beta_k^{HF}$ are:
$$(\beta_k^{HF})_{(1|1|0)}:\left\{\begin{array}{rcll}
({}_i\rho_k,\, {\bf u}^*\otimes {\bf u})&\mapsto&{}_i\rho_k&(i \geq k+1)\medskip\\
({}_k\sigma_{k-1},\,{\bf v}^*\otimes {\bf u})&\mapsto&{}_k\Id_k\\
({}_i\rho_{k-1},\,{\bf v}^*\otimes {\bf u})&\mapsto&{}_i\rho_k&(i \geq k+1)\\
({}_i\sigma_{k-1},\,{\bf v}^*\otimes {\bf u})&\mapsto&{}_i\sigma_k&(i \geq k+1)\medskip\\
({}_i\rho_{k-1},\,{\bf v}^*\otimes {\bf v})&\mapsto&{}_i\rho_{k-1}&(i\ge k)\\
\end{array}\right.$$
$$\ \ \text{and}\ \ (\beta_k^{HF})_{(0|1|1)}:\left\{\begin{array}{lcll}
({\bf u}^*\otimes {\bf u},\,{}_k\rho_j)&\mapsto&{}_k\rho_j&(j\le k-1)\medskip\\
({\bf v}^*\otimes {\bf u},\,{}_k\sigma_{k-1})&\mapsto&{}_{k-1}\Id_{k-1}\\
({\bf v}^*\otimes {\bf u},\,{}_k\rho_j)&\mapsto&{}_{k-1}\rho_j&(j \leq k-2)\\
({\bf v}^*\otimes {\bf u},\,{}_k\sigma_j)&\mapsto&{}_{k-1}\sigma_j&(j \leq k-2)\medskip\\
({\bf v}^*\otimes {\bf v},\,{}_{k-1}\rho_j)&\mapsto&{}_{k-1}\rho_j&(j \leq k-2)\\
\end{array}\right.$$
\end{lemma}

\begin{figure}
\begin{center}
\resizebox{4.5in}{!}{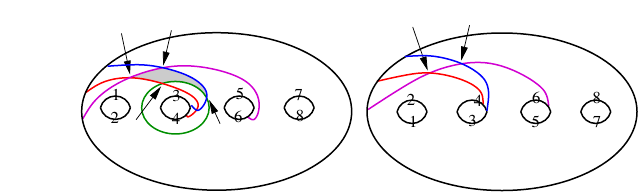}
\end{center}
\caption{The above diagram verifies that
\smallskip\newline\centerline{
$(\beta_k^{HF})_{(1|1|0)}({}_i\rho_k, {\bf u}^*\otimes {\bf
u})={}_i\rho_k$ \ and \ $(\beta_k^{HF})_{(1|1|0)}({}_i\sigma_k, 
{\bf u}^*\otimes {\bf u})=0$}\smallskip\newline
for $i>k$. By definition, \smallskip\newline
$\beta^{HF}_{k\ (1|1|0)}:\mathrm{Hom}(Q_i^\Sigma,
Q_k^\Sigma)\otimes \mathrm{Hom}(Q_k^\Sigma,P_k^\Sigma)
\otimes \mathrm{Hom}(P_k^\Sigma,Q_k^\Sigma)\to \mathrm{Hom}(
Q_i^\Sigma,Q_k^\Sigma)$\smallskip\newline
counts rigid holomorphic 4-gons with successive edges, in 
counterclockwise order, on perturbed copies of $Q_i^\Sigma$,
$Q_k^\Sigma$ (denoted ${(Q_k^\Sigma)}_1$), $P_k^\Sigma$, and
$Q_k^\Sigma$ again (denoted ${(Q_k^\Sigma)}_2$). The only contribution comes
from the shaded region.}
\label{fig:BetakSample2}
\end{figure}

\begin{figure}
\begin{center}
\resizebox{4.5in}{!}{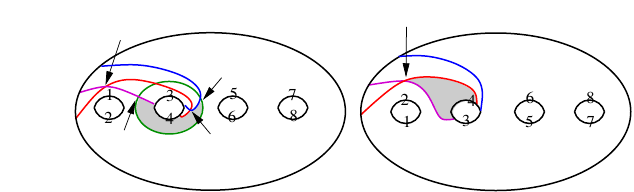}
\end{center}
\caption{The above diagram verifies that
\smallskip\newline\centerline{
$(\beta_k^{HF})_{(1|1|0)}({}_k\rho_{k-1}, {\bf v}^*\otimes {\bf
u})=0$ \ and \ $(\beta_k^{HF})_{(1|1|0)}({}_k\sigma_{k-1}, 
{\bf v}^*\otimes {\bf u})={}_k\Id_k$.}}
\label{fig:BetakSample3}
\end{figure}

\proof By definition, $(\beta_k^{HF})_{(1|1|0)}$ counts
rigid holomorphic 4-gons in $\Sigma$ whose successive edges, in
counterclockwise order, lie on suitably perturbed copies of the following
Lagrangians: $Q_i^\Sigma$; either $Q_k^\Sigma$ (for ${\bf u}^*$) or 
$Q_{k-1}^\Sigma$ (for ${\bf v}^*$); $P_k^\Sigma$; and either $Q_k^\Sigma$
(for ${\bf u}$) or $Q_{k-1}^\Sigma$ (for~${\bf v})$. The count depends
on the perturbations, so we have to be more specific. 

Since we are working in the Fukaya category $\mathcal{F}(\pi_\Sigma)$, the
various arcs must be perturbed by Hamiltonian isotopies which ensure
that their end points are suitably ordered along $\partial\Sigma$; these
perturbations are responsible for the intersection points corresponding 
to the generators ${}_i\rho_k$ and ${}_i\sigma_k$ (resp.\ 
${}_i\rho_{k-1}$, ${}_i\sigma_{k-1}$),
which we take to lie close to the boundary of $\Sigma$. By contrast, the
intersection points corresponding to the generators
${\bf u}^*,{\bf u}$ and
${}_k\Id_k$ normally all lie at the $k$-th branch point of $\pi_\Sigma$, 
and perturbations are needed to avoid triple intersections. 
As mentioned above, we achieve this by choosing a Hamiltonian which pushes
$P_k^\Sigma$ slightly towards the ``$\rho$'' side of the surface.
Likewise for ${\bf v}^*,{\bf v}$ and ${}_{k-1}\Id_{k-1}$.

With this understood, the calculation simply becomes a matter of drawing
the relevant diagrams and looking for immersed four-gons with locally convex
corners. The first two cases are shown on 
Figures \ref{fig:BetakSample2} and \ref{fig:BetakSample3}; the others
are similar.
\endproof

As a consistency check, it is not hard to verify that the map $\beta_k^{HF}$
is indeed an $A_\infty$-homomorphism, namely for all $a_1,a_2\in B^{HF}$ and
$m\in P_k^{HF}\otimes{}_kP^{HF}$ we have the identities
$$\beta_{k\ (1|1|0)}^{HF}(a_1a_2,m)+
\beta_{k\ (1|1|0)}^{HF}(a_1,a_2m)+
a_1\beta_{k\ (1|1|0)}^{HF}(a_2,m)=0,$$
$$\beta_{k\ (0|1|1)}^{HF}(m,a_1a_2)+
\beta_{k\ (0|1|1)}^{HF}(ma_1,a_2)+
\beta_{k\ (0|1|1)}^{HF}(m,a_1)a_2=0,$$
$$a_1\beta_{k\ (0|1|1)}^{HF}(m,a_2)+
\beta_{k\ (0|1|1)}^{HF}(a_1m,a_2)+
\beta_{k\ (1|1|0)}^{HF}(a_1,m)a_2+
\beta_{k\ (1|1|0)}^{HF}(a_1,ma_2)=0.$$

\section{A spectral sequence from the Khovanov-Seidel to the bordered Floer algebra} \label{sec:Filtalg}
In Sections \ref{sec:Khalgmod} and \ref{sec:HFalgmod} we showed how to use the data of a basis, $\widetilde{\mathcal{Q}}$, to construct 
\begin{itemize}
	\item a graded algebra, $B^{Kh}$, using a construction of Khovanov-Seidel in \cite{MR1862802} and 
	\item a (graded) algebra $B^{HF}$, using ideas of Lipshitz-Ozsv{\'a}th-Thurston in \cite{GT08100687} as generalized by Zarev in \cite{GT09081106} and reinterpreted by the first author in \cite{GT10014323}.  
\end{itemize}
	
In this section, we establish the existence of a spectral sequence connecting $B^{Kh}$ and $B^{HF}$.  Explicitly, we prove:

\begin{theorem} \label{thm:filtdiffalg} Let \[B^{Kh} := H_*\left(\bigoplus_{i,j=0}^m {\rm Hom}_A(Q_i,Q_j)\right)\] be the homology of the Hom algebra associated to the basis $\widetilde{\mathcal{Q}}$ and let $B^{HF} := \mathcal{A}\left(\mathcal{Z}_\mathcal{Q},1\right)$ be the $1$--moving strands algebra associated to the arc diagram, $\mathcal{Z}_\mathcal{Q}$.  There exists a filtration on $B^{HF}$ whose associated graded algebra is isomorphic, as an ungraded algebra, to $B^{Kh}$.  Accordingly, one obtains a spectral sequence whose $E^1$ page is isomorphic to $B^{Kh}$ and whose $E^\infty$ page is isomorphic to $B^{HF}$.
\end{theorem}

\begin{remark} The observant reader will at this point notice that the spectral sequence described in the statement of Theorem \ref{thm:filtdiffalg} must be somewhat unusual, since $B^{HF}$ is not a {\em dg algebra} but an {\em algebra}; hence, the induced differential on the associated graded page  is necessarily trivial and the associated spectral sequence on \F--vector spaces collapses immediately.  This should perhaps not be surprising, as we have $\mbox{dim}\left({}_iB^{Kh}_j\right) = \mbox{dim}\left({}_iB^{HF}_j\right)$ for each $i, j \in \{0, \ldots, m\}$.  On the other hand, $B^{Kh}$ and $B^{HF}$ are not isomorphic as algebras.  The filtration serves only to alter the multiplicative structure on the underlying algebra and not to change the dimensions of the underlying \F--vector spaces. \end{remark}

We pave the way for a proof of Theorem \ref{thm:filtdiffalg} by focusing first on a ``toy model" given by the following two lemmas.  Though not logically necessary for the proof of Theorem \ref{thm:filtdiffalg}, we include them 
in order to motivate the definition of the filtration yielding the spectral sequence from $B^{Kh}$ and $B^{HF}$.

\begin{lemma}  \label{lem:ZtwoequivS1} There exists a filtered differential algebra, $\mathcal{C}$, whose associated graded homology algebra is isomorphic to $H^*(S^1)$ and whose total homology algebra is isomorphic to  $H^*(S^0)$.  Furthermore, the associated graded complex and the total complex of $\mathcal{C}$ are formal $A_\infty$ algebras.
\end{lemma}

\begin{proof} We construct $\mathcal{C}$ using a $\Ztwo$--equivariant cochain complex for $H^*(S^1)$.  Specifically, identify $S^1$ with the unit circle in $\C$ and give it the structure of a simplicial complex by placing two $0$--simplices labeled {\bf a} and {\bf b} at $-1$ and $1$, respectively, and two $1$--simplices labeled {\bf A} and {\bf B} along the arcs $\left\{e^{i\theta}|\theta \in [\pi, 0]\right\}$ and $\left\{e^{i\theta}|\theta \in [-\pi, 0]\right\}$, respectively, as in Figure \ref{fig:Z2EquivS1}.  Let ${\bf a^*}$ (resp. ${\bf b^*}$, ${\bf A^*}$, ${\bf B^*}$) represent the $\Ztwo$  cochain that assigns $1$ to ${\bf a}$ (resp., {\bf b}, {\bf A}, {\bf B}) and $0$ to all other simplices in the basis.

\begin{figure}
\begin{center}
\resizebox{2in}{!}{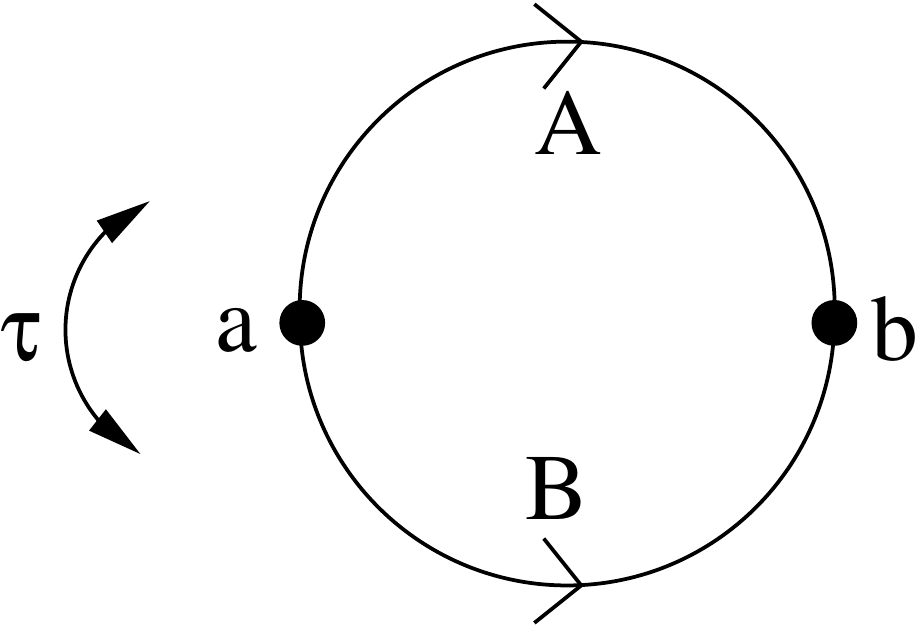}
\end{center}
\caption{A \Ztwo--equivariant chain complex for $S^1$.}
\label{fig:Z2EquivS1}
\end{figure}

The filtered differential algebra, $\mathcal{C}$, is generated by ${\bf a^*}$, ${\bf b^*}$, ${\bf A^*}$, and ${\bf B^*}$ with multiplication given by the cup product on cochains (cf. \cite{MR1867354}):

\begin{equation} \label{eqn:CupProd}
\begin{array}{c|cccc} 
	\cup		& {\bf a^*}	& {\bf b^*}	& {\bf A^*}	& {\bf B^*}\\
			\hline
	{\bf a^*}	& {\bf a^*} & 0		& {\bf A^*} & {\bf B^*}\\
	{\bf b^*}	& 0		& {\bf b^*}	& 0		& 0\\
	{\bf A^*}	& 0		& {\bf A^*} & 0		& 0\\
	{\bf B^*} 	& 0		& {\bf B^*} & 0		& 0\\
	\end{array}
\end{equation}

There are two commuting differentials, $\delta$ and $\partial_\tau$, on $\mathcal{C}$, giving $\mathcal{C}$ the structure of a differential algebra: 
\begin{itemize}
	\item $\delta$ is the standard coboundary map on the simplicial cochain complex (hence satisfies the Leibniz rule with respect to the cup product multiplication), and 
	\item $\partial_\tau =$ {\textbb 1} + $\tau$, where $\tau$ is the involution on the cochain complex induced by complex conjugation on \C.  One easily checks that $\partial_\tau$ satisfies the Leibniz rule with respect to the cup product multiplication.
\end{itemize}

We have the following two-step filtration $ \mathcal{F}_{-1} \subseteq  \mathcal{F}_0 \subseteq \mathcal{F}_1$: \[ 0 \subseteq \ker(\partial_\tau) \subseteq \mathcal{C}\] on $(\mathcal{C}, \delta + \partial_\tau)$.  This gives $\mathcal{C}$ the structure of a filtered algebra, since $\mathcal{F}_i \cdot \mathcal{F}_j \subseteq \mathcal{F}_{i+j}$ for all $i, j$.\footnote{The only non-trivial check that must be performed is that $\mathcal{F}_0 \cdot \mathcal{F}_0 \subseteq \mathcal{F}_0$, but this follows from the fact that $\partial_\tau$ satisfies the Leibniz rule.}  Furthermore, the associated graded complex is $(\mathcal{C},\delta)$, with homology $H^*(S^1)$ and the homology of the total complex $(\mathcal{C},\delta+\partial_\tau)$ is the cohomology of the fixed point set of $\tau$, i.e., $H^*(S^0)$.
 
We now use Proposition \ref{prop:InducedAinfty} to compute the $A_\infty$ structure on the associated graded complex of $\mathcal{C}$, defining maps $\iota: H^*(S^1) \rightarrow (\mathcal{C},\delta)$, $p: (\mathcal{C},\delta) \rightarrow H^*(S^1)$ and $h: (\mathcal{C},\delta) \rightarrow (\mathcal{C},\delta)$ satisfying the conditions in Equation \ref{eqn:Ainftymap}.

Let {\bf 1} denote the generator of $H^0(S^1)$ and {\bf x} denote the generator of $H^1(S^1)$.  Then we define 
\begin{eqnarray*}
	\iota({\bf 1}) &:=& {\bf a^*} + {\bf b^*}\\
	\iota({\bf x}) &:=& {\bf A^*},
\end{eqnarray*}
\begin{eqnarray*}
	p({\bf a^*}) &:=& {\bf 1}\\
	p({\bf A^*}) = p({\bf B^*}) &:=& {\bf x}\\
	p({\bf b^*}) &:=& 0,
\end{eqnarray*}
and
\begin{eqnarray*}
	h({\bf B^*}) &:=& {\bf b^*}\\
	h({\bf a^*}) = h({\bf b^*}) = h({\bf A^*}) &:=& 0\\
\end{eqnarray*}

An application of Lemma \ref{lem:FormalCond} then implies that the associated graded algebra is formal.

We proceed similarly for $(\mathcal{C},\delta + \partial_\tau)$.  Let $\rho, \sigma$ denote the two generators of $H^*(S^0)$ corresponding to the two connected components of $S^0$.  We define:
\begin{eqnarray*}
	\iota(\rho) &:=& {\bf a^*} + {\bf A^*}\\
	\iota(\sigma) &:=& {\bf b^*} + {\bf A^*},
\end{eqnarray*}
\begin{eqnarray*}
	p({\bf a^*}) &:=& \rho\\
	p({\bf b^*}) &:=& \sigma,\\
	p({\bf A^*}) = p({\bf B^*}) &:=& 0
\end{eqnarray*}
and
\begin{eqnarray*}
	h({\bf B^*}) &:=& {\bf A^*}\\
	h({\bf a^*}) = h({\bf b^*}) = h({\bf A^*}) &:=& 0,
\end{eqnarray*}

Once again, an application of Lemma \ref{lem:FormalCond} implies that the total algebra of $\mathcal{C}$ is formal.
\end{proof}

As noted in the proof of Lemma \ref{lem:ZtwoequivS1}, we have simple descriptions of $H^*(S^1)$ and $H^*(S^0)$ as \F--algebras: \[H^*(S^1) \cong \F[{\bf x}]/{\bf x}^2\] and \[H^*(S^0) := \mbox{Span}_\F\langle \rho,\sigma\rangle,\] with multiplication given by 
\begin{eqnarray*}
	m_2(\rho \otimes \rho) &=& \rho\\
	m_2(\sigma \otimes \sigma) &=& \sigma\\
	m_2(\rho \otimes \sigma) = m_2(\sigma \otimes \rho) &=& 0.
\end{eqnarray*}

Furthermore, the filtration on the filtered differential algebra $\mathcal{C}$ defined in the proof of Lemma \ref{lem:ZtwoequivS1} induces a filtration on $H^*(S^0)$.  Accordingly, we have:

\begin{lemma} \label{lem:ZtwoequivS1collapse}  Consider the following filtration, $\cF_{-1} \subseteq \cF_0 \subseteq \cF_1$, on $H^*(S^0)$:  \[0 \subseteq \mbox{Span}_\F\langle\rho + \sigma\rangle \subseteq H^*(S^0).\]  With respect to this filtration, $H^*(S^0)$ is a well-defined filtered (differential) algebra with associated graded algebra isomorphic to $H^*(S^1)$.
\end{lemma}

\begin{proof}  The claim follows immediately from the observation that the $A_\infty$ quasi-isomorphism $\iota: H^*(S^0) \rightarrow \mathcal{C}$ guaranteed by Lemma \ref{lem:FormalCond} is filtered, hence induces a filtered $A_\infty$ quasi-isomorphism.

However, we find it instructive to give a more direct proof.

First, $H^*(S^0)$ is easily seen to be a well-defined filtered ($A_\infty$) algebra (Definition \ref{defn:FiltAlg}) with respect to the above choice of filtration.  The only non-trivial check that must be performed is that $m_2((\rho + \sigma) \otimes (\rho + \sigma)) \subseteq \cF_0$, which follows since $1 := \rho + \sigma$ is the identity element of $H^*(S^0)$.  Recalling that the multiplication on the associated graded is given by \[m_2: \cF_{r}/\cF_{r-1} \otimes \cF_{s}/\cF_{s-1} \rightarrow \cF_{r+s}/\cF_{r+s-1},\] we see immediately that $1$ is also the multiplicative identity in $\mbox{gr}(H^*(S^0))$, since it lies in filtration level $0$.  

The underlying \F--vector space of the associated graded algebra $\mbox{gr}(H^*(S^0))$ can be described by: \[\cF_n/\cF_{n-1} := \left\{\begin{array}{cl}
		\mbox{Span}_\F\langle 1\rangle & \mbox{if $n=0$,}\\
		\mbox{Span}_\F\langle \rho \rangle & \mbox{if $n=1$,}\\
		0 & \mbox{otherwise.}
		\end{array}\right.\]

Furthermore, \[m_2(\rho \otimes \rho) = \rho = 0 \in \cF_2/\cF_1.\]  Hence, $\mbox{gr}(H^*(S^0))$ is isomorphic to $H^*(S^1)$, by identifying $1, \rho \in \mbox{gr}(H^*(S^0))$ with $1, {\bf x} \in H^*(S^1)$.
\end{proof}



We now proceed to the proof of Theorem \ref{thm:filtdiffalg}.

\begin{proof}[Proof of Theorem \ref{thm:filtdiffalg}]  
Recalling (see Remark \ref{rmk:BHFMatrix}) that $B^{HF}$ is isomorphic to the algebra of lower triangular $(m+1) \times (m+1)$ matrices over $H^*(S^0)$ with only $0$'s and $1$'s on the diagonal, we define the desired filtration, $\cF_{-1} \subseteq \cF_0 \subseteq \cF_1$, on $B^{HF}$ as follows:
\[0 \subseteq \left\{M \in B^{HF} \,\,\vline \,\, \phi_{i,j} \in \{0,1\} \,\,\forall\,\,i>j\right\} \subseteq B^{HF}\]
		
We now claim that the associated graded algebra, $\mbox{gr}\left(B^{HF}\right)$, is isomorphic to $B^{Kh}$.  To see this, note that
\[\cF_n/\cF_{n-1} := \left\{\begin{array}{cl}
		\{M \in B^{HF} \,\,\vline\,\, \phi_{i,j} \in \{0,1\} \,\,\forall\,\,i>j\} & \mbox{when $n=0$,}\\
		\{M \in B^{HF} \,\, \vline \,\, \phi_{i,j} \in \{0,\rho\} \,\,\forall\,\,i > j, \mbox{ and } d_k = 0 \,\,\forall\,\,k\} & \mbox{when $n=1$, and}\\
		0 & \mbox{otherwise.}\end{array}\right.\]

In particular, gr($B^{HF}$) is isomorphic to the algebra of $(m+1) \times (m+1)$ lower triangular matrices over gr($H^*(S^0)$) with only $0$'s and $1$'s on the diagonal, where the filtration on $H^*(S^0)$ is the one described in Lemma \ref{lem:ZtwoequivS1collapse}.  Hence, Lemma \ref{lem:ZtwoequivS1collapse} tells us that gr($B^{HF}$) is isomorphic to $B^{Kh}$ as an \F--algebra, as desired.
\end{proof}

\section{A spectral sequence from the Khovanov-Seidel to the bordered Floer bimodules} \label{sec:Filtdiffmod}

In analogy to Theorem \ref{thm:filtdiffalg}, we prove the following theorem relating the Hom modules described in Section \ref{sec:Khalgmod} to the bordered Floer modules described in Section \ref{sec:HFalgmod}.

Recall that $\widetilde{\mathcal{Q}}$ is the basis (of $\partial$--admissible bigraded curves in normal form) pictured in Figure \ref{fig:Qbasis}.

\begin{theorem} \label{thm:filtdiffmodgen} Let $\sigma \in B_{m+1}$ be a braid, $\cM_\sigma^{Kh}$ the bimodule associated to the pair $(\widetilde{\mathcal{Q}},\sigma)$ in Section \ref{sec:Khalgmod}, and $\cM_\sigma^{HF}$ the bordered Floer bimodule associated to the pair $(\mathcal{Q},\sigma)$ in Section \ref{sec:HFalgmod}.  There exists a filtration on $\cM_\sigma^{HF}$ whose associated graded bimodule is isomorphic (as an ungraded $A_\infty$ bimodule over $B^{Kh}$) to $\cM_{\sigma}^{Kh}$.  Accordingly, one obtains a spectral sequence whose $E^1$ page is isomorphic to $\cM_{\sigma}^{Kh}$ and whose $E^\infty$ page is isomorphic to $\cM_{\sigma}^{HF}$.
\end{theorem}

Note that Theorem \ref{thm:filtdiffalg} is Theorem \ref{thm:filtdiffmodgen} in the special case $\sigma =$ Id.  The proof of Theorem \ref{thm:filtdiffmodgen} proceeds in two steps.  We begin by giving an explicit construction of the filtration in the special case where $\sigma$ is one of the elementary Artin braid generators, $\{\sigma_k^\pm| k = 1, \ldots, m\}$ (Proposition \ref{prop:filtdiffmodelem}).  Then in the general case, $\sigma = \sigma_{k_1}^\pm \cdots \sigma_{k_n}^\pm$, we explain how to construct a filtration and appropriate spectral sequence on the $A_\infty$ module formed as the $A_\infty$ tensor product 
\[\cM^{HF}_{\sigma^\pm_{k_1}} \widetilde{\otimes}_{B^{HF}} \ldots \widetilde{\otimes}_{B^{HF}} \cM^{HF}_{\sigma^\pm_{k_n}}.\]

\begin{proposition} \label{prop:filtdiffmodelem} Let $\sigma_k^\pm \in B_{m+1}$ be an elementary Artin braid generator,  $\cM_{\sigma_k^\pm}^{Kh}$ the bimodule associated to the pair $(\widetilde{\mathcal{Q}},\sigma_k^\pm)$ in Section \ref{sec:Khalgmod}, and $\cM_{\sigma_k^\pm}^{HF}$ the bordered-Floer bimodule associated to the pair $(\mathcal{Q},\sigma_k^\pm)$ in Section \ref{sec:HFalgmod}.   There exists a filtration on $\cM_{\sigma_k^\pm}^{HF}$ whose associated graded bimodule is isomorphic (as an ungraded $A_\infty$ bimodule over $B^{Kh}$) to $\cM_{\sigma_k^\pm}^{Kh}$.  Accordingly, one obtains a spectral sequence whose $E^1$ page is isomorphic to $\cM_{\sigma_k^\pm}^{Kh}$ and whose $E^\infty$ page is isomorphic to $\cM_{\sigma_k^\pm}^{HF}$.

\end{proposition}

\begin{proof}[Proof of Proposition \ref{prop:filtdiffmodelem}]
Guided by the models $\cM_{\sigma_k^\pm}^{Kh}$ and $\cM_{\sigma_k^\pm}^{HF}$ constructed in Sections \ref{sec:Khalgmod} and \ref{sec:HFalgmod}, we turn now to constructing filtrations on the filtered bimodules $\cM_{\sigma_k^\pm}^{HF}$ (over the filtered algebra $B^{HF}$) with the desired properties.

We begin by defining, for each $k \in \{0, \ldots, m\},$ filtrations on $P_k^{HF}$ and ${}_kP^{HF}$.  Since:
\begin{enumerate}
	\item we have already defined (Theorem \ref{thm:filtdiffalg}) a filtration on $B^{HF}$,
	\item the tensor product of two filtered $A_\infty$ modules inherits the structure of a filtered $A_\infty$ module,
	\item the mapping cone of two filtered $A_\infty$ modules inherits the structure of a filtered $A_\infty$ module, and
	\item we have \[\cM_{\sigma_k^+}^{HF} := MC\left(\beta_k^{HF}: (P_k^{HF} \otimes {}_kP^{HF}) \rightarrow B^{HF}\right)\] and \[\cM_{\sigma_k^-}^{HF} := MC\left(\gamma_k^{HF}: B^{HF} \rightarrow (P_k^{HF} \otimes {}_kP^{HF})\{-1\}\right),\]
\end{enumerate}
this will induce a filtration on each $\cM_{\sigma_k^\pm}^{HF}$, as desired.

Recalling that $P_k^{HF} := \mbox{Span}_\F\langle{\bf u}^*,{\bf v}^*\rangle$ (resp., ${}_kP^{HF} := \mbox{Span}_\F\langle{\bf u},{\bf v}\rangle$), we define the filtration, $\mathcal{F}_{-1} \subseteq \mathcal{F}_0 \subseteq \mathcal{F}_1$, on $P_k^{HF}$ to be $0 \subseteq \mbox{Span}\langle {\bf v}^*\rangle \subseteq P_k^{HF}$ (resp., on ${}_kP^{HF}$ to be $0 \subseteq \mbox{Span}\langle {\bf u}\rangle \subseteq {}_kP^{HF}$).  

Verification that $\beta_k^{HF}$ and $\gamma_k^{HF}$ are filtered $A_\infty$ morphisms with respect to this choice of filtration is a straightforward check of a small number of cases, and is left to the reader.

We now must show that the associated graded (homology) of $\cM_{\sigma_k^\pm}^{HF}$ is isomorphic to $\cM_{\sigma_k^\pm}^{Kh}$ as a $\left(\mbox{gr}(B^{HF}) = B^{Kh}\right)$--bimodule.

Since we have already shown (in the proof of Theorem \ref{thm:filtdiffalg}) that the multiplication on $\mbox{gr}(B^{HF})$ matches the multiplication on $B^{Kh}$, all that remains to show is:
\begin{enumerate}
	\item the multiplication of $\mbox{gr}(B^{HF})$ on $\mbox{gr}\left(P_k^{HF} \otimes {}_kP^{HF}\right)$ matches the multiplication of $B^{Kh}$ on $P_k^{Kh} \otimes {}_kP^{Kh}$ and
	\item the maps induced by $\gamma_k^{HF}$ and $\beta_k^{HF}$ on $\mbox{gr}(B^{HF})$ and $\mbox{gr}\left(P_k^{HF} \otimes {}_kP^{HF}\right)$ match the maps $\gamma_k^{Kh}$ and $\beta_k^{Kh}$.
\end{enumerate}

Seeing that the multiplication of $\mbox{gr}(B^{HF})$ on $\mbox{gr}\left(P_k^{HF} \otimes {}_kP^{HF}\right)$ matches the multiplication of $B^{Kh}$ on $P_k^{Kh} \otimes {}_kP^{Kh}$ is a simple check of a small number of cases, bearing in mind that under the isomorphism $\mbox{gr}(B^{HF}) \leftrightarrow B^{Kh}$, we have the identification ${}_i\rho_j \leftrightarrow {}_i{\bf x}_j$.

The map induced by $\gamma_k^{HF}$ on $\mbox{gr}(B^{HF})$ is quickly seen to match the map $\gamma_k^{Kh}$, since $\gamma_k^{HF}$ is a filtered morphism with no higher terms, and the descriptions of $\gamma_k^{Kh}$ (Proposition \ref{prop:MCgammak}) and $\gamma_k^{HF}$ (Lemma \ref{lem:MCgammakHF}) are identical.

Verifying that the map induced by $\beta_k^{HF}$ on $\mbox{gr}\left(P_k^{HF} \otimes {}_kP^{HF}\right)$ matches the map $\beta_k^{Kh}$ is a bit more involved but, again, requires only a handful of checks.  We perform a couple here, leaving the rest to the reader.

Lemma \ref{lem:MCbetakHF} tells us that when $i \geq k+1$: \[\left(\beta_k^{HF}\right)_{(1|1|0)}\left[{}_i\rho_k \otimes ({\bf u}^* \otimes {\bf u})\right] := {}_i\rho_k.\]  But viewed as elements of the associated graded, we have ${}_i\rho_k \in \cF_1/\cF_0(B^{HF})$ and ${\bf u}^* \otimes {\bf u} \in \cF_1/\cF_0\left(P_k^{HF} \otimes {}_kP^{HF}\right)$, and thus the induced associated graded map is:\[\left(\beta_k^{HF}\right)_{(1|1|0)}\left[{}_i\rho_k \otimes ({\bf u}^* \otimes {\bf u})\right] := {}_i\rho_k = 0 \in \cF_2/\cF_1(B^{HF}).\]  Under the identification $\left({}_i\rho_k \in \mbox{gr}\left(B^{HF}\right)\right) \leftrightarrow \left({}_i{\bf x}_k \in B^{Kh}\right)$, this agrees with Proposition \ref{prop:MCbetak}, which says:
\[\betatwoleft\left[{}_i{\bf x}_k \otimes ({\bf u}^* \otimes {\bf u})\right] := 0.\]

Lemma \ref{lem:MCbetakHF} also tells us that when $j \leq k-1$:
\[\left(\beta_k^{HF}\right)_{(0|1|1)}\left[({\bf v}^* \otimes {\bf u}) \otimes {}_k(\rho+\sigma)_j\right] := {}_{k-1}(\rho+\sigma)_j.\]  Since $({\bf v}^* \otimes {\bf u}), {}_k(\rho+\sigma)_j,$ and ${}_{k-1}(\rho + \sigma)_j$ are all in $\cF_0/\cF_{-1}$, the induced map on the associated graded is still:
\[\left(\beta_k^{HF}\right)_{(0|1|1)}\left[({\bf v}^* \otimes {\bf u}) \otimes {}_k(\rho+\sigma)_j\right] := {}_{k-1}(\rho+\sigma)_j.\]

Under the identification $\left({}_i\Id_j := {}_i(\rho+\sigma)_j\in \mbox{gr}(B^{HF})\right)  \leftrightarrow {}_i\Id_j \in B^{Kh}$, this agrees with Proposition \ref{prop:MCbetak} which says:
\[\betatworight[({\bf v}^* \otimes {\bf u}) \otimes {}_k\Id_j] := {}_{k-1}\Id_j.\]
\end{proof}

\begin{proof}[Proof of Theorem \ref{thm:filtdiffmodgen}] Now that we have a filtration on the $A_\infty$ bimodule $\cM^{HF}_{\sigma_k^\pm}$ yielding a spectral sequence from $\cM^{Kh}_{\sigma_k^\pm}$ to $\cM^{HF}_{\sigma_k^\pm}$ for each elementary Artin generator, $\sigma_k^\pm$, we would like to construct a filtered $A_\infty$ bimodule $\cM^{HF}_\sigma$ and corresponding spectral sequence $\cM^{Kh}_\sigma \rightarrow \cM^{HF}_\sigma$ for every $\sigma \in B_{m+1}$.

We begin with a decomposition $\sigma = \sigma_{k_1}^\pm \cdots \sigma_{k_n}^\pm$ and define  \[\cM^{HF}_{\sigma} := \cM^{HF}_{\sigma_{k_1}^\pm} \widetilde{\otimes}_{B^{HF}} \ldots \widetilde{\otimes}_{B^{HF}} \cM^{HF}_{\sigma_{k_n}^\pm},\]  which has the structure of a filtered $A_\infty$ bimodule, by Lemma \ref{lem:Ainftensorprod}.

We then check that the associated graded complex of $\cM^{HF}_{\sigma}$ is equivalent to $\cM^{Kh}_\sigma$ in $D_\infty\left(B^{Kh}\right)$, i.e.:
\begin{eqnarray*}
\mbox{gr}\left(\cM^{HF}_{\sigma}\right) &\sim& \cM^{Kh}_\sigma\\
\mbox{gr}\left(\cM^{HF}_{\sigma_{k_1}^\pm} \widetilde{\otimes}_{B^{HF}}\,\, \ldots \widetilde{\otimes}_{B^{HF}}\,\, \cM^{HF}_{\sigma_{k_n}^\pm}\right) &\sim& \cM^{Kh}_{\sigma_{k_1}^\pm} \widetilde{\otimes}_{B^{Kh}}\,\, \ldots \widetilde{\otimes}_{B^{Kh}}\,\ \cM^{Kh}_{\sigma_{k_n}^\pm}
\end{eqnarray*}
in $D_\infty(B^{Kh})$.

Lemma \ref{lem:assgrtensor} tells us that
\[\mbox{gr}\left(\cM^{HF}_{\sigma_{k_1}^\pm} \widetilde{\otimes}_{B^{HF}}\,\, \ldots \widetilde{\otimes}_{B^{HF}}\,\, \cM^{HF}_{\sigma_{k_n}^\pm}\right) \sim \mbox{gr}\left(\cM^{HF}_{\sigma_{k_1}^\pm}\right)\widetilde{\otimes}_{\mbox{gr}(B^{HF})}\,\, \ldots \widetilde{\otimes}_{\mbox{gr}(B^{HF})}\,\, \mbox{gr}\left(\cM^{HF}_{\sigma_{k_n}^\pm}\right)\] as bimodules over $\mbox{gr}(B^{HF})$.  Therefore, they are equivalent in $D_\infty\left(B^{Kh}\right)$, since $\mbox{gr}(B^{HF})$ is isomorphic to $B^{Kh}$ (Theorem \ref{thm:filtdiffalg}).  Furthermore, we also know (Proposition \ref{prop:filtdiffmodelem}) that $\mbox{gr}\left(\cM^{HF}_{\sigma_{k_i}^\pm}\right) \sim \cM^{Kh}_{\sigma_{k_i}^\pm}$ in $D_\infty(B^{Kh})$, so we have
\begin{eqnarray*}
\mbox{gr}\left(\cM_{\sigma}^{HF}\right) &=& \mbox{gr}\left(\cM^{HF}_{\sigma_{k_1}^\pm} \widetilde{\otimes}_{B^{HF}}\,\, \ldots \widetilde{\otimes}_{B^{HF}}\,\, \cM^{HF}_{\sigma_{k_n}^\pm}\right)\\
 &\sim& \cM^{Kh}_{\sigma_{k_1}^\pm} \widetilde{\otimes}_{B^{Kh}} \,\,\ldots \widetilde{\otimes}_{B^{Kh}} \,\,\cM^{Kh}_{\sigma_{k_n}^\pm} = \cM^{Kh}_\sigma,
\end{eqnarray*}
as desired.
\end{proof}

\section{An Example} \label{sec:Example}
We have now constructed, for every braid $\sigma$, a filtration on the Floer bimodule $\mathcal{M}_\sigma^{HF}$ whose associated graded bimodule is quasi-isomorphic to the Khovanov-Seidel bimodule $\mathcal{M}_\sigma^{Kh}$. It is natural to wonder about the $A_\infty$ operations induced by Proposition \ref{prop:InducedAinfty} on the pages of the corresponding spectral sequence.

As a first point, we note that \cite[Prop. 4.9]{MR1862802} (see also \cite[Lem. 4.1]{FaithfulMCG}) implies that the ranks of $H_*(\mathcal{M}_\sigma^{Kh})$ and $H_*(\mathcal{M}_\sigma^{HF})$ always agree. On the other hand, the effect of the filtration on the higher $A_\infty$ operations is in general nontrivial. For example, the $m_2$ products in $B^{Kh}$ (Remark \ref{rmk:BKhMatrix}) and $B^{HF}$ (Remark \ref{rmk:BHFMatrix}) differ.

A more interesting manifestation of the non-triviality of the operations induced by the filtration can be found by examining $HH_*(\mathcal{M}_\sigma^{HF})$, the Hochschild homology of $\mathcal{M}_\sigma^{HF}$.  As in the proofs of Lemma \ref{lem:assgrtensor} and Theorem \ref{thm:filtdiffmodgen}, we have an induced filtration on $HH_*(\mathcal{M}_\sigma^{HF})$, and the associated graded homology can be identified with $HH_*(\mathcal{M}_\sigma^{Kh})$.

Moreover, it is proved in \cite[Thm. 14]{GT10030598} that $HH_*(\mathcal{M}_\sigma^{HF})$ is isomorphic to the next-to-top Alexander grading of 
\[\left\{\begin{array}{cl}
	\widehat{HFK}(\Sigma(\widehat{\sigma}),\widetilde{K}_B) & \mbox{when the braid index of $\sigma$ is even, and}\\
	\widehat{HFK}(\Sigma(\widehat{\sigma}),\widetilde{K}_B) \otimes V & \mbox{when the braid index of $\sigma$ is odd}
	\end{array}\right.\] where:
\begin{itemize}
	\item $\widehat{HFK}(\Sigma(\widehat{\sigma}),\widetilde{K}_B)$ denotes the ``hat" version of the knot Floer homology of the preimage, $\widetilde{K}_B$, of the braid axis, $K_B$, in the double-branched cover of $S^3$ over $\widehat{\sigma}$, the closure of $\sigma$, and
	\item $V= \mathbb{F}_{(0,0)} \oplus \mathbb{F}_{(-1,-1)}$ is a ``standard" $2$--dimensional vector space (the subscripts on the generators indicate their (Alexander, Maslov) bigrading).
\end{itemize}

\begin{remark}
Note that the extra factor of $V$ arises in the odd braid index case because $B^{HF}$ is a strands algebra associate to a {\em twice}--pointed matched circle. The pairing theorem then implies that the Hochschild homology of $\mathcal{M}_\sigma^{HF}$ coincides with the sutured Floer homology of the double--branched cover of $A \times I$ branched over the odd-index braid. As spelled out in \cite[Ex. 2.4]{MR2253454}, a sutured Heegaard diagram for a knot complement with $4$ (rather than $2$) meridional sutures corresponds to a $4$--pointed Heegaard diagram for the knot. The extra pair of basepoints has the effect of tensoring the knot Floer complex with $V$ as in \cite[Thm. 1.1]{CombHFK}.
\end{remark}

In \cite{HochHom}, we prove an analogous result on the Khovanov-Seidel side: namely that the $0$th Hochschild homology of $\mathcal{M}_\sigma^{Kh}$ is isomorphic to the next-to-top filtration grading of the so-called {\em sutured annular Khovanov homology}, $\mbox{SKh}(\widehat{\sigma} \subset A \times I)$, of the closure of $\sigma$ in the solid torus complement of $K_B$, considered as a product sutured manifold $A \times I$. This invariant of the isotopy class of $\widehat{\sigma} \subset (A \times I)$ was first defined in \cite{MR2113902} and studied extensively in \cite{GT07060741} (where it is denoted $H(\widehat{\sigma})$) and in \cite{AnnularLinks} (where it is denoted $Kh^*(\widehat{\sigma})$). By the {\em filtration} grading, we mean the $k$--grading of \cite{GT07060741} and \cite{AnnularLinks}. The following facts are well-known to the experts (see e.g. \cite[Thm. 8.1]{GT07060741}, \cite[Sec. 4]{JacoFest}, \cite[Thm. 3.1]{AnnularLinks}).

\begin{proposition} \label{prop:SKh} Let $\widehat{\sigma} \subset (A \times I)$ be the annular closure of the $n$--strand braid $\sigma$. Letting $\mbox{SKh}(\widehat{\sigma};k)$ denote the sutured annular Khovanov homology of $\widehat{\sigma}$ in filtration grading $k$, we have:
\begin{enumerate}
	\item $\mbox{SKh}(\widehat{\sigma};k) = 0$ unless $k \in \{-n, -(n-2), \ldots, n-2, n\}$,
	\item $\mbox{SKh}(\widehat{\sigma};n) = \F$.
	\item $\mbox{SKh}(\widehat{\sigma};k) \cong \mbox{SKh}(\widehat{\sigma};-k)$ for all $k \in \Z$,
	\item There is a spectral sequence whose $E^1$ page is $\mbox{SKh}(\widehat{\sigma})$ and whose $E^\infty$ page is $\mbox{Kh}(\widehat{\sigma})$, the ordinary $\Z/2\Z$ Khovanov homology of $\widehat{\sigma}$.
\end{enumerate}
\end{proposition}

\begin{proof}

Statements (1) and (2) are immediate consequences of the correspondence between generators of the chain complex underlying $\mbox{SKh}(\widehat{\sigma})$ and enhanced Kauffman states (see, e.g., \cite[Sec. 4.2]{JacoFest}).

To understand (3), we once again use the identification between enhanced Kauffman states and generators of $\mbox{CKh}(\widehat{\sigma})$, the chain complex underlying $\mbox{SKh}(\widehat{\sigma})$.  One then constructs inverse chain maps \[\mbox{CKh}(\widehat{\sigma};k) \leftrightarrow \mbox{CKh}(\widehat{\sigma};-k)\] by reversing the orientations of those circular components of the enhanced Kauffman state representing nontrivial elements of $H_1(A)$.

Statement (4) is \cite[Lem. 1]{GT07060741}.
\end{proof}

Consider now the $3$--braid $\sigma = (\sigma_1\sigma_2)^5$ whose closure is the positive $(3,5)$ torus knot $T_{3,5}$. The double-branched cover, $\Sigma(\widehat{\sigma})$, is the Heegaard Floer $L$--space integer homology sphere $\Sigma(2,3,5)$, and the preimage of the braid axis $\widetilde{K}_B \subset \Sigma(2,3,5)$ is a genus one fibered knot whose corresponding open book has positive monodromy, hence is compatible with a Stein fillable contact structure \cite{LoiPierg}. It follows that its Heegaard-Floer contact invariant  \cite{MR2153455} is nonzero, so \cite[Prop. 3.1]{BaldwinGOne} implies that $\widehat{HFK}(\Sigma(2,3,5), \widetilde{K}_B)$ has rank one in the next-to-top Alexander grading, hence (recalling also \cite[Thm. 1.1]{MR2153455}) $HH_*(\mathcal{M}_\sigma^{HF})$ has rank $2$.

On the other hand, Proposition \ref{prop:SKh}(4), tells us that $\mbox{rk}(\mbox{SKh}(\widehat{\sigma}))$ is bounded below by $\mbox{rk}(\mbox{Kh}(\widehat{\sigma}))$, and \cite{KnotAtlas} tells us that the rank of the $\Z/2\Z$ Khovanov homology is $14$. The main result of \cite{HochHom} combined with Proposition \ref{prop:SKh} now implies: \begin{eqnarray*}
	\mbox{rk}(HH_*(\mathcal{M}_\sigma^{Kh})) &\geq& \mbox{rk}(HH_0(\mathcal{M}_\sigma^{Kh}))\\
	&=& \mbox{rk}(\mbox{SKh}(\widehat{\sigma};1))\\
	&=& \frac{1}{2}\left(\mbox{rk}(\mbox{SKh}(\widehat{\sigma})) - 2\right)\\
	&\geq& \frac{1}{2}\left(\mbox{rk}(\mbox{Kh}(\widehat{\sigma})) - 2\right)\\
	&=& 6,
\end{eqnarray*} so the $A_\infty$ structures on $\mathcal{M}_\sigma^{Kh}$ and $\mathcal{M}_\sigma^{HF}$ must differ.
\bibliography{QuiverAlgebras}
\end{document}